\documentclass[11pt]{amsart}
\usepackage[margin=1in]{geometry}
\numberwithin{equation}{section}

\usepackage{amssymb,amsmath,amsthm}
\usepackage{bbm}
\usepackage{xspace}
\usepackage{marvosym} 
\usepackage[sans]{dsfont} 
\usepackage{etoolbox}
\usepackage{csquotes} 
\usepackage[dvipsnames]{xcolor}
\usepackage[bottom,first]{draftcopy}
\usepackage{changebar}

\usepackage[textsize=tiny, backgroundcolor=blue!20!white, linecolor=blue!20!white,]{todonotes}
\usepackage[colorlinks]{hyperref}
\usepackage{graphicx}
\usepackage{mathtools}
\usepackage{cleveref}
\usepackage{enumitem}
\usepackage{cancel}
\usepackage{color}
\usepackage[normalem]{ulem}
\usepackage{comment}

\allowdisplaybreaks[4]
\newtheorem{thm}{Theorem}[section]
\newtheorem{lem}[thm]{Lemma}
\newtheorem{cor}[thm]{Corollary}
\newtheorem{sublem}[thm]{Sub-lemma}
\newtheorem{prop}[thm]{Proposition}

\newtheorem{exmp}[thm]{Example}
\theoremstyle{remark}
\newtheorem{rem}[thm]{Remark}
\newtheorem*{ack}{Acknowledgement}

\newcommand\eps{\eps}

\newcommand{\fG}{\mathfrak{G}}

\newcommand\One{{\bf 1}}
\newcommand{\Id}{\mathsf{Id}}
\newcommand{\BV}{{\operatorname{BV}}}
\newcommand\V{{\mathsf V_{\alpha,\beta,\gamma}}}
\newcommand{\VVc}[3]{{\mathsf V}_{#1,#2,#3}}

\DeclareMathOperator*{\esssup}{ess\,sup}

\DeclareMathOperator*{\osc}{osc}

\newcommand{\wh}[1]{\widehat{#1}}
\newcommand{\vs}{\vspace{5pt}}

\makeatletter
\def\@tocline#1#2#3#4#5#6#7{\relax
  \ifnum #1>\c@tocdepth 
  \else
    \par \addpenalty\@secpenalty\addvspace{#2}%
    \begingroup \hyphenpenalty\@M
    \@ifempty{#4}{%
      \@tempdima\csname r@tocindent\number#1\endcsname\relax
    }{%
      \@tempdima#4\relax
    }%
    \parindent\z@ \leftskip#3\relax \advance\leftskip\@tempdima\relax
    \rightskip\@pnumwidth plus4em \parfillskip-\@pnumwidth
    #5\leavevmode\hskip-\@tempdima
      \ifcase #1
       \or\or \hskip 1em \or \hskip 2em \else \hskip 3em \fi%
      #6\nobreak\relax
    \hfill\hbox to\@pnumwidth{\@tocpagenum{#7}}\par
    \nobreak
    \endgroup
  \fi}
\makeatother

\def\eps{{\varepsilon}}

\def\Prob{{\mathbb{P}}}

\def\EXP{{\mathbb{E}}}

\def\complex{\mathbb{C}}

\def\naturals{\mathbb{N}}

\def\reals{\mathbb{R}}

\def\cB{\mathcal{B}}

\def\cC{\mathcal{C}}

\def\cD{\mathcal{D}}

\def\cL{\mathcal{L}}

\def\cM{\mathcal{M}}

\def\cO{\mathcal{O}}

\def\cX{\mathcal{X}}

\def\fF{\mathfrak{F}}

\def\fN{\mathfrak{N}}

\def\fn{\mathfrak{n}}

\makeatother

\def\beq{\begin{equation}}
\def\eeq{\end{equation}}

\begin{document}

\title[Limit Theorems for unbounded observables]{Limit Theorems for a class of unbounded observables with an application to ``Sampling the Lindel\"of hypothesis"}
\author{Kasun Fernando and Tanja I.~Schindler}

\address{Kasun Fernando\\
Department of Mathematics\\
Brunel University London}
\email{{\tt kasun.fernando@brunel.ac.uk}}

\address{Tanja I.\ Schindler\\
Faculty of Mathematics and Computer Science\\ Jagiellonian University in Krakow, Poland and\\
Department of Mathematics and Statistics, University of Exeter, UK}
\email{{\tt tanja.schindler@uj.edu.pl, t.schindler@exeter.ac.at}}

\begin{abstract}
We prove the Central Limit Theorem (CLT), the first order Edgeworth Expansion and a Mixing Local Central Limit Theorem (MLCLT) for Birkhoff sums of a class of unbounded heavily oscillating observables over a family of full-branch piecewise $C^2$ expanding maps of the interval. As a corollary, we obtain the corresponding results for Boolean-type transformations on $\reals$. The class of observables in the CLT and the MLCLT on $\reals$ include the real part, the imaginary part and the absolute value of the Riemann zeta function. Thus obtained CLT and MLCLT for the Riemann zeta function are in the spirit of the results of Lifschitz \& Weber \cite{LW} and Steuding \cite{steuding} who have proven the Strong Law of Large Numbers for \textit{sampling the Lindel\"of hypothesis}.\vspace{-10pt} 
\end{abstract}

\subjclass[2020]{37A50, 60F05, 37A44, 11M06}
\keywords{central limit theorem, mixing local limit theorem, Edgeworth expansion, ergodic limit theorems, unbounded observables, expanding interval maps, Riemann zeta functions, Lindel\"of hypothesis, quasicompact transfer operators, Keller-Liverani perturbation theory}

\maketitle

\tableofcontents

\section{Introduction}

Expanding maps of the unit interval are the most elementary class of dynamical systems that exhibit chaotic behaviour. There is  extensive body of literature on limit theorems for Birkhoff sums of expanding maps as summarised in \cite{Den, hennion_limit_2001, FL}. In particular, in \cite{RE}, the central limit theorem (CLT) is established for observables with bounded variation ($\BV$) over piecewise uniformly expanding maps whose inverse derivative is also $\BV$. Further, in \cite{FL}, Edgeworth expansions describing the error terms in the CLT are established in the case of $\BV$ observables over $C^2$ uniformly expanding maps that are \textit{covering}. In both cases, the results are limited to bounded observables since the observables considered are $\BV$.

One standard technique of establishing limit theorems for dynamical systems is the Nagaev-Guivarc'h spectral approach which was first introduced by Nagaev in the Markovian setting in \cite{NG1}, and later, adapted to deterministic dynamical systems by Guivarc'h in \cite{GH}. The key idea is to code the characteristic function using iterated twisted transfer operator (one can think of this as the deterministic counterpart of the dual of the Markov operator in the Markovian setting) and to analyse the spectral data of these families of operators in a suitable Banach space. More precise formulations of this idea can be found in \cite{Gora, hennion_limit_2001, G}. 

Though transfer operator techniques to handle unbounded observables are available, see for example, \cite{HervePene, Butterley, Liverani1, FP}, there are only a few results for limit theorems for unbounded observables over uniformly expanding maps of the interval; see for example,\ \cite{arb_sma} and \cite{chen_zhang} of which the latter, however, does not  use transfer operator techniques.
In \cite{Butterley, Liverani1, arb_sma}, the goal was to obtain bounds for the spectral radius and the essential spectral radius of the (twisted) transfer operators associated with expanding maps acting on their corresponding Banach spaces. While the first two works did not address limit theorems, in \cite{arb_sma}, a CLT and an almost sure invariance principle was proven. However, to the best of our knowledge, nothing is known about the first order Edgeworth expansion (the quantitative CLT) or the mixing local central limit theorem (MLCLT) in this setting.

Notably, \cite{HervePene} introduced a general framework for establishing the first order Edgeworth expansion in a Markovian context that is nearly optimal, comparable to conditions in the independent identically distributed (IID) setting. This framework was further extended in \cite{FP} to obtain expansions of all orders in both the CLT and the MLCLT, with applications to systems modelled by shifts of finite type or Young towers and unbounded observables with nearly optimal order of integrability. Despite its potential, this generalised theory has not been applied to expanding maps until this work. 

In this paper, we introduce a class of Banach spaces that are not contained in $L^\infty$ and for which the conditions introduced in \cite{HervePene, FP} can be verified in the context of $C^2$ uniformly expanding maps of the interval. In Remark \ref{rem: Banach spaces}, we compare the class of Banach spaces we introduce with other Banach spaces in the literature that include unbounded observables and are known to possess a spectral gap for the associated transfer operator.

The observables $\chi:(0,1)\to \reals$ we focus on and which belong to our Banach space are \textit{unbounded heavily oscillating observables} characterised by the conditions $$
 |\chi|\lesssim x^{-a}(1-x)^{-a} \quad \text{and} \quad
  \max\{|\chi'(x+)|, |\chi'(x-)|\}\lesssim x^{-b}(1-x)^{-b}\,
$$ for some $a,b > 0$. The permissible ranges for $a$ and $b$ vary depending on the specific interval map and the limit theorem of interest.   

As the underlying transformations we consider are full-branched $C^2$ uniformly expanding maps of the interval and we see that the non removable singularities of $\chi$ are always at a fixed point of the interval map. (This, however, could easily be extended to any periodic point.) The behaviour of such maps can be considered as particularly interesting:~once an orbit lands close to a fixed point, a few subsequent iterates might stay relatively close to the fixed point and the Birkhoff sum might be very large locally. Such situations can cause the system to behave qualitatively different from the  
IID setting, see for example, \cite[Theorem 1.10]{kessschimean}.  

Further, we show that the general framework developed in \cite{FP} for limit theorems involving unbounded observables can be applied to our class of observables. By adapting the ideas in \cite{FP} to our context, in \Cref{sec:AbsRes}, we identify a set of sufficient conditions on both the system and the observables which ensure the validity of limit theorems. In particular, we establish the CLT, its first order correction -- the first order Edgeworth expansion, and a 
MLCLT for the Birkhoff sums of $\chi$.

Indeed, we consider a sequence of increasing Banach spaces (all of them containing unbounded observables) on each of which the twisted transfer operators corresponding to full-branch $C^2$ expanding maps satisfy Doeblin-Fortet Lasota-Yorke (DFLY) inequalities and other good spectral properties. These properties, in turn, lead to the establishment of the limit theorems for this class of expanding maps. In the previous literature (including \cite{FP}), the conditions in the general framework were not verified in a context similar to ours, and there lies the key novelty of this work. Though our results regarding the introduced Banach spaces are tailored to prove limit theorems in the spirit of \cite{FP}, similar, general techniques using a chain of Banach spaces were established in \cite{HervePene} and also used in \cite{Pene}. 
Having this in mind, the intermediate technical results of this paper regarding the precise details of the chain of Banach spaces (relegated to \Cref{subsec: Banach space corr}) might be helpful to prove further limit theorems in the future.

As an immediate application, we deduce limit theorems for a Boolean-type transformation $\phi: \mathbb{R} \to \mathbb{R}$ given by $\phi(x)=1/2(x-1/x)$ if $x\neq0$ and observables that heavily-oscillate at $\pm \infty$. This application makes use of a smooth conjugacy between the doubling map on the unit interval and $\phi$.
In particular, this has applications to \textit{sampling the Lindel\"of hypothesis}, a line of research in analytic number theory that deals with understanding the properties of the Riemann zeta function on the critical strip. We elaborate on this in \Cref{sec:Lindelof}.
\newline

The structure of the paper is as follows:~\Cref{sec:Results} is dedicated to preliminaries and main results:~in \Cref{sec:Notation}, we introduce the relevant notation and common definitions that we will use throughout the paper, in \Cref{sec:DS}, we state precisely the class of expanding maps we consider, in \Cref{sec:IntroSpaces}, we introduce our Banach spaces, in \Cref{sec:IntThms}, we state our main results for the interval maps, and in \Cref{subsec: Boolean trafo}, we state the corresponding results for the Boolean transformation on $\reals$ and their implications to sampling the Lindel\"of hypothesis is discussed in \Cref{sec:Lindelof}. In \Cref{sec:AbsRes}, we recall known abstract results in \cite{HervePene, FP} tailored (with justifications) to our setting.  
The spectral properties of twisted transfer operators acting on these spaces including the DFLY inequality are established in \Cref{sec:Operators}.  In \Cref{sec:Proofs}, we collect the proofs of our main results. In particular, the proofs of the limit theorems for interval maps appear in \Cref{sec:LimitThms} and in \Cref{sec:BLimitThms} we prove the corresponding results for the Boolean-type transformation. Finally, we have relegated some technical results to the Appendices; in particular, an in-depth discussion about our Banach spaces appears in \Cref{subsec: Banach space corr}.

\section{Main Results}\label{sec:Results}
\subsection{Preliminaries} \label{sec:Notation}
Let $X$ be a metric space with a reference Borel probability measure $m$, and let $T :X\to X$ be a non-singular dynamical system, i.e., for all $U \subseteq X$ Borel subsets $m(U)=0$ holds if and only if $m({ T}^{-1} U)=0$ holds. We denote by $\cM_1(X)$ the set of Borel probability measures on $X$. Let $\nu \in \cM_1(X)$.
For $p\geq 1$, by $L^p(\nu)$, we denote the standard Lebesgue spaces with respect to $\nu$, i.e., $$L^p(\nu) = \{ h: X \to X\,|\, h\, \text{is Borel measurable},\, \nu(|h|^p)<\infty \}$$ where the notation $\nu(h)$ refers to the integral of a function $h$ with respect to a measure $\nu$ and the corresponding norm is denoted by $\|\cdot\|_{L^p(\nu)}$. When $\nu=m$, we often write, $L^p$ instead of $L^p(m)$ and  $\|\cdot\|_{p}$ instead of $\|\cdot\|_{L^p(m)}$. 

For us, an observable is a real valued function $f \in L^{2}$ for which we consider the Birkhoff sums (also commonly referred to as ergodic sums), 
\begin{equation}\label{eq:BirkhoffSum}
    S_{n}({ f,T}) = \sum_{k=0}^{n-1} { f \circ 
    T}^{k}
\end{equation}
which we denote by $S_{n}(f)$ when the dynamical system $T$ is fixed. 

We say $\wh{ T}:L^1\to L^1$ is the transfer operator of ${ \wh T}$ with respect to $m$, if for all ${ f}\in L^1$ and ${ f}^* \in L^\infty$,
\begin{equation}\label{eq:Duality}
    m(\wh{ T}({ f})\cdot { f}^*) = m({ f}\cdot { f}^* \circ { T} ).
\end{equation}
Let ${\overline{m}} \in \cM_1(X)$ be absolutely continuous with respect to $m$ with density $\rho_{\overline{m}}$.
Then, from \eqref{eq:Duality}, it follows that
\begin{equation}\label{eq:CharFunc}
    \EXP_{\overline{m}}(e^{isS_{n}({ f})})=m\left( {{\wh{ T}}}^n_{ is}(\rho_{\overline{m}})\right)
\end{equation}
where $\EXP_{{\overline{m}}}$ is the expectation with respect to the law of $S_n$ where the initial point $x$ is distributed according to $ {\overline{m}}$ and 
\begin{equation}\label{eq:cL_theta_is}
    {\wh{ T}}_{is}(\cdot)=\wh{ T}(e^{is { f}}\, \cdot ),\ s \in \reals \,,
    \end{equation}
see, for example, \cite[Chapter 4]{hennion_limit_2001}. Eventually, we are interested in the asymptotics of quantities of the form ${\overline{m}}(S_n(f) \leq z_n)$ and $\EXP_{{\overline{m}}}(V_n(S_n(f)))$ as $n \to \infty$ where $z_n \in \reals$ and $V_n : \reals \to \reals$ are from a suitable class of observables, and to obtain these asymptotics we exploit the relation \eqref{eq:CharFunc}.

We denote
\begin{align*}
  A  = \lim_{n \to \infty} \EXP_{\overline{m}}\left(\frac{S_{n}(f, T)}{n}\right)\quad\,\,\,\text{and}\quad\,\,\, 
     \sigma^2
     =\lim_{n\to\infty} \EXP_{\overline{m}}\left(\frac{S_{n}(f,T)-n\,A}{\sqrt{n}}\right)^2 
\end{align*}
for the asymptotic mean and the asymptotic variance of Birkhoff sums, $S_{n}(f, T)$, respectively. Under the assumptions we impose on $T$ in \Cref{sec:DS}, $A$ and $\sigma^2$ are independent of the choice of ${\overline{m}}$ because each initial measure $\overline m$ converges weakly to the unique unique absolutely continuous invariant probability measure (acip) under the action of $\wh{T}$,  and we may focus on zero average observables by considering $\overline{f} := f - A$ instead of $f$.

We call $f$ to be $T-$\textit{cohomologous to a constant} in the function space $\fF$ if there exist $\ell \in \fF$ 
and a constant $c$ such that 
$$f= \ell \circ T - \ell + c\,$$ 
and $T-$\textit{coboundary} in $\fF$ if there exists $\ell\in \fF 
$ such that 
 $$f= \ell \circ T - \ell\,.$$
We say $f$ is \textit{non-arithmetic} in $\fF$ if it is not $T-$cohomologous in $\fF$ to a sublattice-valued function, i.e., if there exists no triple $(\gamma , B,\theta)$ with $\gamma : X \to \reals$, $B$ a closed proper subgroup of $\reals$ and a constant $\theta $ such that $ f + \gamma - \gamma \circ T \in \theta + B$.

Given a Banach space $\cB_1$, the $\complex-$valued continuous linear functionals are denoted by $\cB_1^{\,\prime}$ and given another Banach space $\cB_2$, $\cL(\cB_1, \cB_2)$ denotes the space of bounded linear operators from $\cB_1$ to $\cB_2$. 
When $\cB_1=\cB_2$, we write $\cL(\cB_1 ,\cB_1)$ as $\cL(\cB_1)$. When $\cB_1 \subset \cB_2$, $\cB_1 \hookrightarrow \cB_2$ denotes continuous embedding of Banach spaces, i.e., there exists $\mathfrak c >0$ such that $\|\cdot\|_{\cB_2} \leq \mathfrak c \|\cdot\|_{\cB_1}$.  

Given a set $D \subseteq X$,  its complement $X \setminus D$ is denoted by $D^c$, and $\mathring{D}$ denotes its interior. Given a function $f: D \to \reals$ set $f_+:=\max\{f,0\}$ and $f_-:=\max\{-f,0\}$. Given $g : D \to \reals$, $g \lesssim f$ denotes that there exists constant $K>0$ such that $g(x) \leq K f(x),$ for all $x \in D$. Let $Q_1,Q_2$ be $\reals^{+}_0$ valued functionals acting on a class of functions $\fG_1$ and $\fG_2$, the inequality $Q_1(g) \lesssim Q_2(h)$ for all $g \in \fG_1$ and $h \in \fG_2$ is written to denote that there exists $K$ independent of the choices of $g$ and $h$ such that $Q_1(g)\leq K Q_2(h)$. 
Finally, given two numbers $a,b \in \reals$, $a \approx b$ means that $0 \leq a-b \leq 10^{-3}$. 

We denote the standard Gaussian density and the corresponding distribution function by
\begin{equation*}\label{eq:standardGauss}
\fn(x)=\frac{1}{\sqrt{2\pi}}e^{-x^2/2}\,\,\quad \text{and}\quad\,\, \fN(x)=\int_{-\infty}^x \fn(y)\, \mathrm{d}y\,,
\end{equation*}
respectively.

\subsection{The classes of dynamical systems}\label{sec:DS}

Let $I=[0,1]$ and $\lambda$ the Lebesgue measure (on $\reals$) and $\lambda_I$ its restriction to $I$. We use  $\lambda_I$ as the reference measure on $I$ and let $I=\bigcup_{j=0}^{k-1} [c_j,c_{j+1}]$ be a partition of $I$ with $c_0=0$ and $c_{k}=1$. We consider the class of maps $\psi: I \to I$ satisfying the following conditions. 
\begin{enumerate}[leftmargin=*]
    \item There are  $\psi_{j+1}:[c_{j},c_{j+1}] \to I $ such that for all $j$, $\psi_{j+1} \in C^2\,,$ 
    $|\psi'_{j+1}|>1\,,$ 
     $\mathrm{Range}(\psi_{j+1})=I$ and $$\psi_{j+1}|_{(c_{j},c_{j+1})} =\psi |_{(c_{j},c_{j+1})}\,.$$ 
    \item For all $j$, 
    the derivative of $\psi^{-1}_{j+1}$ is uniformly ${\vartheta}-$H\"older, i.e., there exists $c$ such that for all $j$, for all $\eps>0$, for all $z\in I$ and for all $x,y \in B_\eps(z):=[z-\eps,z+\eps]\cap [0,1]\,,$
    $$|(\psi^{-1}_{j+1})'(x)-(\psi^{-1}_{j+1})'(y)|\leq c|(\psi^{-1}_{j+1})'(z)|\eps^{\vartheta}
    \,.$$
\end{enumerate}
\begin{rem}
The full branch assumption was made in order to simplify our calculations. Later on, we will be particularly interested in the doubling map which is conjugate to the Boolean type transformation, a transformation over the real line.

Since the maps studied here are $C^2$, Markov and topologically mixing, each map has one and only one acip and it is exact \cite[Theorem 6.1.1]{Gora}. We denote this acip by $\pi$. Since $\psi'_{j+1}$ are $C^1\,,$ there exists $\eta_{+}<\infty$ such that
$$\max_{j} \|\psi'_{j+1}\|_{\infty} = \eta_{+}.$$
Also, since $|\psi'_{j+1}|>1\,,$ there exists $\eta_->1$ such that $$\max_{j} \|(\psi^{-1}_{j+1})'\|_{\infty} = 1/\eta_-.$$
\end{rem}

\begin{rem}
If $f \in \mathfrak{F}$ is a $\psi$-coboundary in the space of measurable functions, then it is a $\psi$-coboundary in the class of piecewise $C^2$ functions, see \cite{Liv1, Jen}.
\end{rem}

Without loss of generality we assume that $\psi'>0$ and we have 
\begin{equation}\label{eq:TransferOp}
    \wh\psi_{is}(\varphi)(x)=\sum_{j=0}^{k-1} \frac{e^{is\chi(\psi^{-1}_{j+1}x)}}{\psi'(\psi^{-1}_{j+1}x)}\varphi(\psi^{-1}_{j+1}x)\,,
\end{equation}
see, for example, \cite{hennion_limit_2001} for a proof of this fact.
\subsection{The Banach spaces}\label{sec:IntroSpaces} For a measurable function $f\colon I\to \reals$ and a Borel subset $S$ of $I$, we define the oscillation on $S$ by
\begin{align*}
 \osc\left(f,S\right)\coloneqq \esssup_{x,y\in S} |f(x)- f(y)| 
\end{align*}
 and we set $\osc(f,\emptyset)\coloneqq 0$.
For a complex valued functions $f$ we generalize the definition to 

\begin{align*}
 \osc\left(f,S\right)\coloneqq \osc(\mathfrak{R}f,S) +\osc(\mathfrak{I}f,S), 
\end{align*} 
where $\mathfrak{R}f$ and $\mathfrak{I}f$ refer to real and imaginary parts of $f$, respectively. Also, in the case of a complex valued function $f$, up to a constant, this is equivalent to the more intuitive definition
\begin{align*}
 \overline\osc\left(f,S\right)\coloneqq \esssup_{x,y\in S} |f(x)- f(y)| \,.
\end{align*}
This can be easily seen:~We have $|f(x)-f(y)|\leq |\Re f(x)-\Re f(y)|+|\Im f(x)- \Im f(y)|$, and thus, $\overline{\osc}(f,S)\leq \osc(f,S)$.  On the other hand, we have $\osc(f,S)\leq 2\max\{\overline{\osc}(\Re f, S), \overline{\osc}(\Im f, S)\}\leq 2\,\overline{\osc}(f,S)$. In what follows, we use $\osc$ as the standard definition.

For $\alpha \in \reals$, define, $R_{\alpha}$, an operator on the space of measurable functions by
$$R_{\alpha}f(x)\coloneqq\left\{ \begin{array}{ll}
	x^{\alpha}\cdot\left(1-x\right)^{\alpha}\cdot f(x) & \mbox{if } |f(x)|< \infty \\
    0 & \mbox{otherwise,} 
	\end{array} 
	\right.\,$$
denote by $B_{\eps}(x)$ the $\eps$-ball around $x$ in $I$, and define a seminorm 
\begin{align*}
 \left|f\right|_{\alpha,\beta}
 \coloneqq \sup_{\eps\in(0,\eps_0]} \eps^{-\beta}\int \osc\left(R_{\alpha}f, B_{\eps}(x)\right)\,\mathrm{d}\lambda_I(x)\,, 
\end{align*}
where $\eps_0$ is sufficiently small (to be chosen later).
 Let
 \begin{align*}
  \left\|\cdot\right\|_{\alpha,\beta,\gamma}\coloneqq \left\|\cdot\right\|_\gamma+\left|\cdot\right|_{\alpha,\beta}
 \end{align*}
  and set
  \begin{align*}
  L^\gamma \coloneqq \left\{f\colon I\to\mathbb{C} \colon \left\|f\right\|_{\gamma}<\infty\right\}\,, \qquad \mathsf{V}_{\alpha,\beta, \gamma}\coloneqq \left\{f\colon I\to\mathbb{C} \colon \left\|f\right\|_{\alpha,\beta,\gamma}<\infty\right\}. 
 \end{align*}
Finally, by $\mathsf{V}^{\,\prime}_{\alpha,\beta, \gamma}$ we denote 
the set of $\complex-$valued continuous linear functionals on $\mathsf{V}_{\alpha,\beta, \gamma}$. 

\begin{rem}
It is shown in \Cref{subsec: Banach space corr} that for $\alpha\in [0,1)$, $\beta \in (0,1]$ and $\gamma \geq 1$, $\V$ is a Banach space. Similar real Banach spaces were considered in \cite{keller_generalized_1985, blank_discreteness_1997, kessebohmer_intermediately_2019}. In all these cases, their spaces correspond to our spaces with $\alpha = 0$, and hence, are embedded in $L^\infty$; see \Cref{lem:CptIncl1}. 

Due to the dampening operation $R_\alpha$, which was first introduced by the second author in \cite{STI}, the observables in $\V$ may be unbounded and oscillate heavily near $0$ and $1$. It was used to establish a CLT. However, there was a critical mistake in the proof:~the normed vector space considered there in order to study the spectrum of the transfer operator is not complete.  In what follows, we not only correct this mistake but also establish a MLCLT for \eqref{eq:BSum}.

Moreover, we remark that depending on the application one could consider different damping operators and use the ideas presented here to prove other limit theorems.

It is clear that due to the structure of $R_{\alpha}$ the oscillating singularity can occur precisely at the fixed points $0$ and $1$. In case that the map $\psi$ has more than two branches, and hence, fixed points $0=a_0< a_1<\ldots< a_{k-2}< a_{k-1}=1$, the proofs could easily be generalised to a Banach space with additional dampening at $a_1,\ldots, a_{k-2}$. In this case we would consider $R_{\alpha}f(x)=\prod_{j=0}^{k-1}|x-a_j|^{\alpha}f(x)$. However, these calculations would make the proof even more technical, and hence, we decided to restrict ourselves to the observables only oscillating at $0$ and at $1$. 

Allowing for singularities of $f$ at a point $y$ which is not a fixed point of $\psi$, i.e.\ introducing a damping function $|x-y|^{\alpha}$ would still result in a Banach space. However, under 
$\psi$, the singularity of $f$ would move and $\VVc{\alpha}{\beta}{\gamma}$ would no longer be closed under the action of the transfer operator $\widehat{\psi}$, a condition fundamental for proofs using transfer operator techniques. 
In the case of $y$ being a periodic point of period $d$, considering $\psi^d$ instead of $\psi$ should work.
\end{rem}

\begin{rem}\label{rem: Banach spaces}
In the literature there are a number of Banach spaces which also allow for unbounded observables, e.g.\ \cite{Butterley, Liverani1, arb_sma}. Those are seemingly more general than the Banach space we introduce because they do not  have the restriction that the singularity can only occur on a fixed point. 

However, the norm of the Banach spaces are defined in an implicit way, e.g., as $$\|h\|=\sup_{ \|\varphi\|_{\alpha} \leq 1}\left|\int_0^1 \varphi'(x) h(x)\,\mathrm{d}x\right|$$ in \cite{Liverani1} with $\varphi$ out of a certain space $(\mathcal{B}_{\alpha},\|\cdot\|_{\alpha})\,,$ 
 and similarly, in an implicit way in \cite{Butterley, arb_sma}.
It is not clear, and not an easy task to check whether the observables we are interested in or even more elementary observables like $x^{-c}\sin(1/x)$, $c>0$ belong to Banach spaces in the literature \cite{Liverani2}. From this point of view, the proposed Banach spaces are interesting because the conditions are relatively easy to check. Moreover, for the method used in our paper, a sequence of Banach spaces is necessary. This would introduce additional technical difficulties, if we were to use other Banach spaces in the literature. 
\end{rem}

\subsection{Results for the unit interval}\label{sec:IntThms}
Now, we are ready to state the limit theorems for $S_n(\chi):=S_n(\chi, \psi)$ over dynamical systems $\psi$ defined as in \Cref{sec:DS}.  Though we do not state this explicitly, it will later turn out that the $\chi$ specified in the following theorems belongs to an appropriate $\V$.

We first state the CLT in the stationary case.

\begin{thm}\label{thm:IntervalCLTsimple}
Suppose $\chi$ is continuous and the right and left derivatives of $\chi$ exist on $\mathring{I}$, $\chi$ is not a coboundary and there exist constants $a,b>0$ such that
 \begin{align}
 |\chi(x)|\lesssim x^{-a}(1-x)^{-a}\quad\text{ and }\quad
  \max\{|\chi'(x+)|, |\chi'(x-)|\}\lesssim x^{-b}(1-x)^{-b}\,. \label{eq: bound ab}
 \end{align}
Assume 
\begin{align}
 a<\min\left\{\vartheta, \frac{1}{b},\frac{1}{2}\right\}\cdot\min\left\{1,\frac{\log\eta_-}{\log \eta_+}\right\}.\label{eq: CLT main cond}
\end{align}
Then, the following Central Limit Theorem holds:
  \begin{equation}\label{eq:IntervalCLT}
   \pi\left(\frac{S_{n}(\chi) - n\,\pi(\chi)}{\sigma\sqrt{n}} \leq x\right) - \fN(x)= o(1),\quad\text{for all}\quad x \in \reals\quad\text{as}\quad n\to \infty\,.
    \end{equation}
\end{thm}

Now, we discuss sufficient conditions for the MLCLT.
\begin{thm}\label{thm:IntervalLLT}
Suppose $\chi$ is continuous and the right and left derivatives of $\chi$ exist on $\mathring{I}$, $\chi$ is not arithmetic and there exist constants $a,b>0$ such that \eqref{eq: bound ab} and \eqref{eq: CLT main cond} are true. Then, $S_n(\chi)$ satisfies the 
following MLCLT:for all $0<\alpha_0<\alpha_1<\beta$, $M\geq 1$,\,  $U\in \VVc{\alpha_0}{\beta}{M}\,,$ $V: \reals \to \reals$  a compactly supported continuous function, ${\overline{m}}\in\cM_1(I)$ being absolutely continuous wrt\ $\lambda_I$, and $W \in L^1$ such that $(W \cdot {\overline{m}} ) \in \VVc{\alpha_1}{\beta}{M}^{\,\prime}\,,$ we have 
    \begin{equation}\label{eq:IntervalLLT}
   \lim_{n \to \infty} \sup_{\ell \in \reals }\left|\sigma\sqrt{2\pi n}\,\EXP_{{\overline{m}}}(U \circ \psi^n \,V(S_n(\overline\chi)-\ell)\, W)  - e^{-\frac{\ell^2}{2n\sigma^2}} \,\EXP_{{\overline{m}}}(W)\,\EXP_\pi(U)\int V(x) \, \mathrm{d}x \right|= 0\,.
   \end{equation}
\end{thm}

\begin{rem}
  In particular, it is possible to choose ${\overline{m}}=\pi$ for all $W\in L^{\bar M}$ where $M^{-1}+\bar M^{-1} = 1$. In fact, under our assumptions, there exists $\rho \in \BV$ such that $\pi = \rho \lambda_I$; see, for example,  \cite{L1}. Therefore, $|W\cdot \pi (h)| = \left|\int (h W)\rho\, \mathrm{d}\lambda_I \right|\leq \|\rho\|_{\infty} |Wh|_{L^1} \leq  \|\rho\|_{\infty} \|W\|_{\bar M} \|h\|_{M} \leq C \|h\|_{\alpha_1, \beta, M}$ with $C=\|\rho\|_{\infty} \|W\|_{\bar M}$, and hence, $W\cdot \pi \in \VVc{\alpha_1}{\beta}{M}^{\,\prime}$ as required.
\end{rem}

Next, we discuss the first order asymptotics of the CLT with no assumptions on the stationarity. In particular, under the conditions of the theorem, we have the CLT for initial measures that are not necessarily invariant. 

\begin{thm}\label{thm:IntervalEdgeExp}
Suppose $\chi$ is continuous and the right and left derivatives of $\chi$ exist on $\mathring{I}$, $\chi$ is arithmetic and there exist constants $a,b>0$ such that \eqref{eq: bound ab} and 
\begin{align}\label{eq: Edge main cond}
  3\min\{2a, \max\{a,a+b-2\}\}<\min\left\{\vartheta, \frac{1}{b},\frac{1}{2}\right\}\cdot\min\left\{1,\frac{\log\eta_+}{\log \eta_-}\right\}\,.
\end{align} 
are true. Then, $S_n(\chi)$ satisfies the first order Edgeworth expansion, i.e.,  for all ${\overline{m}}\in \cM_1(I)$ being absolutely continuous wrt $\lambda_I$ there exists a quadratic polynomial $P$ whose coefficients depend on the first three asymptotic moments of $S_n(\chi)$ but not on $n$ 
 such that 
  \begin{equation*}
\sup_{ x\in \reals} \left|{\overline{m}}\left(\frac{S_n(\chi)-n\,\pi(\chi)}{\sigma\sqrt{n}} \leq x\right) - \fN(x) - \frac{P(x)}{\sqrt{n}}\fn(x)\right| = o(n^{-1/2}),\quad\text{as}\quad n\to \infty\,.
    \end{equation*}
\end{thm}

\begin{rem}
    Note that from \eqref{eq: Edge main cond} and \eqref{eq: bound ab} with the corresponding choices of $a$ and $b$ it follows that $\chi \in L^3$. So, $\EXP_{{\overline{m}}}(|S_n(\chi)|^3) < \infty$ for each $n$. Our proof shows that the third asymptotic moment 
    $$\lim_{n\to \infty}\EXP_{{\overline{m}}}\left(\frac{S_n(\chi)-n \,\pi(\chi)}{\sqrt{n}}\right)^3$$ does, indeed, exist. 
\end{rem}

Finally, we provide a concrete example of a class of observables that satisfies our conditions. 

\begin{exmp}\label{ex: osc function}
 Let $\chi(x)=x^{-c}\sin(1/x)$ and define $\tilde{\eta}=\min\left\{1,\frac{\log\eta_-}{\log\eta_+}\right\}$.
 \begin{enumerate}[leftmargin=*, label=$(\arabic*)$]
     \item If $0\leq c<\min\big\{\sqrt{1+\tilde{\eta}}-1, \vartheta \tilde{\eta}\big\}$, 
     then $S_n(\chi)$ satisfies the CLT and MLCLT.
     \item If $0\leq c<\min\{\sqrt{1+\tilde{\eta}/6}-1,{\vartheta}\tilde{\eta}/6\}$, then $S_n(\chi)$  admits the first order Edgeworth Expansion.
 \end{enumerate}
 If $\psi$ is the doubling map, i.e.\ $\psi(x)=2x\mod 1$, then the conditions simplify in the following way:
 \begin{enumerate}[leftmargin=25pt, label=$(\arabic*a)$]
     \item If $c<\sqrt{2}-1\,(\approx 0.414)$, then $S_n(\chi)$  satisfies the CLT and MLCLT. 
     \item If $c<\sqrt{7/6}-1\,(\approx 0.080)$, then $S_n(\chi)$ admits the first order Edgeworth Expansion.
 \end{enumerate}
\end{exmp}

\subsection{The application to the Boolean-type transformation}\label{subsec: Boolean trafo}
For the following we define the Boolean-type transformation $\phi\colon \mathbb{R}\to\mathbb{R}$ as
\begin{align}
 \phi(x)\coloneqq\begin{cases} \frac{1}{2}\left(x-\frac{1}{x}\right)&\text{if }x\neq 0\\
               0&\text{if }x=0
              \end{cases}\label{eq: def T}
\end{align}
and $\mu$ the $\phi$-invariant probability measure absolutely continues with respect to Lebesgue and defined by
\begin{align}
 \mathrm{d}\mu\left(x\right)\coloneqq \frac{1}{\pi\cdot\left(x^2+1\right)}\,\mathrm{d}\lambda(x).\label{eq: def measure}
\end{align}
We are interested in limit theorems for Birkhoffs sums  $\widetilde{S}_n(h) :=S_n(h,\phi)$ where $h: \mathbb{R}\to\mathbb{R}$. 
To study these systems we go back to an easier system which fulfills all our properties of the last section. 

Let $\psi\colon I\to I$ be given by $\psi\left(x\right)\coloneqq 2x \mod 1$
and $\xi\colon I\to\mathbb{R}$ be given by $\xi\left(x\right)\coloneqq \cot\left(\pi x\right)$.
Note that $\xi$ is almost surely bijective. An elementary calculation yields that the dynamical systems 
$\left(\mathbb{R}, \mathcal{B}_{\mathbb{R}}, \mu, \phi\right)$
and 
$\left(I, \mathcal{B}_{I}, \lambda_I, \psi\right)$
are isomorphic via $\xi$, i.e.\
\begin{align*}
\left(\phi\circ\xi\right)\left(x\right)=\left(\xi\circ\psi\right)\left(x\right),
\end{align*}
for all $x\in I$
and additionally $\xi$ and $\xi^{-1}$ are measure preserving, i.e.\
for all $B\in\mathcal{B}_{\mathbb{R}}$ it holds that $\mu\left(B\right)=\lambda_I\left(\xi^{-1} B\right)$
and for all $B\in\mathcal{B}_{I}$ it holds that 
$\lambda_I\left(B\right)=\mu\left(\xi B\right)$. To simplify the notation, we define $\widetilde{\sigma}^2:=\sigma^2(h,\phi)$. 

Hence, instead of studying the Birkhoff sum 
$\sum_{n=0}^{N-1}\left(h\circ \phi^n\right)\left(x\right)$
with $x\in \mathbb{R}$
we can study the sum 
$\sum_{n=0}^{N-1}\left(h\circ \xi\circ\psi^n\right)\left(y\right)$,
for $y\in I$. 
Since the transformations $\phi$ and $\psi$ are isomorphic we conclude that 
\begin{align}\label{eq:LambdaChi}
 \mu\left(\sum_{n=0}^{N-1}\left(h\circ \phi^n\right)\left(x\right)\in B\right)
=\lambda_I\left(\sum_{n=0}^{N-1}\left(h\circ \xi\circ\psi^n\right)\left(y\right)\in B\right),
\end{align}
for all sets $B\in\mathcal{B}_{\mathbb{R}}$. Formally, we define $
\chi \colon I\to\mathbb{R}$ by $\chi\left(x\right)\coloneqq \left(h\circ \xi\right)\left(x\right)$ 
and consider then the Birkhoff sum $S_n(\chi).$ Then our task reduces to transferring the conditions we have for $\chi$ 
to conditions for $h$. 

Let $\fF$ be the class of functions  $h:\mathbb{R}\to\mathbb{R}$
  such that the left and right derivatives exist and there exist $u,v \geq 0$ 
 fulfilling  \begin{align}
  h\left(x\right)&\lesssim \left|x\right|^{u}\qquad\text{ and }\qquad
  \max\left\{\left|h'\left(x-\right)\right|, \left|h'\left(x+\right)\right|\right\}\lesssim \left|x\right|^{v}\label{eq: boundedness condition}
 \end{align}
 and $u(2+v)<1\,$. Analogously to $\overline{f}\,,$ we define $\overline{h}=h-\,\mu(h)\,.$ 
Easy examples for functions $h\in\fF$ are $h(x)= x^{a}$ with $0\leq a< (\sqrt{5}-1)/2$ or $h(x)= x^{a} \sin(x^b)$ with $a,b>0$ and $a(1+a+b)<1$. 

Under the non-coboundary condition on $\phi$, we have the CLT:
\begin{prop}\label{prop:CLTonReals}
Suppose $h \in \fF$ is not $\phi-$cohomologous to a constant. Then, the following CLT holds:~
  \begin{equation}\label{eq:RealCLT}
   \mu\left(\frac{\widetilde{S}_{n}(h) - n\,\mu(h)}{\widetilde{\sigma}\sqrt{n}} \leq x\right) - \fN(x)= o(1),\quad\text{for all}\quad x \in \reals\quad\text{as}\quad n\to \infty\,
    \end{equation}
with $\widetilde{\sigma}^2 \in (0,\infty)$.
\end{prop}

Under a non-arithmeticity condition on $\phi$, we have the MLCLT:
\begin{prop}\label{prop:MLCLTonReals}
 Let $h\in \fF$ be non-arithmetic. Let $0<\alpha_0<\alpha_1<\beta$ and $M\geq 1$.\, Then, the following MLCLT holds:~for  $V: \reals \to \reals$ compactly supported and continuous, $U$ such that $U \circ \xi \in \VVc{\alpha_0}{\beta}{M}\,,$  $W$ such that $W\circ \xi \in L^1$ for all ${\overline{m}}\in\cM_1(\reals)$ being absolutely continuous with respect to $\lambda$ such that $(W\circ \xi  \cdot \xi_*{\overline{m}}) \in \VVc{\alpha_1}{\beta}{M}^{\,\prime}\,\,,$ 
 we have
   \begin{equation}\label{eq:RealLLT}
   \lim_{n \to \infty} \sup_{\ell \in \reals }\left|\sigma\sqrt{2\pi n}\,\EXP_{{\overline{m}}}(U \circ \psi^n \,V(\widetilde S_n(\overline{h} )-\ell)\, W)  - e^{-\frac{\ell^2}{2n\widetilde{\sigma}^2}} \,\EXP_{{\overline{m}}}(W)\,\EXP_{\mu}(U)\int V(x) \, \mathrm{d}x \right|= 0\,.
   \end{equation}
\end{prop}

Finally, we state a set of sufficient conditions that implies the Edgeworth Expansions for $\phi$.
\begin{prop}\label{prop:EdgeworthonReals}
 Let $h:\mathbb{R}\to\mathbb{R}$
 be such that the left and right derivatives exist and there exist  $u,v \geq 0$ 
 fulfilling \eqref{eq: boundedness condition} and 
 \begin{align}
  \min\{2u(2+v), (u+v)(2+v) \}<1/3\label{eq: Edge Boole main cond}
 \end{align}
 and $h$ is not arithmetic. Then there exists a quadratic polynomial $P$ whose coefficients depend on the first three asymptotic moments of $\widetilde{S}_n(h)$ but not on $n$ such that for all ${\overline{m}}\in\cM_1(\reals)$ being absolutely continuous with respect to $\lambda$ we have
  \begin{equation}\label{eq:Edgeworth}
\sup_{ x\in \reals} \left|{\overline{m}}\left(\frac{ \widetilde{S}_n(h)-n\,\mu(h)}{\widetilde\sigma\sqrt{n}} \leq x\right) - \fN(x) - \frac{P(x)}{\sqrt{n}}\fn(x)\right| = o(n^{-1/2}),\quad\text{as}\quad n\to \infty\,.
    \end{equation}
\end{prop}

\begin{rem}
The condition \eqref{eq: Edge Boole main cond} forces that $0 \leq u < 1$ and $u < v$.
\end{rem}

\subsection{Sampling the Lindel\"of Hypothesis}\label{sec:Lindelof} 
In this section, we apply the results from the last subsection to the context of sampling the Lindel\"of Hypothesis. Let $\zeta : \complex\setminus \{1\} \to \complex$ be the Riemann zeta function defined by $$\zeta(s):=\sum_{n=1}^\infty n^{-s},\,\,\,\Re(s) >1$$ and by analytic continuation elsewhere except $s=1$. The Lindel\"of hypothesis states that the Riemann zeta function does not grow too quickly on the critical  line $\Re z = 1/2$. 
More precisely, it is conjectured that 
\begin{equation}\label{eq:zetadecay}
    \zeta_{1/2}(t):=\zeta\left(\frac{1}{2}+it\right) = \cO(t^{\eps})\,,\quad t\to{\pm} \infty
\end{equation} 
for all $\eps > 0\,,$ i.e., $\lim_{t\to\pm\infty}|\zeta_{1/2}\left(t\right)|/t^\eps<\infty$. To date, the best estimates are due to Bourgain in \cite{Bo} where it is proved that this is true for all $\eps> 13/84 \approx 0.154$. It is worth noting that the Riemann hypothesis implies the Lindel\"of hypothesis and the latter is a substitute for the former in some applications.

The conjecture is related to the value distribution of $\zeta_{1/2}(t)$ as $t\to \pm\infty$. In order to obtain more information about this tail behaviour, one can study ergodic averages of $\zeta_{1/2}$ sampled over the orbits of heavy-tailed stochastic processes. This approach to Lindel\"of hypothesis was initiated by Lifschitz and Weber in \cite{LW}. In particular, they prove that when $\{Y_j\}_{j\geq 0}$ are independent Cauchy distributed random variables and $X_k = \sum_{j=0}^{k-1} Y_j$ (the Cauchy random walk), then for all $b>2$,
$$\frac{1}{n}\sum_{k=0}^{n-1}\zeta_{1/2}(X_k) = 1 + o\left(\frac{(\log n)^b}{\sqrt{n}}\right)\,,\, n \to \infty,$$
almost surely, where we denote $a_n=o(b_n)$ if $\lim_{n\to\infty}|a_n|/b_n=0$. This work was later generalised by Shirai, see \cite{Shi}, where
$X_k$ was taken to be a symmetric $\alpha$-stable process with $\alpha \in [1, 2)$.
Since $X_k$ are heavy tailed, i.e., $\EXP(|X_k|^p)=\infty$ when $p = \lceil \alpha \rceil$ (including the Cauchy case $\alpha=p=1$), the $\alpha$-stable process samples large values with high probability. So, this result illustrates that the values of $\zeta_{1/2}(t)$ are small on average even for large values of $|t|$. 

Similarly, in the deterministic setting, the Birkhoff sums 
\begin{equation}\label{eq:BSum}
   \sum_{k=0}^{n-1}\zeta_{1/2}(\phi^kx)
\end{equation}
where
$\phi:\reals \to \reals$ is the Boolean-type transformation given in \eqref{eq: def T},
are studied in \cite{steuding}. Since $\phi$ preserves the ergodic probability measure $\mu$ given in \eqref{eq: def measure}
(the law of a standard Cauchy random variable) and $\zeta_{1/2}$ is integrable  with respect to $\mu$, it follows from Birkhoff's pointwise ergodic theorem that for almost every (a.e.) $x \in \reals$
\begin{equation}\label{eq:SLLN}
    \lim_{n\to \infty}\frac{1}{n}\sum_{k=0}^{n-1}\zeta_{1/2}(\phi^kx) = \int \zeta_{1/2}(x) \frac{\mathrm{d}x}{\pi(1+x^2)} = \zeta_{1/2}(3/2) - 8/3   \approx -0.054\,.
\end{equation}
This too illustrates that most of the values of $\zeta_{1/2}$ are not too large.

Sampling the Lindel\"of hypothesis has two other theoretical underpinnings. On the one hand, it is known that the Lindel\"of hypothesis is true if and only if for all $m\in \naturals$ and for a.e.\,\,$x\in \reals$, the following limit exists 
\begin{equation*}
    \lim_{n\to \infty}\frac{1}{n}\sum_{k=0}^{n-1}|\zeta_{1/2}(\phi^kx)|^{2m}=\int |\zeta_{1/2}(x)|^{2m} \frac{\mathrm{d}x}{\pi(1+x^2)}.
\end{equation*}
On the other hand, the Riemann hypothesis is true if and only if for a.e.\,\,$x \in \reals$ 
\begin{equation*}
    \lim_{n\to \infty}\frac{1}{n}\sum_{k=0}^{n-1}{\log}|\zeta_{1/2}((\phi^kx)/2)| = 0.
\end{equation*}
In both cases, evidence can be gathered numerically,  
see \cite[Theorems 4.1 and 4.2]{steuding} for details.

The results by Steuding have also been generalised, both by replacing $\zeta$ and replacing $\phi$:~in \cite{elas}, Elaissaoui and Guennoun used $\log|\zeta|$ as the observable and a slight variation of $\phi$ as the transformation, and in \cite{leesuriajaya}, Lee and Suriajaya studied different classes of meromorphic functions such as Dirichlet $L-$functions or Dedekind $\zeta$ functions while taking $\phi$ to be an affine version of the Boolean-type transformation. Maugmai and Srichan gave further generalizations of these results, see \cite{MS}. It must also be mentioned that these transformations $\phi$ have been studied earlier in a solely ergodic theoretic context by Ishitani(s) in \cite{I1, I2}.

In what follows, we will use the results of the last subsection to further understand the value distribution of the Birkhoff averages  on the critical strip given by \eqref{eq:BSum} around their asymptotic mean $A=\zeta_{1/2}(3/2) - 8/3$ 
 and also study the Birkhoff averages of $\zeta(s + i\,\cdot\,)$ for other values $s\in (0,1)$ of the critical strip.
In particular, we will state a CLT and MLCLT for the above setting.

Recall $\fF$ the class of functions $h:\mathbb{R}\to\mathbb{R}$
  such that the left and right derivatives exist and there exist $u,v \geq 0$ 
 fulfilling  \begin{align*}
  h\left(x\right)&\lesssim \left|x\right|^{u}\qquad\text{ and }\qquad
  \max\left\{\left|h'\left(x-\right)\right|, \left|h'\left(x+\right)\right|\right\}\lesssim \left|x\right|^{v}
 \end{align*}
 and $u(2+v)<1\,$. Since $|\Re\,\zeta{ (s+i\,\cdot)}|^a, |\Im\,\zeta{ (s+i\,\cdot)}|^a, |\zeta{ (s+i\,\cdot)}|^a \in  \fF$ for some suitable choices of $s$ and $a$, we obtain two corollaries that improve the existing results on sampling the Lindel\"of hypothesis. 

\begin{cor}\label{cor:CLTforRZ}
    Let $s \in (3 - 2 \sqrt{2}, 1)$ and define $h : \reals \to \reals$ as follows.
\begin{itemize}
    \item $h(x) = \Re\,\zeta(s+ix)\,,$
    \item $h(x) = \Im\,\zeta(s+ix)\,, $ or
    \item $h(x) = |\zeta(s+ix)|$\,.
\end{itemize}
If $h$ is not $\phi-$cohomologous to a constant, then the CLT, \eqref{eq:RealCLT} holds. Moreover, if $h$ is non-arithmetic, then the MLCLT, \eqref{eq:RealLLT}, holds.
\end{cor}

\begin{rem}
See \cite[Section 2.5]{STI} for a discussion where it is shown using numerics that for $\zeta_{1/2}$ all of the above choices of $h$ are not coboundaries. Similarly, for a fixed value of $s$, one can numerically check whether $h$ is not a $\psi-$coboundary by calculating the sum of values of $\chi = h \circ \xi $ over some appropriate periodic orbit of the doubling map and showing that it is not equal to $0$. 

Moreover, by the Cram\'er-Wold theorem the first two statements in Corollary \ref{cor:CLTforRZ} also imply a complex central limit theorem. 
\end{rem}

\begin{cor}\label{cor:CLTforpowerRZ}
   Let $h : \reals \to \reals$ be as follows.
\begin{itemize}
   \item $h = |\mathfrak{R}\,\zeta_{1/2}|^a\,,$
      \item $h = |\mathfrak{I}\,\zeta_{1/2}|^a, $ or
      \item $h = |\zeta_{1/2}|^a$
\end{itemize}
where $1 \leq  a < 84/13(\sqrt{2}-1) \, (\approx 2.677) $. 
If $h$ is not $\phi-$cohomologous to a constant, then the CLT, \eqref{eq:RealCLT} holds. Moreover, if $h$ is non-arithmetic, then the MLCLT, \eqref{eq:RealLLT}, holds.
\end{cor}

\begin{rem}
The best estimates in the literature for $\varepsilon$ in \eqref{eq:zetadecay} (for example, \cite{Bo}) are not sufficient to conclude that the Riemann zeta function, more precisely $\Re\zeta_{1/2}$, $\Im \zeta_{1/2}$, and $|\zeta_{1/2}|$, satisfy the conditions of \Cref{prop:EdgeworthonReals} on the existence of the first order Edgeworth expansion, albeit a slight improvement of results in \cite{Bo} will provide us what is required. In fact, our theorem shows that if the Lindel\"of hypothesis is true, then the first order Edgeworth expansion has to hold.
\end{rem}

\begin{rem}
On the one hand, the Lindel\"of hypothesis states that $|\zeta_{1/2}(x)|\lesssim x^{\eps}$ holds for all $\eps>0$, and hence, if it is true, the above statement of Corollary \ref{cor:CLTforpowerRZ} has to hold for any $a>0$. On the other hand, sampling $|\zeta(s+i\phi^k(x))|^a$ with larger values of $a$ and obtaining normally distributed samples provides further evidence that the Lindel\"of hypothesis is indeed true.
The same holds for the first order Edgeworth expansion: Under the condition that $h$ is non-arithmetic with $h$ as in Corollary \ref{cor:CLTforRZ} and assuming the Lindel\"of hypothesis holds, also a first order Edgeworth expansion has to hold. A numerical simulation is not part of this paper. However, observing convergence or not gives a further hint that the Lindel\"of hypothesis holds or not.
\end{rem}

\section{Review of Abstract Results for Limit Theorems}\label{sec:AbsRes}
One known technique used to establish limit theorems for ergodic sums with unbounded observables is a combination of the Keller-Liverani perturbation result (see \cite{KL}) applied to a sequence of Banach spaces as in \cite{HervePene, FP, Pene}. We have stated elementary criteria for the CLT and the MLCLT to exist below as propositions adapted from \cite[Corollary 2.1, Theorem 5.1]{HervePene} to our setting.  

\begin{prop}\label{prop:AbsCLT}
   Let $T: X \to X$ be a 
   dynamical system that has an ergodic invariant probability
   measure $\widetilde m$. Let ${f} \in L^2(\widetilde m)$ be such that $\widetilde m({f})=0$ and $\sum_{n\geq 0} \wh{ T}^n ({f}) $ converges in $L^2(\widetilde m)$. Then, we have the following CLT. 
   \begin{equation}\label{eq:AbsCLT}
  \lim_{n\to\infty} \widetilde m\left(\frac{S_{n}({f})}{\sqrt{n}} \leq x\right) = \fN\left(\frac{x}{\sigma}\right),\quad\text{for all}\quad x \in \reals\quad\text{as}\quad n\to \infty\,,
    \end{equation}
where $\sigma^2=\sigma^2(f,T)$ can be written as
$$\sigma^2 = \EXP_{\widetilde m}({f}^2)+ 2\sum_{k=1}^\infty \EXP_{\widetilde m}({f}\cdot {f}\circ { T}^k) \in [0,\infty)\,.$$ 
Here $\sigma = 0$ if and only if $f$ is a $T$-coboundary and in this case $\fN(x/\sigma):={\bf 1}_{[0,\infty)}$ and $\frac{S_{n}({ f})}{\sqrt{n}} \to \delta_0$ in distribution as $n \to \infty$. 
\end{prop}
\begin{proof}
This follows due to Gordin \cite{Gordin}. See \cite[Corollary 2.1, Proposition 2.4]{HervePene} for details.
\end{proof}

\begin{prop}\label{prop:AbsLLT}
Let $T: X\to X$ be a non-singular
   dynamical system  wrt a probability measure $m$. Suppose $T$ has an ergodic invariant probability measure $\widetilde m$ absolutely continuous wrt $m$ and that
there exist two, not necessarily distinct, Banach spaces  $\cX$ and $\cX^{(+)}$ such that
\begin{equation}\label{eq:TwoSpaces}
    \cX\hookrightarrow\cX^{(+)}\hookrightarrow L^1(\pi)
\end{equation}
each containing ${\bf 1}_X$ and satisfying the following:
\begin{enumerate}[topsep=0pt, label=\normalfont{(\Roman*)}]
\item \label{cond1A}  For all $s \in \reals$, ${\wh { T}}_{is} \in \cL(\cX) \cap \cL(\cX^{(+)})\,.$
\item \label{cond2A} The map  $s \mapsto {\wh { T}}_{is}\in\mathcal L(\cX,\cX^{(+)})$ is continuous on  $\reals\,.$
\item \label{cond3A}  Either $\cX = \cX^{(+)}\,,$ or
there  exist $\kappa \in (0,1)$ and $\delta>0$ such that for all $$ z \in D_{\kappa}:=\{z \in \complex||z|>\kappa,|z-1|>(1-\kappa)/2\},$$  and for all $s \in (-\delta, \delta)$ we have $$(z\Id-{\wh { T}}_{is})^{-1} \in \cL(\cX)\,\quad\,\text{and}\,\quad\,\sup_{|s|<\delta} \sup_{z \in D_\kappa} \|(z\Id-{\wh { T}}_{is})^{-1}\|_{\cX\to \cX}< \infty \,.$$
\item \label{cond4A} $\lim_{n\to \infty}\|\wh{ T}^n(\cdot) - \widetilde m(\cdot){\bf 1}_X\|_{\cX_0 \to \cX_0 } = 0\,.$
\item \label{cond5A} The CLT, \eqref{eq:AbsCLT} holds with $\sigma >0$.
\item \label{cond6A} For all $s \neq 0$, the spectrum of the operators ${\wh { T}}_{is}$ acting on $\cX$ is contained in the open unit disc, $\{z\in \complex\ |\ |z|<1\}\,.$  
\end{enumerate} 
Then, for all $U\in \cX\,,$ $V: \reals \to \reals$ a compactly supported continuous function, ${\overline{m}} \in \cM_1(X)$ being absolutely continuous wrt $m$ and $W \in L^1$ such that $(W \cdot  {\overline{m}} ) \in \cX^{(+)}\,'\,,$  we have
    \begin{equation}\label{eq:LLT}
   \lim_{n \to \infty} \sup_{\ell \in \reals }\left|\sigma\sqrt{2\pi n}\,\EXP_{{\overline{m}}}(U \circ { T}^n \,V(S_n(\overline\chi)-\ell)\, W)  - e^{-\frac{\ell^2}{2n\sigma^2}} \,\EXP_{{\overline{m}}}(U)\,\EXP_{\widetilde m}(W)\int V(x) \, \mathrm{d}x \right|= 0\,.
   \end{equation}
\end{prop}
\begin{proof}
This follows from a modified version of \cite[Theorem 5.1]{HervePene}. 
The condition (CLT) there is assumed here in \ref{cond5A}. 

Also, the Condition ($\widetilde{K}$) there follows from our assumptions { \ref{cond1A}} through  \ref{cond4A} because $(K1)$ is \ref{cond4A}, $(\widetilde{K1})$ is \ref{cond2A}, and finally, $(\widetilde{K2})$ can be replaced by \ref{cond3A} (see \Cref{rem:ReplaceIV}). 

Our assumptions \ref{cond2A} and \ref{cond6A} yield that on any compact set $K \subset \reals\setminus\{0\} $, there exist $\rho \in (0,1)$ and $C_K>0$ such that 
$$\sup_{s \in K} \|\wh{ T}^n_{is}\|_{\cX \to \cX^+} \leq C_K \rho^n,$$ 
for all $n \in \naturals$ (see, for example, \cite[Proposition 1.13]{FP} for a proof). This replaces the non-lattice condition $(S)$ there. 

So, for all $U\in \cX\,,$ $V: \reals \to \reals$ a compactly supported continuous function and $W \in L^1$ such that $(W \cdot {\overline{m}}) \in \cX^{(+)}\,'\,,$ we have the MLCLT due to \cite[Theorem 5.1]{HervePene}.
\end{proof}

Finally, we state a result that gives us sufficient conditions for the 
first order Edgeworth expansion. It is adapted from \cite{{HervePene}, FP} to our setting (compare with \cite[Proposition 7.1, Propostion A.1]{{HervePene}} and \cite[Corollary 1.8, Proposition 1.12]{FP}). 
\begin{prop}\label{prop:AbstractThm}
Let $T: X \to X$ be a non-singular dynamical system wrt a probability measure $m$. Suppose $T$ has an ergodic invariant probability measure $\widetilde m$ absolutely continuous wrt $m$ and that there exists a sequence of, not necessarily distinct, Banach spaces 
\begin{equation}\label{eq:SpaceChain}
    \cX_0\hookrightarrow\cX^{(+)}_0\hookrightarrow\cX_1\hookrightarrow\cX^{(+)}_1\hookrightarrow \cX_{2}\hookrightarrow \cX^{(+)}_{2} \hookrightarrow\mathcal X_{3}\hookrightarrow \cX_{3}^{(+)}
\end{equation}
each containing ${\bf 1}_X$, $\cX_{3}^{(+)}\hookrightarrow L^1$ and satisfying the following:
\begin{enumerate}[topsep=0pt, label=\normalfont{(\Roman*)}]
\item \label{cond1}  For each space $\cC$ in \eqref{eq:SpaceChain}$,$ $s \in \reals$, ${\wh { T}}_{is} \in \cL(\cC)\,.$
\item \label{cond2} For all $a= 0,1,2,3$, the map  $s \mapsto {\wh { T}}_{is}\in\mathcal L(\cX_a,\cX_{a}^{(+)})$ is continuous on  $\reals\,.$
\item \label{cond3} For all $a= 0,1,2$, 
the map $s \mapsto {\wh { T}}_{is}\in\mathcal L(\cX_a^{(+)},\cX_{a+1})$ is $C^{1}$ on 
$(-\delta,\delta)\,.$
\item\label{cond4}  Either all spaces in \eqref{eq:SpaceChain} are equal$,$ or
 there  exist $\kappa \in (0,1)$ and $\delta>0$  such that for all $$ z \in D_{\kappa}:=\{z \in \complex||z|>\kappa,|z-1|>(1-\kappa)/2\},$$ for all $s \in (-\delta, \delta)$ and for each space $\cC$ in \eqref{eq:SpaceChain}$,$ $$(z\Id-{\wh { T}}_{is})^{-1} \in \cL(\cC)\,\quad\,\text{and}\,\quad\,\sup_{|s|<\delta} \sup_{z \in D_\kappa} \|(z\Id-{\wh { T}}_{is})^{-1}\|_{\cC\to \cC}< \infty \,.$$
\item \label{cond5} $\wh { T}$ has a spectral gap of $(1-\kappa)$ on each space $\cC$ in \eqref{eq:SpaceChain}.
\item \label{cond6} For all $s \neq 0$, the spectrum of the operators ${\wh { T}}_{is}$ acting on either $\cX_0$ or $\cX_0^{(+)}$ is contained in the open unit disc, $\{z\in \complex\ |\ |z|<1\}\,.$  
\item \label{cond7} The sequence $$\left\{\sum_{k=0}^{n-1} {\overline{ f}} \circ { T}^k\right\}_{n \in \naturals}$$ where $ \overline{ f} := { f} - A$ has an $L^2-$weakly convergent subsequence\,.
\item \label{cond8} $f$ is not $T-$cohomologous to a constant.
\end{enumerate}
Then  for all ${\overline{m}} \in \cM_1(X)$ being absolutely continuous wrt $m$,
    there exists a quadratic polynomial $P$ whose coefficients depend on the first three asymptotic moments of $S_n(\chi)$ such that the following asymptotic expansion holds;
    \begin{equation}\label{eq:EdgeExp}
    \sup_{ x\in \reals} \left|\widetilde m\left(\frac{S_{n}(\overline{ f})}{\sigma\sqrt{n}} \leq x\right) - \fN(x) - \frac{P(x)}{\sqrt{n}}\fn(x)\right| = o(n^{-1/2}),\quad\text{as}\quad n\to \infty\,.
    \end{equation}
\end{prop}
\begin{rem}\label{rem:ReplaceIV}
In \cite{HervePene} and \cite{FP}, instead of the condition \ref{cond4} above, the following stronger condition of a uniform DFLY inequality is assumed. 
\begin{displayquote}
    Either all spaces in \eqref{eq:SpaceChain} are equal, or
there exist $\widetilde C>0$,  $\widetilde\kappa_1\in (0,1)$ and $p_0 \geq 1$ such that, for every $\cC$ in \eqref{eq:SpaceChain},
\begin{equation}\label{eq:EqualHyp}
    \forall h\in\cC,\quad \sup_{|s|<\delta}\Vert \wh { T}_{is}^nh\Vert_{\cC}\le
\widetilde C\left(\widetilde\kappa_1^n\Vert h\Vert_{\cC}+\Vert h\Vert_{L^{p_0}(\bar\nu)}\right)\, .
\end{equation}
\end{displayquote}
However, the proof of the key theorem, \cite[Proposition 1.11]{FP}, is based on \cite[Proposition A, Corollary 7.2]{HervePene} which uses the hypothesis $\cD(m)$ in \cite[Appendix A]{HervePene} that contains the much weaker condition \ref{cond4} instead of condition \eqref{eq:EqualHyp}. Therefore, all the results in \cite{FP} based on \cite[Proposition 1.11]{FP} including \cite[Proposition 1.12]{FP} remain true with this replacement. We refer the reader to \cite{HervePene} for more details. 
\end{rem}

\begin{rem}
For an elementary illustration of the proof of the CLT based on the classical Nagaev-Guivarc'h approach, we refer the reader to \cite{G} where the $C^2$ regularity of $s \mapsto \wh { T}_{is}$ along with the spectral gap of $\wh { T}$ on a single Banach space (instead of a chain) is used. This corresponds to the $C^2$ regularity of the characteristic function in the IID case. When it comes to the MLCLT in the IID setting, a non-lattice assumption is necessary. In our case, the equivalent assumption is \ref{cond6}. 
\end{rem}

\begin{proof}[Proof of Proposition \ref{prop:AbstractThm}]
We apply results in \cite{FP} restricted to a single dynamical system with $r=1$ there, i.e., when Assumptions (0) and (A)[1](1-2) in \cite[Section 1.2]{FP} are trivially true. This case is, thus, similar to the $r=1$ case of \cite[Proposition 1.12]{FP} which implies \cite[Corollary 1.8]{FP} which, in turn, gives the first order Edgeworth expansion. This is because our assumptions above imply Assumptions (A)[1] and (B) in \cite[Section 1.2]{FP}, \textit{except} for (A)[1](4) which is equivalent to \eqref{eq:EqualHyp}. However, as discussed in \Cref{rem:ReplaceIV},  \cite[Corollary 1.8]{FP} remains true because the key ingredient of the proof in \cite{FP} is our assumption (IV) (implied by the much stronger (A)[1](4)). 
\end{proof}

\section[Twisted Transfer Operators]{Twisted Transfer Operators \texorpdfstring{$\wh \psi_{is}$}{hat psi is}}\label{sec:Operators}

\subsection{Properties of twisted transfer operators}

We first prove $L^\gamma$ norm estimates for $\wh \psi_{is}$.

\begin{lem}\label{lem:LpEst} For all $\gamma >1$,  $s\in \reals$ and $\varphi \in L^{\gamma}$, there exists a constant $C_{\gamma} > 1$ that depends only on $\psi$ and $\gamma$ such that 
$$\|\wh \psi_{is} (\varphi)\|_1 \leq \|\wh \psi_{is} (\varphi)\|_{\gamma} \leq C_\gamma \|\varphi\|_\gamma\, .$$
\end{lem}
\begin{proof}
The first inequality follows from a direct application of H\"older's inequality. The second one is a straightforward application of Minkowski's inequality. 
\begin{align*} \left(\int|{\wh\psi}_{is}(\varphi)|^\gamma \, \mathrm{d}\lambda_I\right)^{1/\gamma} &\leq \left(\int \wh\psi(|\varphi|)^\gamma \, \mathrm{d}\lambda_I\right)^{1/\gamma}\\ &= \left(\int \left( \sum_{j=0}^{k-1}  \frac{|\varphi|\circ \psi^{-1}_{j+1}}{|\psi'\circ \psi^{-1}_{j+1}|}\, \right)^{\gamma} \mathrm{d}\lambda_I \right)^{1/\gamma} \\ &\leq \sum_{j=0}^{k-1}\left( \int\left( \frac{|\varphi|\circ \psi^{-1}_{j+1}}{|\psi'\circ \psi^{-1}_{j+1}|}\right)^{\gamma}\,  \mathrm{d}\lambda_I \right)^{1/\gamma}\\ &= \sum_{j=0}^{k-1} \left( \int\left(\frac{|\varphi|}{|\psi'|}\right)^{\gamma}\One_{[c_j,c_{j+1}]}|\psi'|\, \mathrm{d}\lambda_I\right)^{1/\gamma} \\ &\leq \frac{k}{\eta_{-}^{1-\gamma}}\left(\int |\varphi|^{\gamma}\,\mathrm{d}\lambda_I\right)^{1/\gamma}\,.
\end{align*}
Put $C_\gamma = k\cdot \eta_{-}^{\gamma-1}$. Then 
$$\|{\wh \psi}_{is}(\varphi)\|_{\gamma} \leq C_\gamma \|\varphi\|_\gamma\, .
$$
\end{proof}

Next, we have the following result on the required regularity of the transfer operators. 

\begin{cor}\label{cor:RegOpr} 
Let $0 \leq \alpha_0,
\alpha^*,  \alpha^{**},\beta\leq 1 $ and $ \gamma_0, \gamma\geq 1$. 
Put \begin{align*}
    \alpha_1 &= \alpha_0+\alpha^* &
    \alpha_2 & =\alpha_1+\max\{\alpha^{**},\alpha^*\}\\ 
    1 \leq \gamma_1 &\leq \gamma_0  & 1 \leq \gamma_2 &\leq (\gamma_1^{-1}+\gamma^{-1})^{-1}
\end{align*}
and consider the chain of Banach spaces
\begin{equation}\label{SpaceChain}
    \VVc{ \alpha_0}{ \beta}{\gamma_0}\hookrightarrow \VVc{ \alpha_1}{ \beta}{\gamma_1} \hookrightarrow  \VVc{ \alpha_2}{ \beta}{\gamma_2}\,.
\end{equation}
Suppose that for all $s \in \reals$,  $|e^{is\chi}|_{0, \beta} < \infty$. Then
\begin{enumerate}
\item[$(1)$]  for $s\in \reals$, $\wh \psi_{is}$ is a bounded linear operator on each of the Banach spaces in \eqref{SpaceChain}.
\end{enumerate}
Suppose, in addition, that $\lim_{s \to 0}|{1-e^{is\chi}}|_{\alpha^*,\beta} = 0.$ Then 
\begin{enumerate}
\item[$(2)$] $s \mapsto \wh \psi_{is}$ is continuous as a function from $\reals$ to $\cL(\VVc{ \alpha_0}{ \beta}{\gamma_0},\VVc{ \alpha_1}{ \beta}{\gamma_1})$.
\end{enumerate}
Finally, suppose that
    \begin{equation*}
      \lim_{s \to 0}\left|\frac{e^{is\chi}-1-is\chi}{s}\right|_{\alpha^{**},\beta}=0\quad\text{and}\quad \|\chi\|_{\gamma} < \infty.
    \end{equation*}
Then,
\begin{enumerate}
\item[$(3)$] $s \mapsto \wh \psi_{is}$ is continuously differentiable as a function from $\reals$ to $\cL(\VVc{ \alpha_1}{ \beta}{\gamma_1} , \VVc{ \alpha_2}{ \beta}{\gamma_2})$.
\end{enumerate}
\end{cor}

\begin{proof}
Since $\wh\psi$ is a bounded linear operator on each of the Banach spaces in \eqref{SpaceChain}  (in particular, due to the DFLY inequality below), the corollary follows from \Cref{lem:MultOp} and \Cref{lem:LpEst}. 
\end{proof}

\subsection{DFLY Inequalities}
In this section, we prove DFLY inequalities for the family $\wh \psi_{is}$. First, we state and prove two preparatory lemmas. Throughout this section, we assume that $\chi$ is continuous and the right and left derivatives of $\chi$ exist on $\mathring{I}$ and that there exists a constant $b>0$ such that 
\begin{align}
  \max\{|\chi'(x+)|, |\chi'(x-)|\}\lesssim x^{-b}(1-x)^{-b}\,.\label{eq: cond 3 cond a}
\end{align}
\begin{lem}\label{lem:DF1} Let $\alpha,\beta \in (0,1)$ and let $\bar\gamma \in [ 1, 1/\alpha)$. Suppose the constant $b>0$ in \eqref{eq: cond 3 cond a} is such that 
\begin{align}\label{eq:GammaAlphaBetaB}
    \min\big\{\bar\gamma^{-1}+(\alpha-\beta)b\,,\, \bar\gamma^{-1}+\alpha -\beta b\big\} >0\,.
\end{align} 
Then, there exists $C_{\eps_0}>0$ independent of $\bar{\gamma}$ 
such that 
\begin{equation}\label{eq:EXPoscBound}
    \sup_{\eps\in (0,\eps_0]}\eps^{-\beta}\left\|R_{\alpha}\osc\left(e^{i s\chi}, B_{\eps}(\cdot)\right)\right\|_{\bar\gamma} \leq C_{\eps_0}
\end{equation}
  for all $s \in \reals$. 
\end{lem}
\begin{rem}
We note that, if $b>1,$ then $\bar\gamma^{-1}+\alpha -\beta b >0\implies \bar\gamma^{-1}+(\alpha-\beta)b>0$, and if $b<1$ then $\bar\gamma^{-1}+(\alpha-\beta)b>0 \implies \bar\gamma^{-1}+\alpha -\beta b >0.$
\end{rem}
\begin{proof}[Proof of \Cref{lem:DF1}]
Since $e^{i s\chi}$ is $2\pi$ periodic in $s$, we will estimate $$\sup_{s\in [0,2\pi]}\sup_{\eps\in (0,\eps_0]}\eps^{-\beta}\|R_{\alpha}\osc(e^{i s\chi}, B_{\eps}(\cdot))\|_{\bar\gamma} \,.$$
Note that
\begin{align*}
 \sup_{\eps\in (0,\eps_0]}\|R_{\alpha}\osc(e^{i s\chi}, B_{\eps}(\cdot))\|_{\bar\gamma}\cdot \eps^{-\beta}
 &\leq \sup_{\eps\in (0,\eps_0]} \left(\int_0^{1/2} \left(R_{\alpha}\osc(e^{i s\chi}, B_{\eps}(x))\right)^{\bar\gamma}\,\mathrm{d}\lambda_I(x)\right)^{1/\bar\gamma}\cdot \eps^{-\beta}\\
 &\qquad+\sup_{\eps\in (0,\eps_0]} \left(\int_{1/2}^1 \left(R_{\alpha}\osc(e^{i s\chi}, B_{\eps}(x))\right)^{\bar\gamma}\,\mathrm{d}\lambda_I(x)\right)^{1/\bar\gamma}\cdot \eps^{-\beta}.
\end{align*}
We will only estimate the first summand as the estimation of the second follows analogously.
Using the definition $\osc(h,A)=\osc(\Re h,A)+\osc(\Im h,A)$ and $|e^{it_1}-e^{it_2}|\leq \min\{2,|t_1-t_2| \}$,
we note that for any measurable set $A$ we have $\osc\left(e^{is\chi}, A\right)\leq \min\{4, 4s/\pi \osc(\chi, A)\}$.
Due to \eqref{eq: cond 3 cond a} there exists $C>0$ such that for all $s>0$, for all $\eps>0$ and all $x\in [\eps,1/2]$ we have
\begin{align*}
 \osc(e^{i s\chi}, B_{\eps}(x))
 &\leq \frac{8|s|\eps}{\pi}\sup_{y\in B_{\eps}(x)}\max\{|\chi'(y+)|, |\chi'(y-)|\}
 \leq \frac{8C|s|\eps}{\pi} (x-\eps)^{-b}.
\end{align*}
We have that $8C|s|\eps (x-\eps)^{-b} /\pi \leq 4$ if and only if
\begin{align*}
 x\geq \left(\frac{2C|s|\eps}{\pi}\right)^{1/b}+\eps=: K_\eps > \eps.
\end{align*}
Since $K_\eps> \eps$, on $[K_\eps, 1/2]$, we use $\frac{8C|s|\eps}{\pi} (x-\eps)^{-b}$, and on $[0, K_\eps)$, we use $4$ as upper bounds for $\osc(e^{i s\chi}, B_{\eps}(x))$, to obtain
\begin{align}
 \MoveEqLeft\sup_{\eps\in (0,\eps_0]} \left(\int_0^{1/2} \left(R_{\alpha}\osc(e^{i s\chi}, B_{\eps}(x))\right)^{\bar\gamma}\,\mathrm{d}\lambda_I(x)\right)^{1/\bar\gamma}\cdot\eps^{-\beta}\notag\\
 &\leq \sup_{\eps\in (0,\eps_0]} \left( 4{K_\eps}\sup_{[0,K_\eps]}R_\alpha\One \cdot \eps^{-\beta}+\left(\int_{K_\eps}^{1/2} \left(\frac{8C|s|\eps^{1-\beta}}{\pi} R_\alpha\One \cdot (x-\eps)^{-b}\right)^{\bar\gamma}\mathrm{d}\lambda_{I}(x)\right)^{1/\bar\gamma}\right)\notag\\
 &\leq \sup_{\eps\in (0,\eps_0]} 4{K_\eps}^{1+\alpha}\eps^{-\beta}+\sup_{\eps'\in (0,\eps_0]}\frac{8C|s|\eps^{1-\beta}}{\pi}\,\left(\int_{K_\eps}^{1/2} \left( x^{\alpha} (x-\eps)^{-b}\right)^{\bar\gamma}\mathrm{d}\lambda_{I}(x)\right)^{1/\bar\gamma}.\label{eq: estim osc a}
\end{align}
For the first summand of \eqref{eq: estim osc a}, we have that there exists $\tilde C_{\eps_0}>0$ such that
\begin{align*}
 \sup_{\eps\in (0,\eps_0]} 4{K_\eps}^{1+\alpha}\eps^{-\beta}
 &\leq 8\sup_{\eps\in (0,\eps_0]} \max\left\{\left(\frac{2C|s|}{\pi}\right)^{(1+\alpha)/b}\eps^{(1+\alpha)/b-\beta},\eps^{1+\alpha-\beta} \right\}\\ &\leq \tilde C_{\eps_0}(1+|s|^{(1+\alpha)/b}) < \infty
\end{align*}
which follows from the fact that $\beta < (1/\bar\gamma+\alpha)/b < (1+\alpha)/b$ and $\beta\leq 1$.

For the second summand of \eqref{eq: estim osc a}, we use $\bar\gamma < 1/\alpha$ and $(x+\eps)^{\alpha \bar\gamma} \leq x^{\alpha \bar\gamma}+\eps^{\alpha \bar\gamma}$ to compute
\begin{align*}
\MoveEqLeft\sup_{\eps\in (0,\eps_0]}\frac{8C|s|\eps^{1-\beta}}{\pi}\,\left(\int_{K_\eps}^{1/2} \left( x^{\alpha} (x-\eps)^{-b}\right)^{\bar\gamma}\mathrm{d}\lambda_{I}(x)\right)^{1/\bar\gamma}\\
 &\leq \frac{8C|s|}{\pi}\sup_{\eps\in (0,\eps_0]} \eps^{1-\beta}\left(\int_{\left(\frac{2Cs\eps}{\pi}\right)^{1/b}}^{1/2}(x+\eps)^{\alpha\bar\gamma}  x^{-b\bar\gamma}\mathrm{d}\lambda_{I}(x)\right)^{1/\bar\gamma}\\
 &\leq \frac{8C|s|}{\pi}\sup_{\eps\in (0,\eps_0]} \left(
 \eps^{1-\beta}\left(\int_{\left(\frac{2Cs\eps}{\pi}\right)^{1/b}}^{1/2}x^{\bar\gamma(\alpha-b)}\mathrm{d}\lambda_{I}(x)\right)^{1/\bar\gamma} +\eps^{1+\alpha-\beta}\left(\int_{\left(\frac{2Cs\eps}{\pi}\right)^{1/b}}^{1/2}x^{-b\bar\gamma}\mathrm{d}\lambda_{I}(x)\right)^{1/\bar\gamma}\right)\\
 &\lesssim  |s|\sup_{\eps\in (0,\eps_0]} \left( \eps^{1-\beta} \max\left\{\frac{1}{2},\left(\frac{2Cs\eps}{\pi}\right)^{1/b} \right\}^{1/\bar\gamma+\alpha-b}+\eps^{1+\alpha-\beta} \max\left\{\frac{1}{2},\left(\frac{2Cs\eps}{\pi}\right)^{1/b} \right\}^{1/\bar\gamma-b}\right)\\
 &\lesssim |s|\sup_{\eps\in (0,\eps_0]} \left(\max\left\{\eps^{1-\beta},|s|^{1/(\bar\gamma b)+\alpha/b-1}\eps^{1/(\bar\gamma b)+\alpha/b-\beta}\right\}+\max\left\{\eps^{1+\alpha-\beta},|s|^{1/(\bar\gamma b)-1}\eps^{1/(\bar\gamma b)+\alpha-\beta}\right\}  \right) \\
 &\leq \tilde{C}_{\eps_0}|s|(1+|s|^{1/(\bar\gamma b)+\alpha/b-1}+|s|^{1/(\bar\gamma b)-1})
\end{align*}
for some constant $\tilde{C}_{\eps_0}>0$. This follows from the assumption that $1/(\bar\gamma b)+\alpha/b-\beta>0$ and $1/(\bar\gamma b)+\alpha-\beta>0$. 

Finally, combining this with the first step and using symmetry, we have that
\begin{align*}
  \sup_{s\in [0,2\pi]}\sup_{\eps\in (0,\eps_0]}\|R_{\alpha}\osc(e^{i s\chi}, B_{\eps}(\cdot))\|_{\bar\gamma}\cdot \eps^{-\beta} &\leq \tilde{C}_{\eps_0} \sup_{s\in [0,2\pi]}(1+|s|+|s|^{1/(\bar\gamma b)+\alpha/b}+|s|^{1/(\bar\gamma b)}+|s|^{(1+\alpha)/b}) \\ &\leq   C_{\eps_0}
\end{align*}
for some $C_{\eps_0}>0$ which is independent of $\bar\gamma\geq 1$.
\end{proof}

For the following for all $j=0,
 \dots, k-1$, let $ \bar R_{j+1}:[c_j,c_{j+1}] \to{ \mathbb{R}}$ be given by
 $$ \bar{R}_{j+1}=\frac{(R_\alpha\One) \circ  \psi_{j+1}}{R_\alpha\One}$$
and the following lemma is independent of the choice of $\chi$. 

\begin{lem}\label{lem: Rj+1 Bdd}
 $\bar R_{j+1}$ is bounded\footnote{In fact, they are $\alpha$-H\"older continuous. See \Cref{sec:alphaHolder}.} for all $j$. Further, let $0<\eps<\delta$ and $\alpha\in (0,1)$. Then, for all $j$, there is a constant $C$ which is independent of $\eps$ and $\delta$ such that
\begin{equation}\label{eq:SupEst}
    \sup_{x \in [c_j+\delta+\eps, c_{j+1}-\delta-\eps]}\Big( (R_\alpha \One)(x) \sup_{B_{\eps}(x)}|\bar R'_{j+1}|\Big) \leq C\cdot \delta^{\alpha-1}.
\end{equation}
\end{lem}
\begin{proof}

First, we notice that for all $j$
\begin{align}
    \bar  R_{j+1}(x) &=\frac{\psi_{j+1}(x)^\alpha(1-\psi_{j+1}(x))^\alpha}{x^\alpha(1-x)^\alpha}
    \leq \max\left\{\frac{(\psi_{j+1}(x)-0)^\alpha}{x^\alpha}, \frac{(1-\psi_{j+1}(x))^\alpha}{(1-x)^\alpha}\right\}\notag\\
    &\leq \max\left\{\frac{((x-0)\eta_+)^{\alpha}}{x^{\alpha}},\frac{((1-x)\eta_+)^{\alpha}}{(1-x)^{\alpha}}\right\}
    \leq \eta_+^\alpha,\label{eq: bar Rj+1 <gamma}
\end{align}
where the first inequality holds true, because at most one of the arguments in the maximum can be larger than $1$.
Hence, for all $j$, $\bar R_{j+1}$ is bounded. 

We know from \eqref{en: proof strategy 1} in the proof of Lemma \ref{lem: Rj+1 Hoelder} that $\bar{R}_{1}'$ is bounded at $0$ and $\bar{R}_{k-1}'$ is bounded at $1$. We can infer from the representation in \eqref{eq: bar Rj+1'} 
that 
there exist $K_3',K_3>0$ such that
\begin{align}
 |\bar{R}'_{j+1}(x)|\leq \frac{K_3'}{\left(\psi_{j+1}(x) (1-\psi_{j+1}(x))\right)^{1-\alpha}}
 \leq \frac{K_3}{\left((x-c_j) (c_{j+1}-x)\right)^{1-\alpha}},\label{eq: bar R' estim}
\end{align}
for all $j\in \{1,\ldots,k-2\}$.
This can be deduced as follows:~We assume we are in the interval $[\delta_0, 1-\delta_0]$ with $\delta_0$ as in \eqref{en: proof strategy 1} of the proof of Lemma \ref{lem: Rj+1 Hoelder}. Then the subtrahend of \eqref{eq: bar Rj+1'} has to be bounded as it only has a pole at $0$ and $1$. Furthermore, considering the minuend it is easy to notice that the factor $\alpha\psi_{j+1}'(x)(1-2\psi_{j+1}(x))/(x(1-x))^{\alpha}$ has to be bounded on $[\delta_0, 1-\delta_0]$ as well. This leaves the remaining factor as in the middle term of \eqref{eq: bar R' estim}.

In order to verify the second inequality we notice that $\psi_{j+1}(x)\in [\eta_-(x-c_j), \eta_+(x-c_j)]$ which follows from the fact that $\lim_{\eps\to 0}\psi_{j+1}(c_j+\eps)=0$ and from the bound on the derivative. With a similar argumentation, using that $\lim_{\eps\to0}\psi_{j+1}(c_{j+1}-\eps)=1$ we obtain $1-\psi_{j+1}(x)\in [\eta_-(c_{j+1}-x), \eta_+(c_{j+1}-x)]$.

In addition, from the proof of \Cref{lem: Rj+1 Hoelder}
\begin{align*}
 |\bar{R}'_{1}(x)|\leq \frac{K_3}{(c_{1}-x)^{1-\alpha}}
 \qquad \text{ and } \qquad
 \bar{R}'_{k}(x)\leq \frac{K_3}{(x-c_{k-1})^{1-\alpha}}.
\end{align*}
Hence, 
\begin{align}
    \sup_{x \in [c_j+\delta+\eps, c_{j+1}-\delta-\eps]}{\Big(} (R_\alpha \One)(x)\sup_{B_{\eps}(x)}|\bar R'_{j+1}| \Big) &\lesssim 
\begin{cases} \sup\frac{1}{[(x\pm \delta - c_j)(c_{j+1}-x\pm \delta)]^{1-\alpha}}  & j \notin\{ 0, k-1\}\\  \sup\frac{1}{(c_{1}-x\pm \delta)^{1-\alpha}} & j = 0 \\ \sup\frac{1}{(x\pm \delta- c_{k-1})^{1-\alpha}}  & j = k-1  \end{cases} \,\\ &\lesssim \delta^{\alpha-1}\,.\notag
\end{align}
\end{proof}

Now, we are ready to prove the main lemma. 
\begin{lem}\label{lem:UniDoeblinFortetspec}
Let  $0 \leq \alpha < \beta < \min\,\{1/2,{\vartheta}, 1/b\}$ be such that 
$$\kappa:=\frac{\eta^{\alpha}_+}{\eta^{\beta}_-} < 1\,,\qquad 
\text{and}\qquad 
  \max\{|\chi'(x+)|, |\chi'(x-)|\}\lesssim x^{-b}(1-x)^{-b}\,.
$$
Then, for all $1 \leq \gamma < 1/\alpha$ there exist $C, \widetilde{C}>0$ and $\bar \gamma$ with $\gamma < \bar\gamma < 1/\alpha $ such that  for all $s\in \reals$ we have that for all $h\in \V$ and for all $n \in \naturals\,,$
\begin{equation}\label{eq:DFInequality}
   \Vert {\wh \psi}^n_{is} h\Vert_{\alpha,
\beta, \gamma}
\le \widetilde{C} \left(\kappa^n\Vert h\Vert_{\alpha,\beta, \gamma}+ C^n\| h\|_{\bar \gamma}\right)\, .
\end{equation}
\end{lem}
\begin{rem}
In the linear expanding case, i.e., $\eta_+=\eta_- >1$, the condition $\kappa < 1$ reduces to $\beta >\alpha.$ Also, the constant $C$ is independent of $\bar \gamma$.
\end{rem}
\begin{rem}\label{rem:LpInclusion}
Restricting $\bar \gamma$ to $(\gamma, 1/\alpha)$ ensures that $h\in \V$ implies $h\in L^{\bar \gamma}$. To see this, observe that $|R_\alpha h| \lesssim  \One $ which yields that $|h|^{\bar \gamma} \lesssim R_{-\alpha\bar \gamma} \One $, and since $\bar \gamma \alpha < 1\,,$ $R_{-\alpha \bar \gamma} \One$ is integrable. 
\end{rem}

\begin{proof}[Proof of Lemma \ref{lem:UniDoeblinFortetspec}] 
Let $s\in \reals$ and $h \in \V$ be $\reals-$valued. We estimate $|{\wh \psi}_{is} h|_{\alpha,\beta}$:
\begin{align*}
   \osc\big(R_\alpha ({\wh \psi}_{is} h), B_\eps(x)\big)  &= \osc\left( R_\alpha  \sum_{j=0}^{k-1} \left(\frac{e^{is\chi}\cdot h }{|\psi'|} \right)\circ \psi^{-1}_{j+1}\One_{\psi[c_j,c_{j+1}]}\,, B_\eps(x) \right) \\&\leq \sum_{j=0}^{k-1} \osc\left(  R_\alpha \left(\frac{e^{is\chi}\cdot h }{|\psi'|} \right)\circ \psi^{-1}_{j+1} \,, B_\eps(x) \right)
   \\&\leq \sum_{j=0}^{k-1} \osc\left( \frac{ R_\alpha \One\circ \psi_{j+1} }{ R_\alpha \One} \cdot R_\alpha\frac{e^{is\chi } \cdot h }{|\psi'|}, \psi^{-1}_{j+1} B_{\eps}(x) \cap [c_j,c_{j+1}]\right)
   \\&\leq \sum_{j=0}^{k-1} \osc\left( \frac{ R_\alpha \One\circ \psi_{j+1} }{ R_\alpha \One} \cdot R_\alpha\frac{e^{is\chi } \cdot h }{|\psi'|}, B_{\eps/\eta_{-}}(\psi^{-1}_{j+1}x) \cap [c_j,c_{j+1}]\right) 
   \\&= \sum_{j=0}^{k-1} \osc\left( \bar{R}_{j+1}\cdot \frac{e^{is\chi}}{|\psi'|} \cdot R_\alpha h\, , D_{j+1}(x,\eps/\eta_{-})\right), 
\end{align*}
where 
$D_{j+1}(x,\eps):=B_{\eps}(\psi^{-1}_{j+1}x)  \cap [c_j,c_{j+1}]$.
 So, by \cite[Prop.~3.2 (iii)]{BS} there exists $c>0$ such that
\begin{align*}
   \osc\big(R_\alpha ({\wh \psi}_{is} h), B_\eps(x)\big) &\leq \sum_{j=0}^{k-1} \osc\left( R_\alpha h\, , D_{j+1}(x,\eps/\eta_{-})\right) 
   \sup_{D_{j+1}(x,\eps/\eta_{-})} \left|\bar{R}_{j+1}\cdot \frac{e^{is\chi}}{|\psi'|}\right| \\ &\qquad+\sum_{j=0}^{k-1} \osc\left( \left|\bar{R}_{j+1}\cdot \frac{e^{is\chi}}{|\psi'|}\right| ,D_{j+1}(x,\eps/\eta_{-})\right) \inf_{D_{j+1}(x,\eps/\eta_{-})} |R_\alpha h|\,.\\
   &\leq \left(1+c(\eps\eta_{-}^{-1})^{\vartheta}\right) \sum_{j=0}^{k-1} \frac{\osc\left( R_\alpha h\, , B_{\eps/\eta_{-}}(\psi^{-1}_{j+1}x)\right)}{|\psi'|(\psi^{-1}_{j+1}x)}
   \sup_{D_{j+1}(x,\eps/\eta_{-})} \left|\bar{R}_{j+1}\right| \\ &\qquad+\sum_{j=0}^{k-1} \osc\left( \left|\bar{R}_{j+1}\cdot \frac{e^{is\chi}}{|\psi'|}\right|,
    D_{j+1}(x,\eps/\eta_{-})\right) |R_\alpha h(\psi^{-1}_{j+1}x)| \,.
\end{align*}
The last inequality follows from the fact that $\psi^{-1}$ is $C^1$ and its derivative is uniformly ${\vartheta}-$H\"older. 

Hence, using the upper bound \eqref{eq: bar Rj+1 <gamma}, and then using the definition of the transfer operator $\wh\psi$, we have 
\begin{align}\label{eq:OscOp0}
   \osc\big(R_\alpha ({\wh \psi}_{is} h), B_\eps(x)\big) &\leq \left(1+c(\eps\eta_{-}^{-1})^{\vartheta}\right){\eta^{\alpha}_{+}}\wh \psi\big(\osc(R_\alpha h, B_{\eps/\eta_{-}}(\,\cdot\,))\big)(x) \\ &\qquad+\sum_{j=0}^{k-1} |R_\alpha h(\psi^{-1}_{j+1}x)| \osc\left( \left|\bar{R}_{j+1}\cdot \frac{e^{is\chi}}{|\psi'|}\right|, 
    D_{j+1}(x,\eps/\eta_{-})\right)\,.\notag
\end{align}
Taking the integral over the first term in \eqref{eq:OscOp0} and multiplying by $\eps^{-\beta}$ we obtain 
\begin{align}
 \MoveEqLeft\eps^{-\beta}\int \left(1+c(\eps\eta_{-}^{-1})^{\vartheta}\right){\eta}^{ \alpha}_{+}\wh \psi\big(\osc(R_\alpha h, B_{\eps/\eta_{-}}(\,\cdot\,))\big)(x)\,\mathrm{d}\lambda_I(x)\notag\\
 &\leq \eps^{-\beta}\left(1+c(\eps\eta_{-}^{-1})^{\vartheta}\right){\eta}^{\alpha}_{+}\int \wh \psi\big(\osc(R_\alpha h, B_{\eps /\eta_{-}}(\,\cdot\,))\big)(x)\,\mathrm{d}\lambda_I(x)\notag\\
 &= \eps^{-\beta}\left(1+c(\eps\eta_{-}^{-1})^{\vartheta}\right){\eta}^{\alpha}_{+} 
 \int \osc(R_\alpha h, B_{\eps/\eta_{-}}(\,\cdot\,))(x)\,\mathrm{d}\lambda_I(x)\notag\\
 &\leq\left(1+c(\eps\eta_{-}^{-1})^{\vartheta}\right)\eta^{\alpha}_{+}\eta_-^{-\beta}|h|_{\alpha,\beta}\notag\\
 &\leq \left(1+c(\eps_0\eta_{-}^{-1})^{\vartheta}\right)\kappa\|h\|_{\alpha,\beta,\gamma},\label{eq:OscOp0a}
\end{align}
for all $\gamma \geq 1$. 
Next, we analyze the second term in \eqref{eq:OscOp0}. Again, by \cite[Prop.~3.2 (iii)]{BS} we have
\begin{align}
     &\osc\left( \left|\bar{R}_{j+1}\cdot \frac{e^{is\chi}}{|\psi'|}\right|, D_{j+1}(x,\eps/\eta_{-})\right)\nonumber \\
     &\quad\leq \osc\left(  \frac{1}{|\psi'|}\,, B_{\eps/\eta_{-}}(\psi^{-1}_{j+1}x)\right)\left(\esssup_{B_{\eps/\eta_{-}}(\psi^{-1}_{j+1}x)} |\Re\bar{R}_{j+1}e^{is\chi}|+ \esssup_{B_{\eps/\eta_{-}}(\psi^{-1}_{j+1}x)} |\Im\bar{R}_{j+1}e^{is\chi}|\right) \nonumber\\
     &\qquad\quad+\osc\left(  \bar{R}_{j+1}e^{is\chi}\,, D_{j+1}(x,\eps/\eta_{-})\right)\inf_{D_{j+1}(x,\eps/\eta_{-})} \frac{1}{|\psi'|} \nonumber\\
     &\quad\leq c(\eps\eta_{-}^{-1})^{\vartheta}\eta^{\alpha}_{+} \frac{1}{|\psi'|(\psi^{-1}_{j+1}x)}+(1+c(\eps\eta_{-}^{-1})^{\vartheta})\frac{\osc\left(\bar{R}_{j+1}e^{is\chi}\,,D_{j+1}(x,\eps/\eta_{-})\right)}{|\psi'|(\psi^{-1}_{j+1}x)}.\label{eq: osc bar Rj e psi}
\end{align}
Note that 
\begin{align}\label{eq:FirstTerm}
    \eps^{-\beta}c(\eps\eta_{-}^{-1})^\eta\eta^{\alpha}_{+} \int \sum_{j=0}^{k-1}  \frac{|R_\alpha h(\psi^{-1}_{j+1}x)|}{|\psi'|(\psi^{-1}_{j+1}x)}\, \mathrm{d}\lambda_I(x)  \nonumber &=\eps^{-\beta} c(\eps\eta_{-}^{-1})^{\vartheta}\eta^{\alpha}_{+} \int\wh\psi(|R_\alpha h|)\, \mathrm{d}\lambda_I(x)  \\ 
    &=\eps^{-\beta} c(\eps\eta_{-}^{-1})^{\vartheta}\eta^{\alpha}_{+} \int  |R_\alpha h|\, \mathrm{d}\lambda_I(x) \nonumber \\ &\leq K_1 \eps^{{\vartheta} - \beta}\|R_\alpha\One\|_{\gamma_1}\|h\|_{\bar \gamma}
\end{align}
where $\gamma^{-1}_1+\bar\gamma^{-1}=1$, $K_1:= c\eta_{-}^{-{\vartheta}}\eta^{\alpha}_{+}\|R_\alpha\One\|_{\bar\gamma}$ and 
$\beta < {\vartheta}$. So, the contribution from the first summand of  \eqref{eq: osc bar Rj e psi} to \eqref{eq:OscOp0} is under control.

To estimate the contribution from second summand of \eqref{eq: osc bar Rj e psi} to \eqref{eq:OscOp0} we note that for all $j$ and for all $A \subset [c_j,c_{j+1}]$, we have
$$\osc\left(  \bar{R}_{j+1}e^{is\chi}\,,A \right)= \osc\left( \sum_{j=0}^{k-1}\bar{R}_{j+1}e^{is\chi}\One_{[c_j,c_{j+1})}\,,A\right), $$
and therefore, we can bound this contribution by
\begin{align}
    &(1+c(\eps\eta_{-}^{-1})^{\vartheta})\sum_{j=0}^{k-1} \frac{|R_\alpha h(\psi^{-1}_{j+1}x)|}{|\psi'|(\psi^{-1}_{j+1}x)}\osc\left( F\, ,\, B_{\eps/\eta_{-}}(\psi^{-1}_{j+1}x) \right)\nonumber \\
    &\quad\quad=(1+c(\eps\eta_{-}^{-1})^{\vartheta}) \wh \psi \left(|R_\alpha h| \osc\left( F\, ,\, B_{\eps/\eta_{-}}(\,\cdot\,)\right)\right),
\end{align}
where
$$F(x) = e^{is\chi(x)}  \sum_{j=0}^{k-1} \bar R_{j+1}(x) \One_{[c_j,c_{j+1})}(x)=e^{is\chi(x)}  \sum_{j=0}^{k-1} \frac{R_{\alpha } \One \circ \psi_{j+1}(x)}{R_{\alpha } \One (x)} \One_{[c_j,c_{j+1})}(x).$$
This is bounded by
\begin{align}\label{eq:SecondTerminSecondary}
   & 
   (1+c(\eps\eta_{-}^{-1})^{\vartheta}) \int \wh \psi \left(|R_\alpha h| \osc\left( F\, ,\, B_{\eps/\eta_{-}}(\,\cdot\,)\right)\right)(x)\, \mathrm{d}\lambda_I(x)  \nonumber \\ &\qquad= 
   (1+c(\eps\eta_{-}^{-1})^{\vartheta})\int|R_\alpha h|(x) \osc\left( F\, ,\, B_{\eps/\eta_{-}}(x)\right)\, \mathrm{d}\lambda_I(x) \nonumber\\
   &\qquad= 
   (1+c(\eps\eta_{-}^{-1})^{\vartheta})\int|h (x)| \cdot \left(R_\alpha\osc\Big( F\, ,\, B_{\eps/\eta_{-}}(x)\Big)\right)\, \mathrm{d}\lambda_I(x)\,.
\end{align}
To estimate the integral we split it as follows.
\begin{align*}
    &\int|h (x)| \cdot \left(R_\alpha\osc\Big( F\, ,\, B_{\eps/\eta_{-}}(x)\Big)\right)\, \mathrm{d}\lambda_I(x) \\ 
    &\quad=  \left( \sum_{ j=1}^{k-1} \int_{c_j-\eps^\iota-\eps}^{ c_j+\eps^\iota+\eps} + \sum_{ j=1}^{k} \int_{ c_{j-1}+\eps^\iota+\eps}^{c_{j}-\eps-\eps^\iota} + \int_{0}^{ \eps+\eps^\iota } + \int_{ 1-\eps-\eps^\iota}^{1}\right)  |h (x)| \cdot \left(R_\alpha\osc\Big( F\, ,\, B_{\eps/\eta_{-}}(x)\Big)\right)\, \mathrm{d}\lambda_I(x) 
\end{align*}
 where we choose for $\iota$ any number fulfilling 
 \begin{equation}
  \frac{\beta}{1-\alpha}<\iota< \frac{1-\beta}{1-\alpha}. \label{eq: iota cond}
 \end{equation}
Because $\beta<1/2$  such a choice is possible.
Note that for $j=1,\dots k-1\,,$ $x \in { [c_j-\eps^\iota-\eps,c_j+\eps^\iota+\eps]}\,,$ 
\begin{align*}
    \osc\left( F\, ,\, B_{\eps/\eta_{-}}(x)\right) \leq 2 (\sup \bar R_{j}+\sup \bar R_{j+1}) \leq 4K
\end{align*}
and for $x \in { [0,\eps+\eps^\iota) \cup (1-\eps-\eps^\iota, 1]}\,,$
\begin{align*}
    \osc\left( F\, ,\, B_{\eps/\eta_{-}}(x)\right) \leq 2 (\sup \bar R_{0}+\sup \bar R_{ k})\, \leq 4K
\end{align*}
where $K:= \sup_j \sup R_{j+1} < \infty$. So,
\begin{align}
     &\sum_{j=0}^{k} \int_{ (c_j-\eps^\iota-\eps) \vee 0}^{  (c_j+\eps^\iota+\eps) \wedge 1 }  |h (x)| \cdot \left(R_\alpha\osc\Big( F\, ,\, B_{\eps/\eta_{-}}(x)\Big)\right)\, \mathrm{d}\lambda_I(x) \notag\\ &\quad\leq \| h\|_{\bar\gamma}\sum_{j=0}^{k}  \left(\int_{ (c_j-\eps^\iota-\eps) \vee 0}^{(c_j+\eps^\iota+\eps) \wedge 1 } \left(R_\alpha\osc\Big( F\, ,\, B_{\eps/\eta_{-}}(x)\Big)\right)^{\gamma_1}\, \mathrm{d}\lambda_I(x)\right)^{1/\gamma_1}\notag \\
     &\quad\leq K_{\alpha}\eps^{\iota/\gamma_1}\|h\|_{\bar\gamma}\label{eq:NearDiscEst}
\end{align}
where $\gamma^{-1}_1+\bar\gamma^{-1}=1$
and $K_{\alpha}= 4^{\iota/\gamma_1}2^{-2\alpha}K$.  Here, we choose $\bar\gamma$ such that 
\begin{equation}
 \max\left\{\gamma,\, \frac{\iota}{\iota-\beta},\, \frac{1}{1-b\beta+\alpha}\,,\frac{1}{1-b(\beta-\alpha)} \right\}<\bar\gamma<\frac{1}{\alpha}.\label{eq: gamma1 cond}
\end{equation}
We will see later in the proof why this restrictions on $\bar{\gamma}$ are needed.

Now, we show that such a choice is possible. Since we were assuming that $\iota>\beta/(1-\alpha)$, we have $\iota/(\iota-\beta)<1/\alpha$. We note that when $b \leq 1$, $1-b\beta+\alpha \geq 1-b(\beta-\alpha)$, and it is enough to see whether $\alpha < 1-b(\beta-\alpha)$. In fact, this is true because $b(\beta-\alpha)< \beta -\alpha < 1 - \alpha$. On the contrary, when $b> 1$, we have $1-b\beta+\alpha < 1-b(\beta-\alpha)$, and it is enough to see whether $\alpha < 1-b\beta+\alpha$. This is true because $\beta < 1/b$.

To estimate the remaining terms we note, using \eqref{eq: bar Rj+1 <gamma}  and \cite[Prop.~3.2(iii)]{BS}, that for all $j=0,\dots, k-1\,,$ for all $x \in  [c_j+\eps^\iota+\eps, c_{j+1}-\eps^\iota-\eps]$,
\begin{align*}
    \MoveEqLeft\osc\Big( F\, ,\, B_{\eps/\eta_{-}}(x)\Big)\\
    &=  \osc\Big( e^{is\chi}\bar R_{j+1}\, ,\, B_{\eps/\eta_{-}}(x)\Big)\\
    &\leq  \sup_{B_{\eps}(x)} (\Re |e^{is\chi}|+\Im |e^{is\chi}|)  \osc\Big(\bar R_{j+1}\, ,\, B_{\eps/\eta_{-}}(x)\Big)
    + \osc\Big( e^{is\chi}\, ,\, B_{\eps/\eta_{-}}(x)\Big) \inf_{B_{\eps}(x)} \bar R_{j+1}\\
        &\leq  2  \sup_{B_{\eps}(x)}|\bar R'_{j+1}|\,  \frac{\eps}{\eta_-}
    + \osc\Big( e^{is\chi}\, ,\, B_{\eps/\eta_{-}}(x)\Big)\, \eta_+^{\alpha},
\end{align*}
and thus,

\begin{align*}
   &\sum_{j=0}^{k-1}\int_{c_j+\eps^\iota+\eps}^{c_{j+1}-\eps^\iota-\eps} |h (x)| \cdot \left(R_\alpha\osc\Big( F\, ,\, B_{\eps/\eta_{-}}(x)\Big)\right)\, \mathrm{d}\lambda_I(x) \nonumber \\
   & \quad\leq  \frac{2\eps}{\eta_-}\left\||h| \sum_{j=0}^{k-1}\One_{[c_j+\eps^\iota+\eps,c_{j+1}-\eps^\iota-\eps]}\right\|_{1} \sup_{x \in { [c_j+\eps^\iota+\eps, c_{j+1}-\eps^\iota-\eps]}}R_{\alpha}\sup_{B_{\eps}(x)}|\bar R'_{j+1}|\\
   &\qquad\quad+\eta_+^{\alpha}\left\| \sum_{j=0}^{k-1}\One_{[c_j+\eps^\iota+\eps, c_{j+1}-\eps-\eps^\iota]}\cdot |h| \left(R_\alpha \osc\Big( e^{is\chi}\, ,\, B_{\eps/\eta_{-}}(\cdot)\Big)\right)\right\|_{1}\\
   &\quad\leq  \frac{2\eps}{\eta_-}\left\|h\right\|_{1} \sup_{x \in { [c_j+\eps^\iota+\eps, c_{j+1}-\eps^\iota-\eps]}}R_{\alpha}\sup_{B_{\eps}(x)}|\bar R'_{j+1}|+\eta_+^{\alpha}\|h\|_{\bar\gamma}\left\|R_\alpha \osc\Big( e^{is\chi}\, ,\, B_{\eps/\eta_{-}}(\cdot)\Big)\right\|_{\bar\gamma}\,.
\end{align*}

Now, in order to estimate the first summand taking the maximum over $j$ of the supremum in \eqref{eq:SupEst} above  with $\delta=\eps^\iota$ yields that the  outer supremum above is bounded by  $C\eps^{\iota(\alpha-1)}$ for some constant $C>0$. For the second summand, from \eqref{eq: gamma1 cond}, we have that ${\bar\gamma}^{-1} <  1 - b\beta+\alpha$ which implies that $b\beta < 1- {\bar\gamma}^{-1}+\alpha = \bar\gamma^{-1}+\alpha $, and hence, when $b>1$, we have the condition \eqref{eq:GammaAlphaBetaB}. Also from \eqref{eq: gamma1 cond}, ${\bar\gamma}^{-1} <  1 - b(\beta-\alpha)$ which implies that $b\beta <  \bar\gamma^{-1}+b\alpha $, and hence, when $b \leq 1$, we have \eqref{eq:GammaAlphaBetaB}. Therefore, we can apply \Cref{lem:DF1} with $\alpha$, $\beta$, $b$, $\bar\gamma$, $\eps/\eta_{-}$ to conclude
$$\left\|R_\alpha \osc\Big( e^{is\chi}\, ,\, B_{\eps/\eta_{-}}(\cdot)\Big)\right\|_{\bar \gamma} \leq C_{\eps_0} \eps^\beta \eta^{-\beta}_{-}\,$$
where $C_{\eps_0}$ is independent of $\bar\gamma$.
Therefore, for all $s \neq 0$, 
\begin{align}\label{eq:AwayDiscEst}
    \MoveEqLeft \sum_{j=0}^{k-1} \int_{ c_j+\eps^\iota+\eps}^{  c_{j+1}-\eps^\iota-\eps} |h (x)| \cdot \left({  R_{\alpha}}\osc\Big( F\, ,\, B_{\eps/\eta_{-}}(x)\Big)\right)\, \mathrm{d}\lambda_I(x)\notag\\
    &\leq \bar  C_{\eps_0}
    \eps^{\min\{1-\iota(1-\alpha), \beta\}}\|h\|_{\bar\gamma}\,.
\end{align}
Finally, combining \eqref{eq:NearDiscEst} and \eqref{eq:AwayDiscEst}, we estimate \eqref{eq:SecondTerminSecondary}  multiplied by $\eps^{-\beta}$ by
\begin{align}\label{eq:SecondTerm}
 \MoveEqLeft\eps^{-\beta}(1+c(\eps\eta_{-}^{-1})^{\vartheta}) \int \wh \psi \left(|R_\alpha h| \osc\left( F\, ,\, B_{\eps/\eta_{-}}(\,\cdot\,)\right)\right)(x)\, \mathrm{d}\lambda_I(x)  \nonumber\\ 
 &\leq  \eps^{  (\iota /\bar\gamma\wedge (1-\iota (1-\alpha))\wedge \beta)-\beta}C_{\eps_0} 
 \|h\|_{\bar\gamma}  \leq C_{\eps_0} \|h\|_{\bar\gamma}. 
\end{align}
To justify the last inequality, we analyse the exponent of $\eps$. 
By \eqref{eq: gamma1 cond}
and the relation $\gamma_1^{-1}+\bar{\gamma}^{-1}=1$ we have
$\iota/\gamma_1>\iota(1-{\bar\gamma}^{-1})>\iota(1-(\iota-\beta)/\iota))=\beta$.
Furthermore, the second inequality of \eqref{eq: iota cond} implies that $1-\iota (1-\alpha)>\beta$. 

Combining \eqref{eq:OscOp0}, \eqref{eq:OscOp0a}, \eqref{eq:FirstTerm} and \eqref{eq:SecondTerm}, we have 
\begin{align*}
| {\wh \psi}_{is} h|_{\alpha,
\beta}&=\sup_{\eps\in (0,\eps_0)}\int \frac{\osc\big(R_\alpha ({\wh \psi}_{is} h), B_\eps(x)\big)}{\eps^{\beta}}\,\mathrm{d}\lambda_I(x)\\
&\leq \left(1+c(\eps_0\eta_{-}^{-1})^{\vartheta}\right)\kappa\|h\|_{\alpha,\beta,\gamma}+C_{\eps_0}\|h\|_{\bar\gamma},
\end{align*}
for all $\gamma\geq 1$. 
Therefore, for all $\bar\gamma$ chosen appropriately
\begin{align*}
\| {\wh \psi}_{is} h\|_{\alpha,
\beta, \gamma}  &= | {\wh \psi}_{is} h|_{\alpha,\beta} + \|{\wh \psi}_{is} h\|_{\gamma}\\
&\leq \left(1+c(\eps_0\eta_{-}^{-1})^{\vartheta}\right)\kappa\|h\|_{\alpha,\beta,\gamma}+C_{\eps_0}{\|h\|_{\bar\gamma}}+ C_\gamma \|h\|_\gamma\,  \\
&\leq \bar\kappa\|h\|_{\alpha,\beta,\gamma}\bar C \|h\|_{\bar\gamma} 
\end{align*}
where $\bar\kappa = \left(1+c(\eps_0\eta_{-}^{-1})^\eta\right)\kappa < 1$ for sufficiently small $\eps_0$, and  $\bar C=C_{\eps_0}+C_{\gamma}$, where $C_{\gamma}$ is given in Lemma \ref{lem:LpEst}.

Iterating, we obtain the following DFLY inequality:~for all $h\in \V$
\begin{align*}
    \sup_{s}\Vert {\wh \psi}^n_{is} h\Vert_{\alpha,\beta,\gamma}
    &\leq \kappa\Vert {\wh \psi}^{n-1}_{is}h\Vert_{\alpha,\beta,\gamma}+\bar C\| {\wh \psi}^{n-1}_{is}h\|_{\bar\gamma}\,\\
    &\leq \kappa^2\Vert {\wh \psi}^{n-2}_{is}h\Vert_{\alpha,\beta,\gamma}+ \kappa \bar C\|{\wh \psi}^{n-2}_{is}h\|_{\bar\gamma}+\bar C \bar C^{n-1}\| h\|_{\bar\gamma}\,\\
    &\leq \kappa^n\Vert h\Vert_{\alpha,\beta,\gamma}+ \bar C \|h\|_{\bar\gamma} \sum_{j=0}^{n-1} \kappa^j \bar C^{n-1-j}\,\\
    &\leq \kappa^n\Vert h\Vert_{\alpha,\beta,\gamma}+ C \bar C^{n+1} \|h\|_{\bar\gamma}\,
\end{align*}
for some $C>0$.

In the proof above, we assumed that $h$ is $\reals-$valued. When $h=h_1+ih_2$ where $h_j,\, j=1,2$ are $\reals-$valued, using linearity of the operator $$\|{\wh \psi}_{is}^n h\|_{\alpha,\beta,\gamma} \leq \|{\wh \psi^n}_{is} h_1\|_{\alpha,\beta,\gamma}+ \|{\wh \psi^n}_{is} h_2\|_{\alpha,\beta,\gamma}\,,$$ and also, $\| h_j\|_{\alpha,\beta,\gamma} \leq \| h\|_{\alpha,\beta,\gamma}$ and $\|h_j\|_{\bar\gamma} \leq \|h\|_{\bar\gamma}$ for all $j=1,2$. So, applying DLFY inequality proven above in the $\reals-$valued case to $h_1$ and $h_2$, we conclude that DFLY in the general case of $h$ holds up to a constant multiple.
\end{proof}

\section{Proofs of the Main Theorems}\label{sec:Proofs}
Finally, we give the proofs of our main theorems. We start with the theorems from \Cref{sec:IntThms}.

\subsection{Proofs of limit theorems for expanding interval maps}\label{sec:LimitThms}

\begin{proof}[Proof of Theorem \ref{thm:IntervalCLTsimple}]
From \eqref{eq: CLT main cond} we obtain that there exist $\alpha, \beta$ fulfilling 
\begin{align}
 a&< \alpha < \beta \cdot \min\left\{1,\frac{\log\eta_-}{\log \eta_+}\right\}< \min\left\{\vartheta,\frac{1}{b},\frac{1}{2}\right\}\cdot \min\left\{1,\frac{\log\eta_-}{\log \eta_+}\right\}\,.\label{eq: a alpha beta theta}
\end{align} 
Furthermore, since $\alpha>a$, the inequality $\beta<1/b$ which we can deduce immediately from \eqref{eq: a alpha beta theta} that $\frac{1}{b}<\frac{1}{b-a}$.
So, by Lemma \ref{lem:abtoalphabeta}, we obtain $|\chi|_{\alpha,\beta}<\infty$ and also  
$\chi \in \VVc{\alpha}{\beta}{2} \hookrightarrow L^2$. 
Furthermore, from the second inequality of \eqref{eq: a alpha beta theta} we obtain $\eta_+^{\alpha}/\eta_-^{\beta}< 1$.

Since $\psi$ is a piecewise $C^2$ uniformly expanding and a covering map of the interval, it has a unique absolutely continuous invariant mixing probability (acip) with a bounded invariant density; see \cite{L1}.
Let's call this acip $\pi$. Then $L^2 \hookrightarrow L^2(\pi)$ because  
$$ \int |h|^2\, \mathrm{d}\pi = \int |h|^2  \frac{\mathrm{d}\pi}{\mathrm{d}\lambda_I} \, \mathrm{d}\lambda_I \leq \left\|\frac{\mathrm{d}\pi}{\mathrm{d}\lambda_I}\right\|_{\infty} \int |h|^2\, \mathrm{d}\lambda_I\,.$$ 

We claim that $\wh\psi$ has a spectral gap in $\V$ with $\gamma=2$.  In \Cref{sec:CtsIncl}, we show that $\VVc{\alpha}{\beta}{2}$ is continuously embedded in $L^{\bar  \gamma}$ where $\bar\gamma \in (2, 1/\alpha)$ and that the unit ball of $\VVc{\alpha}{\beta}{2}$ is relatively compact in $L^{\bar\gamma}$.  A suitable $\bar\gamma$ exists by the condition $\alpha<1/2$ from \eqref{eq: a alpha beta theta}.
So, the claim follows from \cite[Lemma B.15]{DKL} due to the DFLY inequality \eqref{eq:DFInequality} with $s=0$ 
and \Cref{rem:ToApplyHN}.
 
Now, the CLT (in the stationary case) follows directly from \Cref{prop:AbsCLT} applied to $\chi - \pi(\chi)$. That is, from \eqref{eq:AbsCLT} we have 
\begin{equation*}\label{eq:StationaryCLT} 
     \Prob_{\pi}\left(\frac{S_{n}(\chi) - n\, \pi(\chi)}{\sigma\sqrt{n}} \leq x\right) - \fN(x) = o(1),\quad\text{as}\quad n\to \infty
\end{equation*}
with $\sigma^2>0$ because $\chi$ is not a coboundary.
\end{proof}

Next, we will continue with the proof of Theorem \ref{thm:IntervalEdgeExp} as the proof of Theorem \ref{thm:IntervalLLT} will need similar methods to those of Theorem \ref{thm:IntervalEdgeExp}.
 
\begin{proof}[Proof of Theorem \ref{thm:IntervalEdgeExp}]
\eqref{eq: Edge main cond} implies that there exist $\alpha, \beta$ such that $\alpha>a$ and
\begin{align*}
 3\bar{\alpha}:=3\min\{2\alpha, \max\{\alpha,\alpha+b-2\}\}<\beta\cdot\min\left\{1,\frac{\log\eta_+}{\log \eta_-}\right\}< \min\left\{\vartheta, \frac{1}{b},\frac{1}{2}\right\}\cdot\min\left\{1,\frac{\log\eta_+}{\log \eta_-}\right\}\,.
\end{align*}
Since either $b<a+1$ or $1/b<(1+\alpha-a)/(b-a)$,  
we obtain by Lemma \ref{lem:abtoalphabeta} that $|\chi|_{\alpha,\beta}<\infty$ and additionally we obtain by the last inequality that
\begin{enumerate}[label=$(\alph*)$]
     \item\label{ee: iv}  $0< 3\bar\alpha < \beta < \min \{{\vartheta},1/b, 1/2\},$
  \item\label{ee: v}  $\eta_+^{3\bar\alpha}<\eta_-^{\beta}$.
 \end{enumerate} 
 Hence,
under our assumptions, we have the following:
\begin{enumerate}[label=(\arabic*), leftmargin=*, itemsep=5pt]
        \item \label{it:1}  The second inequality in \eqref{eq: bound ab} and $|\chi|_{\alpha, \beta} < \infty$ imply that $|e^{is\chi}|_{0,\beta}<\infty$ for all $s>0$ (see \Cref{rem:RegObs1}). So, due to {\Cref{cor:RegOpr} (1)}, we have $\wh\psi_{is} \in \cL(\VVc{\tilde{\alpha}}{\beta}{\tilde{\gamma}})$ for all $0<\tilde{\alpha}< \beta$ and $\tilde{\gamma} \geq 1$.
        \item \label{it:2} Since $|\chi|_{\alpha,\beta}<\infty$, from  \Cref{rem:RegObs3}, for all $\alpha^*>0\,$ close to $0$, $$\lim_{s\to 0}|1-e^{is\chi}|_{\alpha^*,\beta}=0\,.$$ Along with {\Cref{cor:RegOpr} (2)}, this yields that for all $0 \leq \alpha_0 < \beta$, $\gamma_0 \geq 1$,  $$s \mapsto \wh\psi_{is} \in \cL(\VVc{\alpha_0}{\beta}{\gamma_0}, \VVc{\alpha_1}{\beta}{\gamma_1})$$ is continuous for $\alpha_1=\alpha^*+\alpha_0$ and $1 \leq \gamma_1\leq \gamma_0.$
        \item  \label{it:3} From the second inequality in \eqref{eq: bound ab} and $|\chi|_{\alpha, \beta} < \infty$\,, for all ${\alpha^{**}>\min\{2\alpha, \max\{\alpha+b-2,\alpha\}\}}\,,$
        $$\lim_{s\to 0}\left|\frac{e^{is\chi}-1-is\chi}{s}\right|_{\alpha^{**}, \beta}=0\,$$ due to  \Cref{rem:RegObs4}. Then, we have that for all $0 \leq \alpha_1 <\beta$, and $\gamma_1 \geq 1\,,$ $$s \mapsto \wh \psi_{is} \in \cL(\VVc{\alpha_1}{\beta}{ \gamma_1},\VVc{\alpha_2}{\beta}{\gamma_2} )$$ is continuously differentiable, for all ${\alpha_2 = \alpha^*+\max\{\alpha^*,\alpha^{**}\}+\alpha_1}$ and $1 \leq \gamma_2 \leq (\gamma^{-1}_1+\gamma^{-1})^{-1}$ due to {\Cref{cor:RegOpr} (2) and (3)}.
    \end{enumerate} 

Next, we define the following chain of spaces in order to invoke \Cref{prop:AbstractThm} with $r=1$:
 \begin{align*}
 \VVc{\alpha_0}{\beta}{\gamma_0
  } &\hookrightarrow \VVc{\alpha_1}{\beta}{\gamma_1}  \hookrightarrow \VVc{\alpha_2}{\beta}{\gamma_2} \hookrightarrow  \VVc{\alpha_3}{\beta}{\gamma_3} \hookrightarrow \VVc{\alpha_4}{\beta}{\gamma_4} \hookrightarrow\VVc{\alpha_5}{\beta}{\gamma_5} \hookrightarrow\VVc{\alpha_6}{\beta}{\gamma_6}\hookrightarrow\VVc{\alpha_7}{\beta}{\gamma_7},
  \end{align*}
where $\alpha_0=0$,
$\alpha_{2j} - \alpha_{2j-1}\geq \min\{2\alpha, \max\{\alpha + b - 2, \alpha\}\}$, for $j=1,2,3$, $\alpha_{2j+1} > \alpha_{2j}$ for $j=0,1,2,3$, and $\alpha_7<\beta$. 
By \ref{ee: iv} such a choice is possible. 
Furthermore, we assume that the $\gamma_j$s are chosen such that $\gamma_0=
M \gg 1$ sufficiently large, 
$\gamma_{2j+1} = \gamma_{2j}$ and $\gamma_{2j}< (\alpha^{-1}+\gamma_{2j-1}^{-1})^{-1}$.

Now, to prove the theorem, we verify the conditions in \Cref{prop:AbstractThm} for the above sequence of Banach spaces. We notice that if for some observable $\varphi$ it holds that $|\varphi|_{\alpha,\beta}<\infty$, then $\|\varphi\|_{\alpha,\beta,\gamma}<\infty$
as long as $\gamma< 1/\alpha$.
We next verify that it is possible to construct valid spaces with the above choice of parameters. 
First, we notice that by \ref{ee: iv} it is  possible to construct $\alpha_0\leq \ldots\leq \alpha_7$ with the above properties that $\alpha_7<\beta$ and thus $\alpha_j<\beta$ for all $j$.
Furthermore, by \ref{ee: iv} we have $\alpha <1/3$. 
Thus, it is possible that $1 \leq  \gamma_{2j} \leq  (\gamma^{-1} + \gamma^{-1}_{2j-1})$ holds together with  
$1/\gamma_j>\alpha_j$.
Moreover, under \ref{ee: v} we have that $\eta_+^{\alpha_j}/\eta_-^{\beta}<1$ holds for all $j$.

With that it becomes immediate from applying the conditions of this theorem on the parameters in the Banach spaces and from the calculations in \ref{it:1}--\ref{it:3} applied to all indices $j$ that conditions \ref{cond1}-\ref{cond3} of \Cref{prop:AbstractThm} are satisfied.

For each $j$, we apply \Cref{lem:UniDoeblinFortetspec} with $\gamma=\gamma_j$ and we choose $\bar\gamma = \bar \gamma_j$  as in the proof of the lemma. In \Cref{sec:CtsIncl}, we show that $\VVc{\alpha_j}{\beta}{\gamma_j}$ is continuously embedded in $L^{\bar\gamma_j}$ and that the unit ball of $\VVc{\alpha_j}{\beta}{\gamma_j}$ is relatively compact in $L^{\bar\gamma_j}$. Also, we recall from \Cref{lem:LpEst} that for all $h \in L^{\bar \gamma_j}$,
$\|\wh \psi_{is} (h)\|_{\bar\gamma_j} 
\leq  C_{\bar\gamma_j} \|h\|_{\bar\gamma_j}$ where $ C_{\bar\gamma_j}>1$. Therefore, $\|\wh\psi^n_{is} (h)\|_{\bar\gamma_j} \leq C_{\bar\gamma_j} \|\wh\psi^{n-1}_{is} (h)\| \leq C_{\bar\gamma_j}^n \|h\|_{\bar\gamma_j}$ which gives us  $\|\wh\psi^n_{is}\|_{L^{\bar \gamma_j} \to L^{\bar \gamma_j}} \leq    C_{\bar{\gamma}_j}^n.$
Choose $\kappa = \max\limits_{0 \leq j \leq 7} \eta^{\alpha_j}_+ \eta^{-\beta}_{-} <1$. Also, by our previous constructions, we have that $\gamma_j < 1/\alpha_j$ for all $j$. So, due to Lemma \ref{lem:UniDoeblinFortetspec}, we have the DFLY inequality:~for all $h\in \VVc{\alpha_j}{\beta}{\gamma_j}$
$$\Vert {\wh \psi^n}_{is} h\Vert_{\alpha_j,\beta,\gamma_j}
\le \widetilde C \left(\kappa^n\Vert h\Vert_{\alpha_j,\beta,\gamma_j}+C^n\| h\|_{\bar\gamma_j}\right)\,
$$
for some $\gamma_j< \bar\gamma_j < 1/\alpha_j$ and $C$ uniform in $j$ and $s$. Therefore, we have the first conclusion, equation (8), of \cite[Theorem 1]{KL} uniformly over all spaces. That is, there exist  $v$ and $w$ such that
$$\sup_{z \in D_{\kappa}}\|(z\Id - \wh\psi_{is})^{-1} h\|_{\VVc{\alpha_j}{\beta}{\gamma_j} \to \VVc{\alpha_j}{\beta}{\gamma_j}} \leq  v\|h\|_{\alpha_j,\beta, \gamma_j} + w\|h\|_{\bar \gamma_j}$$
for all space pairs $\VVc{\alpha_j}{\beta}{\gamma_j} \hookrightarrow L^{\bar \gamma_j}$ and 
$s \in \reals$. This gives \ref{cond4} of \Cref{prop:AbstractThm}. 

The conditions \ref{cond5}--\ref{cond7} of \Cref{prop:AbstractThm} are equivalent to Assumption (B) in \cite[Section I.1.2]{FP} for a single dynamical system, i.e., when Assumptions (0) and (A)(1) in \cite[Section I.1.2]{FP} are trivially true. Moreover, as discussed in \cite{FP}, \cite[Lemma 4.5]{FP} implies Assumption (B).  Therefore, we verify the conditions (with a slight modification) in \cite[Lemma 4.5]{FP} to establish \ref{cond5}--\ref{cond7}:\vspace{-5pt}
\begin{itemize}[leftmargin=*]
    \item We have assumed that $\chi$ is non-arithmetic.
    \item Due to \Cref{rem:ToApplyHN} and the DFLY inequality \eqref{eq:DFInequality}, we can apply \cite[Lemma B.15]{DKL} to conclude that for all $s$ the essential spectral radius of $\wh\psi_{is}$  on $\VVc{\alpha_j}{\beta}{\gamma_j}$ is at most $\kappa$. This is precisely the conclusion of \cite[Proposition 4.3]{FP}.
    \item We know that $\VVc{\alpha_j}{\beta}{\gamma_j}\hookrightarrow L^1$ for all $j$, and that $\|\wh\psi_{is} h\|_1 \leq \|\wh\psi h \|_1 \leq \|h\|_1$ for all $h\in L^1$. So, the spectral radius of $\wh\psi_{is}$ on $L^1\,,$ and hence, on $\VVc{\alpha_j}{\beta}{\gamma_j}$ for all $j$, is at most $1$.
    \item Since $\psi$ is a uniformly expanding, piecewise $C^2$ and a full branch map with finitely many branches, $\psi$ is exact (cf. \cite[Theorem 3]{HofbKell}) and $\psi^{-1}x$ is finite for all $x$.
    \item The Assumption (A)(1) in \cite{FP} is trivially true because there is only a single dynamical system in Figure 2 of \cite{FP}.
\end{itemize}

Hence, \ref{cond5} and \ref{cond6} are true due to the first part of \cite[Lemma 4.5]{FP}. To  establish \ref{cond7}, we need a slight modification of the second part of \cite[Lemma 4.5]{FP}. First, we note that $\chi \in \V \hookrightarrow L^2\,,$ for $\gamma\geq 3\,,$ and $\wh \psi$ has a spectral gap on $\V$. So, we can repeat the argument in the first part of the proof of \cite[Lemma 4.5]{FP} to  conclude that $\sum^{n-1}_{k=0}\bar\chi \circ \psi^k$ is $L^2-$bounded. So, it has an $L^2-$weakly convergent subsequence. This establishes \ref{cond7}. 

Finally, the non-arithmeticity of $\chi$ implies that $\chi$ is not cohomologous to a constant, and hence, we have \ref{cond8} of \Cref{prop:AbstractThm}. 
\end{proof}

\begin{proof}[Proof of Theorem \ref{thm:IntervalLLT}]
To prove this theorem we use \Cref{prop:AbsLLT}. By \Cref{thm:IntervalCLTsimple} we immediately obtain  \ref{cond5A} of \Cref{prop:AbsLLT}.

Next, we define the following chain of spaces. 
 \begin{align*}
 \VVc{\alpha_0}{\beta}{M} &\hookrightarrow \VVc{\alpha_1}{\beta}{M}  \hookrightarrow L^p \hookrightarrow L^1(\pi)
  \end{align*}
with $p \leq  M$ where the choices correspond to
$0\leq \alpha_0<\alpha_1<\beta$ and $\gamma_0=\gamma_1=M \geq 1$ in the proof of \Cref{thm:IntervalEdgeExp}.  
Then, the conditions \ref{cond1A}--\ref{cond4A} and \ref{cond6A}  of \Cref{prop:AbsLLT} follow as in the proof of \Cref{thm:IntervalEdgeExp} due to  {\Cref{cor:RegOpr} (2)} and \cite[Lemma 4.5]{FP}.
\end{proof}
\begin{proof}[Proofs of the results in Example \ref{ex: osc function}]
We first note that 
\begin{align*}
|\chi'(x)|\lesssim x^{-c}(1-x)^{-c}
\end{align*}
and
\begin{align*}
 |\chi'(x)|
 =\left|-cx^{-c-1}\sin\left(\frac{1}{x}\right)-x^{-c-2}\cos\left(\frac{1}{x}\right)\right|
 \lesssim x^{-c-2}(1-x)^{-c-2}.
\end{align*}
So, we obtain $a=c$ and $b=c+2$ in the notation of Theorems \ref{thm:IntervalCLTsimple}, \ref{thm:IntervalLLT} and \ref{thm:IntervalEdgeExp}. In order to prove (1) we note that \eqref{eq: CLT main cond} then simplifies to 
\begin{align*}
 c<\min\left\{\vartheta, \frac{1}{2+c}\right\}\min\left\{1, \frac{\log\eta_-}{\log \eta_+}\right\}.
\end{align*}
So, on the one hand, we have the requirement 
$c<\vartheta \tilde{\eta}$ and on the other hand, we have the condition 
$c< \tilde{\eta}/(c+2)$ which, given that we assume  $c\geq 0$\,, is equivalent to $c<\sqrt{1+\tilde{\eta}}-1$ giving (1). Furthermore, in the doubling map case we have $\vartheta=2$ and $\tilde{\eta}=1$ implying (1a). 

Next, we notice that \eqref{eq: Edge main cond} in our case simplifies to 
\begin{align*}
  3c<\min \left\{{\vartheta}, \frac{1}{2+c}\right\}\min\left\{1,\frac{\log\eta_+}{\log \eta_-}\right\}\,.
\end{align*}
With a similar calculation as above applying  Theorem \ref{thm:IntervalEdgeExp} gives (2) and as above we get (2a).
\end{proof}

\subsection{Proofs of limit theorems for the Boolean-type transformation}
\label{sec:BLimitThms}
Now we give the proofs from \Cref{subsec: Boolean trafo}.
We start with the following technical lemmas:
\begin{lem}\label{lem:EqualMoments}
For all $r\in \naturals$, the $r$th asymptotic moments of both $S_n(\chi)$ and $\widetilde{S}_n(h)$ are equal.
\end{lem}
\begin{proof}
It is enough to show that $\EXP_{\mu}(\widetilde S^{r}_n(h)) = \EXP_{\lambda_I}(S^{r}_n(\chi))$ for all $r$. In fact, due to \eqref{eq:LambdaChi}
\begin{align*}
    \EXP_{\mu}( h \circ \phi^{j_1}\, h \circ \phi^{j_2}\, \cdots \, h \circ \phi^{j_k}) &= \EXP_{\lambda_I}( h \circ \xi \circ \psi^{j_1}\, h \circ \xi \circ  \psi^{j_2}\, \cdots \, h \circ \xi \circ \psi^{j_k}) \\ &=\EXP_{\lambda_I}( \chi \circ \psi^{j_1}\, \chi \circ  \psi^{j_2}\, \cdots \, \chi \circ \psi^{j_k})
\end{align*}
for all $j_1\,, \dots\,, j_k \in \naturals_0$ such that $j_1+\dots+j_k={r}$.
\end{proof}

\begin{lem}\label{lem:uvtoalphabeta}
 Let $h:\mathbb{R}\to\mathbb{R}$
 be such that the left and right derivatives exist and there exist $u,v \geq 0$ fulfilling 
 \begin{align*}
  h\left(x\right)&\lesssim \left|x\right|^{u}\qquad\text{ and }\qquad
  \max\left\{\left|h'\left(x-\right)\right|, \left|h'\left(x+\right)\right|\right\}\lesssim \left|x\right|^{v},
 \end{align*}
and let $\chi:I\to \reals$ be given by $\chi=h\circ\xi$ with $\xi(x) := \cot (\pi x)$, then we have 
  $$|\chi(x)|\lesssim x^{-u}(1-x)^{-u}$$
 and
 $$\max\{|\chi'(x+)|, |\chi'(x-)|\}\lesssim x^{-b}(1-x)^{-b},\quad b=2+v\,.$$
 Further, if 
 \begin{align}\label{eq: cond alpha beta a b 1}
 \alpha > u\,&\,,\nonumber\\ 
  \beta < (1+ \alpha -u)/(2+v-u)\quad &\text{or}\quad 1+v<u\,, \quad\text{and}\\ 
 1 \leq \gamma <\, &1/u \,,\nonumber
 \end{align}
 then
 $\|\chi\|_{\alpha,\beta,\gamma}<\infty\,.$ In particular, if 
$u<1/(2+v-u)$, then there exist $0<\alpha<\beta<1$ such that $|\chi|_{\alpha,\beta}<\infty$.
 \end{lem}
\begin{proof} 
We will apply Lemma \ref{lem:abtoalphabeta}. First, we note that $$\lim_{x\to 0}\xi(x)x=1/\pi\quad\text{and}\quad\lim_{x\to 1}\xi(x)(1-x)=1/\pi\,.$$
 This and \eqref{eq: boundedness condition} imply 
 \begin{align}
  |\chi(x)|\lesssim x^{-u}(1-x)^{-u},
 \end{align}
 and in particular, $\chi \in L^\gamma$ with $1 \leq \gamma<1/u$.  
 
 For simplicity, we assume $\chi$ is differentiable. Otherwise, at a point where $\chi$ is not differentiable, both one-sided derivatives will exist and the following estimates do hold for them. 
 
 Note that we have $\left|{ h}'\left(\xi\left(x\right)\right)\right|\lesssim x^{-v} \left(1-x\right)^{-v}$. 
Using the chain rule $|\chi'(x) |= |{ h}'(\xi (x))||\xi'(x)|\,.$ Since
  $\xi'(x)=-\pi/\sin^2\left(\pi x\right)\,,$ we have that
\begin{align}
 |\chi'(x)|\lesssim x^{-2-v}(1-x)^{-2-v}.\label{eq: F cdot xdelta1}
\end{align}
So, we have $ |\chi'(x)|\lesssim x^{-b}(1-x)^{-b}$ with $b=2+v > 2$. 
The lemma then follows immediately by applying Lemma \ref{lem:abtoalphabeta}.
\end{proof}

With this we are able to prove the results from \Cref{subsec: Boolean trafo}

\begin{proof}[Proof of Proposition \ref{prop:CLTonReals}]
    To prove the statement it is enough to prove its counterpart for $S_n(\chi,\psi)$  where $\chi = h\circ \xi$ and $\psi$ is the doubling map. 
    
    From \Cref{lem:uvtoalphabeta}, we have 
 $$ |\chi(x)|\lesssim x^{-u}(1-x)^{-u}\quad\text{ and }\quad  \max\{|\chi'(x+)|, |\chi'(x-)|\}\lesssim x^{-b}(1-x)^{-b},\quad b=2+v\,.$$
Now, we invoke \Cref{thm:IntervalCLTsimple} with $\psi$, $\eta_+=\eta_-=2$ and  $\log \eta_-/\log\eta_+ =1$. Since $\psi$ is linear, $\vartheta =1$. 
Hence, \eqref{eq: CLT main cond} simplifies to 
$u<1/(2+v)$.
 Also, the assumption that $h$ is not an $L^2(\mu)$ coboundary implies that $\chi$ is not an $L^2(\lambda)$ coboundary. 
    
 Therefore, $\chi$ and $\psi$ satisfy the conditions of \Cref{thm:IntervalCLTsimple}, and hence satisfy the CLT given by \eqref{eq:IntervalCLT} with $$\sigma^2 =\EXP_{\lambda}(\chi^2)+ 2\sum_{k=1}^\infty \EXP_\lambda(\chi\cdot \chi\circ {\psi}^k) \in (0,\infty)\,.$$ From \Cref{lem:EqualMoments}, $\widetilde{\sigma}^2 = {\sigma}^2 $ and $\EXP_{\mu}(h)=\EXP_{\lambda_I}(\chi)$. As a direct consequence of \eqref{eq:LambdaChi}, we obtain the required CLT given by \eqref{eq:RealCLT}. 
\end{proof}

We next prove the MLCLT for a class of observables in $\fF$.
\begin{proof}[Proof of \Cref{prop:MLCLTonReals}]
Our assumption allows us to apply \Cref{thm:IntervalLLT} to the Birkhoff  sum $S_n(\chi)=\sum_{k=0}^{n-1}\chi\circ \psi^k$ with $\chi = h\circ \xi$ and $\psi$ the doubling map 
and conclude
    \begin{equation*}
  \sup_{\ell \in \reals }\left|\sigma\sqrt{2\pi n}\,\EXP_{\xi_* {\overline{m}}}(U \circ \xi  \circ \psi^n \,V(S_n(\overline\chi)-\ell)\, W \circ \xi)  - e^{-\frac{\ell^2}{2n\sigma^2}} \,\EXP_\pi(U\circ 
 \xi)\,\EXP_{\xi_* {\overline{m}} }(W \circ \xi )\int V(x) \, \mathrm{d}x \right|= o(1)\,.
   \end{equation*}
From \Cref{lem:EqualMoments} and the fact that $\xi$ is a conjugacy, we have 
 \begin{equation*}
   \sup_{\ell \in \reals }\left|\widetilde\sigma\sqrt{2\pi n}\,\EXP_{{\overline{m}}}(U \circ \phi^n \,V(\widetilde S_n(\overline{h})-\ell)\, W )  - e^{-\frac{\ell^2}{2n\widetilde \sigma^2}} \,\EXP_\mu(U)\,\EXP_{{\overline{m}}}(W )\int V(x) \, \mathrm{d}x \right|= o(1)\,.
   \end{equation*}
This is because the two LHSs are exactly the same. 
\end{proof}

Now, we prove that corollaries that show the validity of the CLT and MLCLT for the real part, imaginary part and the absolute value of the Riemann zeta function when sampled over the trajectories of $\phi$.  

\begin{proof}[Proof of Corollary \ref{cor:CLTforRZ}]
To apply \Cref{prop:CLTonReals}, we have to show the existence of $u,v$ as in \eqref{eq: boundedness condition}. It is well-known that 
for any $s \in (0,1)$,  for any $\delta > 0$,
\begin{equation}\label{eq:RZdecay}
    \max\left\{|\zeta|(s+ix), |\zeta '|(s+ix)\right\}\lesssim |x|^{(1-s)/2 + \delta}\,;
\end{equation}
see, for example, \cite{ECT}. 

So, we pick $u=v=(1-s)/2 + \delta$ and this is possible when $\left((1-s)/2 + \delta\right) \left((1-s)/2 + \delta +2\right) < 1$ and such $\delta >0$ exists iff $(1-s)(5-s)<4$ iff $s \in (3 - 2 \sqrt{2}, 1)$. So, for such choices of $s$ we can apply \Cref{prop:CLTonReals} and obtain the CLT provided that $h$ is not $\phi-$cohomologous to a constant. The MLCLT follows from \Cref{prop:MLCLTonReals} analogously, when $\phi$ is non-arithmetic.
\end{proof}

\begin{proof}[Proof of Corollary \ref{cor:CLTforpowerRZ}] To apply \Cref{prop:CLTonReals}, we have to show the existence of $u,v$ as in \eqref{eq: boundedness condition}. We assume $a \geq 1$,  set $\tilde{h}(x)=h(x)^{1/a}$.
    Note that  $h'(x) = a \tilde{h}(x)^{a-1}
    \tilde{h}'(x).$
   Since we restrict ourselves to the critical line, $s=1/2$, $|\tilde h(x) | \lesssim |x|^{13/84+\delta} $ and $|\tilde h'(x) | \lesssim |x|^{13/84+\delta}$ for all $\delta>0$, due to \eqref{eq:RZdecay}. So, we can take $u=13a/84+\delta$ and $v= 13(a-1)/84 + 13/84+ \delta = 13a/84 + \delta\,,$ 
   and 
    the condition in  \Cref{prop:CLTonReals} for $u,v$ reduces to 
    $(13a/84)(13a/84+2)  < 1\,.$
    This is equivalent to $1 \leq a < 84/13 (\sqrt{2}-1)$. So, for such choices of $a$, we can apply \Cref{prop:CLTonReals} and obtain the CLT provided that $h$ is not $\phi-$cohomologous to a constant. The MLCLT follows from \Cref{prop:MLCLTonReals} analogously, when $\phi$ is non-arithmetic.
\end{proof}

Finally, we look at the proof for the First Order Edgeworth Expansion for observables over the Boolean-type transformation.

\begin{proof}[Proof of \Cref{prop:EdgeworthonReals}]
We follow the proof of \Cref{prop:CLTonReals} and invoke \Cref{thm:IntervalEdgeExp}.  

Consider $S_n(\chi,\psi)$  where $\chi = \xi \circ h$ and $\psi$ is the doubling map. Remember that from \Cref{lem:uvtoalphabeta}, we have 
 $$ |\chi(x)|\lesssim x^{-u}(1-x)^{-u}\quad\text{ and }\quad  \max\{|\chi'(x+)|, |\chi'(x-)|\}\lesssim x^{-b}(1-x)^{-b},\quad b=2+v\,.$$

Next, to apply \Cref{thm:IntervalEdgeExp} we observe that $\eta_+=\eta_-=2$ and  $\log\eta_-/\log\eta^+=1$ and since $\psi$ is linear $\vartheta=1$. 
Hence, \eqref{eq: Edge main cond} simplifies to \eqref{eq: Edge Boole main cond}.
Also, the assumption that $h$ is not an $L^2(\mu)$ coboundary implies that $\chi$ is not an $L^2(\lambda)$ coboundary. 
\end{proof}
\appendix
\section{The Banach Spaces \texorpdfstring{$\V$}{V}}\label{subsec: Banach space corr}
The spaces $\mathsf{V}_{\alpha,\beta}$ with their particular norm considered in \cite{STI} are not complete, and thus, are not Banach spaces. 
However, with the norm we introduce here, we can construct a family of Banach spaces $\V$, $\alpha \in (0,1), \beta \in (0,1]$ and $\gamma \geq 1$, and use it to correct the proofs in \cite{STI}, and even generalize the results appearing there. 

First, we show that $\left\|\cdot\right\|_{\alpha,\beta,\gamma}$ is indeed  a norm. 
\begin{lem}\label{lem: norms corr}
For all $\alpha\in(0,1)$, $\beta\in(0,1]$ and $\gamma \geq 1$, we have that
 $\left\|\cdot\right\|_{\alpha,\beta,\gamma}$ is a norm. 
\end{lem}
\begin{proof}
 We have for $f,g\in \mathsf{V}_{\alpha,\beta}$ that 
\begin{align*}
 \left|f+g\right|_{\alpha,\beta}
 &=\sup_{\eps\in (0,\eps_0]} \int \frac{\osc\left(R_{\alpha}\left(f+g\right), B_{\eps}(x)\right)}{\eps^{\beta}}\mathrm{d}\lambda_I(x)\\
 &=\sup_{\eps\in(0,\eps_0]} \int \frac{\osc\left(R_{\alpha}f+R_{\alpha}g, B_{\eps}(x)\right)}{\eps^{\beta}}\mathrm{d}\lambda_I(x)\\
 &\leq \sup_{\eps\in (0,\eps_0]} \int \frac{\osc\left(R_{\alpha}f, B_{\eps}(x)\right)}{\eps^{\beta}}\mathrm{d}\lambda_I(x)
 +\sup_{\eps\in(0,\eps_0]} \int \frac{\osc\left(R_{\alpha}g, B_{\eps}(x)\right)}{\eps^{\beta}}\mathrm{d}\lambda_I(x)\\
 &=\left|f\right|_{\alpha,\beta}+\left|g\right|_{\alpha,\beta}
\end{align*}
and thus
 \begin{align*}
  \left\|f+g\right\|_{\alpha,\beta,\gamma}
  =\left\|f+g\right\|_{\gamma}+\left|f+g\right|_{\alpha,\beta} 
  \leq \left\|f\right\|_{\gamma}+\left\|g\right\|_{\gamma}+\left|f\right|_{\alpha,\beta}+\left|g\right|_{\alpha,\beta}
  =\left\|f\right\|_{\alpha,\beta,\gamma}+\left\|g\right\|_{\alpha,\beta,\gamma}.
 \end{align*}
It is obviously true that $\left\|af\right\|_{\alpha,\beta,\gamma}=a\left\|f\right\|_{\alpha,\beta,\gamma}$, for any positive $a$. Since $\left\|\cdot\right\|_{\gamma}$ is already a norm and $\left|f\right|_{\alpha,\beta}=0$ if $f=0$ almost surely, we
know that $\left\|f\right\|_{\alpha,\beta,\gamma}=0$ if and only if $f=0$ almost surely.
\end{proof}

\subsection{Completeness} 
Here we verify that $\mathsf{V}_{\alpha,\beta,\gamma}$ are, in fact, Banach spaces.
\begin{lem}\label{lem: completeness corr}
For $\alpha\in(0,1), \beta\in (0,1]$ and $\gamma \geq 1$,
 $\mathsf{V}_{\alpha,\beta,\gamma}$ is complete. 
\end{lem}
\begin{proof}
Let $\left(f_n\right)$ be a Cauchy sequence with respect to $\left\|\cdot\right\|_{\alpha,\beta,\gamma}$. 
Then, in particular $\left(f_n\right)$ is also a Cauchy sequence with respect to $\left\|\cdot\right\|_{\gamma}$,
we set $f$ as its limit. Also, there exists a subsequence, say $(f_{n_r})$, that converges to $f$ pointwise almost everywhere.   

Since $\left(f_n\right)$ is a Cauchy sequence with respect to $\left\|\cdot\right\|_{\alpha,\beta,\gamma}$,  for each $\delta>0$ we can choose $L>0$ such that $\left\|f_{k}-f_{\ell}\right\|_{\alpha,\beta,\gamma}<\delta$ for all $k,\ell>L$. Let $\delta>0$ and choose $k,\ell$ sufficiently large so that $n_k, n_\ell>L$. Then, 
\begin{align*} 
 \left\|f_{n_k}-f_{n_\ell}\right\|_{\alpha,\beta,\gamma}=\left\|f_{n_k}-f_{n_\ell}\right\|_{\gamma}+\sup_{\eps\in(0,\eps_0]}\frac{\int\osc\left(R_{\alpha}\left(f_{n_k}-f_{n_\ell}\right),B_{\eps}(x)\right)\;\mathrm{d}\lambda_I\left(x\right)}{\eps^{\beta}} < \delta\,.
\end{align*}
Then, by Fatou's Lemma, $\left\|f_{n_k}-f\right\|_{\gamma} \leq \liminf_{\ell \to \infty}\left\|f_{n_k}-f_{n_\ell}\right\|_{\gamma}$ and 
\begin{align*}
    \frac{\int \osc\left(R_{\alpha}\left(f_{n_k}-f\right),B_{\eps}(x)\right)\;\mathrm{d}\lambda_I\left(x\right)}{\eps^{\beta}} & \leq \frac{\int\liminf_{\ell \to \infty}\osc\left(R_{\alpha}\left(f_{n_k}-f_{n_\ell}\right),B_{\eps}(x)\right)\;\mathrm{d}\lambda_I\left(x\right)}{\eps^{\beta}} \\
    &\leq \liminf_{\ell \to \infty}\frac{\int\osc\left(R_{\alpha}\left(f_{n_k}-f_{n_\ell}\right),B_{\eps}(x)\right)\;\mathrm{d}\lambda_I\left(x\right)}{\eps^{\beta}} \\ 
    &\leq \liminf_{\ell \to \infty} \sup_{\eps\in(0,\eps_0]}  \frac{\int\osc\left(R_{\alpha}\left(f_{n_k}-f_{n_\ell}\right),B_{\eps}(x)\right)\;\mathrm{d}\lambda_I\left(x\right)}{\eps^{\beta}}\,.
\end{align*}
As a result, for all $k$ sufficiently large so that $n_k >L$,
\begin{align*} 
 \left\|f_{n_k}-f\right\|_{\alpha,\beta,\gamma}&
 \leq\liminf_{\ell \to \infty}\left\|f_{n_k}-f_{n_\ell}\right\|_{\gamma}+\liminf_{\ell \to \infty} \sup_{\eps\in(0,\eps_0]}  \frac{\int\osc\left(R_{\alpha}\left(f_{n_k}-f_{n_\ell}\right),B_{\eps}(x)\right)\;\mathrm{d}\lambda_I\left(x\right)}{\eps^{\beta}}\\ &\leq \liminf_{\ell \to \infty} \left(\left\|f_{n_k}-f_{n_\ell}\right\|_{\gamma} + \sup_{\eps\in(0,\eps_0]}  \frac{\int\osc\left(R_{\alpha}\left(f_{n_k}-f_{n_\ell}\right),B_{\eps}(x)\right)\;\mathrm{d}\lambda_I\left(x\right)}{\eps^{\beta}}\right) \leq \delta.
\end{align*}
Now, choose $r$ sufficiently large so that $n_r>L$ and $k>L$. Then,
\begin{align*}
 \left\|f_{k}-f\right\|_{\alpha,\beta, \gamma} &\leq \left\|f_{k}-f_{n_r}\right\|_{\alpha,\beta, \gamma} + \left\|f_{n_r}-f\right\|_{\alpha,\beta, \gamma}
 < 2\delta\,.
\end{align*}
Thus, $f\in\mathsf{V}_{\alpha,\beta,\gamma}$ and $\left(f_{n}\right)$ converges to $f$ with respect to $\left\|\cdot\right\|_{\alpha,\beta, \gamma}$ giving completeness.
\end{proof}

Now, we discuss properties of $\V$ that are relevant for the application of \Cref{prop:AbstractThm} to our setting. First, we prove that constant functions belong to the spaces we consider.
\begin{lem}\label{lem:Const}
For $\alpha\in (0,1)$, $\beta \in (0,1]$ and $\gamma \geq 1,$ the constant function, $\One \in \V$.
\end{lem}

\begin{proof}
Since $\|\One\|_\gamma =1$, we only have to show that $ |\One|_{\alpha,\beta}<\infty$. Observe that $R_\alpha\One$ is bounded by $2^{-2\alpha}$, symmetric about $x=1/2$ and strictly increasing on $[0,1/2]$ with a strictly decreasing derivative. Hence, for any $0<\eps \leq \eps_0 <1/4 $,
\begin{align*}
    \int\osc(R_{\alpha}\One , B_\eps(x)) \mathrm{d}\lambda_I(x) &\leq \int_{2\eps}^{1-2\eps}{\osc(R_{\alpha}\One , B_\eps(x))}\, \mathrm{d}\lambda_I(x) + 2^{-2\alpha} \left(\int_{0}^{2\eps}\mathrm{d}\lambda_I(x)+\int_{1-2\eps}^1\, \mathrm{d}\lambda_I(x)\right) \\ &\leq 4\eps \int_{2\eps}^{1/2} \max_{B_\eps(x)} |(R_{\alpha}\One)'| \, \mathrm{d}\lambda_I(x) + 2^{2-2\alpha} \eps
    \\ 
    & = 4\eps \int_{2\eps}^{1/2}  (R_{\alpha}\One)'(x-\eps) \, \mathrm{d}\lambda_I(x) + 2^{2-2\alpha} \eps\\
    &= 4\eps\left(R_{\alpha}\One\big(1/2-\eps\big) - R_{\alpha}\One(\eps)\right) + 2^{2-2\alpha} \eps \leq 2^{3-2\alpha}\eps.
\end{align*}
This implies that $|\One|_{\alpha,\beta}\leq 2^{3-2\alpha}\eps^{1-\beta}_0$. 
\end{proof}

Next, we state two lemmas about the inclusion properties of $\V$.

\begin{lem}\label{lem:CptIncl1}
For $\beta \in (0,1]$ and $\gamma \geq 1\,,$ $$\VVc{0}{\beta}{\gamma}\hookrightarrow \VVc{0}{\beta}{1} \hookrightarrow L^\infty\,.$$ 
\end{lem}

\begin{proof}
This follows from \cite[Proposition 3.4]{BS} applied to the real and imaginary parts of functions in $\VVc{0}{\beta}{1}$ and the fact that $ L^\gamma \hookrightarrow L^1$. 
\end{proof}

\begin{rem}\label{rem:gInVRgInV}
Note that, if $f \in \V$, then $R_\alpha f \in \VVc{0}{\beta}{\gamma}$. So, $\esssup R_\alpha f < \infty$. This fact will be useful in proofs. 
\end{rem}

\begin{lem}\label{lem:CptIncl2}
Suppose $0<\alpha_1\leq \alpha_2 < 1\,,$ $0<\beta_2 \leq \beta_1 \leq 1$ and $1 \leq \gamma_2 \leq  \gamma_1$. Then $$\mathsf V_{\alpha_1,\beta_1, \gamma_1} \hookrightarrow \mathsf V_{\alpha_2,\beta_2, \gamma_2} \hookrightarrow L^1\,.$$ 
\end{lem}

\begin{proof}
Since $\|f\|_{\gamma_2} \leq \|f\|_{\gamma_1}$, it is enough to show that $|f|_{\alpha_2,\beta_2} \lesssim  \|f\|_{\alpha_{1},\beta_{1},\gamma_{1}}$. By applying \cite[Proposition 3.2 (iii)]{BS} to the real and imaginary parts of $f$, we have,
\begin{align*}
    {\osc(R_{\alpha_2} f, B_\eps(x))} &= {\osc(R_{\alpha_2-\alpha_1}\One \cdot  R_{\alpha_1}f, B_\eps(x))}\\
    &\leq \esssup |R_{\alpha_1}f| \cdot \osc(R_{\alpha_2-\alpha_1}\One,B_\eps(x) ) +\osc(R_{\alpha_1}f, B_\eps(x)) \cdot \sup_{B_\eps(x)} R_{\alpha_2-\alpha_1}\One \,,
\end{align*}
and due to \Cref{lem:CptIncl1}, 
$$\esssup |R_{\alpha_1}f| \lesssim |R_{\alpha_1}f|_{0,\beta_1} + \|R_{\alpha_1} f\|_{1} \leq |f|_{\alpha_{1},\beta_{1}} + \|R_{\alpha_1}{\One}\|_{\bar \gamma}\|f\|_{\gamma_1} \lesssim \|f\|_{\alpha_{1},\beta_{1},\gamma_{1}}$$
with $\bar \gamma = (1-\gamma^{-1}_1)^{-1}$. Therefore,
\begin{align*}
    &\eps^{-\beta_2}{\osc(R_{\alpha_2} f, B_\eps(x))} \\&\qquad\lesssim  \eps^{-\beta_1}{\osc(R_{\alpha_2-\alpha_1}\One , B_\eps(x))}\|f\|_{\alpha_{1},\beta_{1},\gamma_{1}}+
    \sup_{ B_{\eps}(x)} R_{\alpha_2-\alpha_1}\One \cdot \eps^{-\beta_1}{\osc(R_{\alpha_1}f, B_\eps(x))}\,.
\end{align*}
Integrating and taking the supremum over $\eps$, 
$$|f|_{\alpha_2,\beta_2}  \lesssim \|f\|_{\alpha_{1},\beta_{1},\gamma_{1}},$$
and the inclusion follows.\end{proof}

\subsection{Continuous inclusion and relative compactness}\label{sec:CtsIncl} To apply Hennion-Nassbaum theory, see \cite{KL, DKL}, we have to show that our \textit{weak} spaces, $L^p$, are continuously embedded in \textit{strong} spaces, $\V$, and that the closed bounded sets in strong spaces are compact with respect to weak norms. 

\begin{lem} Let $\alpha\in(0,1), \beta\in (0,1]$ and $\gamma \geq 1$.
Then for all $\bar \gamma$ such that $\gamma< \bar \gamma < 1/\alpha\,,$ $L^{\bar \gamma}$ is continuously embedded in $\V$. 
\end{lem}
\begin{proof}
Due to \Cref{rem:LpInclusion} and the assumption $\bar\gamma < 1/\alpha$, if $h \in \V\,,$ then $h \in L^{\bar\gamma}$. So, $\V \subseteq L^{\bar\gamma}$. To show that this inclusion is continuous we need to show that if $f_n \to 0$ in $\V$, then $f_n \to 0$ in $L^{\bar \gamma}$. Let $\|f_n\|_{\alpha,\beta,\gamma} \to 0$. Then, $|R_\alpha f_n| \in \VVc{0}{\beta}{1}$ and $\|R_\alpha f_n\|_{0,\beta,1}\to 0$. However, $\VVc{0}{\beta}{1} \hookrightarrow L^\infty$. So, $\|R_\alpha f_n \|_{\infty} \to 0$. Therefore, $ \|f_n^{\bar \gamma}\|_1 \leq  \|R_{-\alpha\bar\gamma}\One\|_1 \|R_\alpha f_n\|_\infty^{\bar\gamma} \to 0 $ proving the claim. 
\end{proof}

\begin{lem} Let $\alpha, \beta, \gamma$ and $\bar \gamma$ be as in the previous lemma. Then, the closed unit ball of $\V$ is compact in $L^{\bar\gamma}$.
\end{lem}

\begin{proof}
Let $\{f_n\}$ be such that $\|f_n\|_{\alpha,\beta,\gamma} \leq 1$. It is enough to show that there is $f \in \V$ such that $\|f\|_{\alpha,\beta,\gamma} \leq 1$  and $\{f_n\}$ converges to $f$ in $L^{\bar\gamma}$ over a subsequence.  
To do this, we recall from \cite[Theorem 1.13]{keller_generalized_1985} that closed subsets of $\VVc{0}{\beta}{\gamma}$ are compact in $L^{\gamma}$. Since $\{R_\alpha f_n \} \subset \VVc{0}{\beta}{\gamma}$ is a bounded sequence, it has an $L^{\gamma}$ convergent subsequence, and in turn, it has a pointwise almost everywhere convergence subsequence. Let's call this subsequence $\{R_\alpha f_{n_k}\}$ and its point-wise limit $f$. 

We claim $f_{n_k} \to R_{-\alpha} f$ in $L^{\bar\gamma}$. Observe that $f_{n_k} \to R_{-\alpha} f $ point-wise almost everywhere, and since $\VVc{0}{\beta}{\gamma} \hookrightarrow L^{\infty}\,,$ $|f_{n_k}| \leq  |R_{-\alpha}\One| |R_\alpha f_{n_k}| \leq C |R_{-\alpha}\One| \in L^{\bar \gamma}\,.$ 
So, $f_{n_{k}}\to R_{-\alpha}f$ in $L^{\bar \gamma}$ if $\alpha \bar\gamma < 1\,.$ 
Moreover, we claim $\|R_{-\alpha}f\|_{\alpha,\beta,\gamma} \leq 1$. To see this, observe that since $L^{\bar \gamma}$ convergence implies $L^\gamma$ convergence, we apply \cite[Lemma 1.12]{keller_generalized_1985} to conclude that $\liminf_{k}| f_{n_k}|_{\alpha,\beta} = \liminf_{k}|R_\alpha f_{n_k} |_{0,\beta} \geq |f|_{0,\beta} = |R_{-\alpha}f|_{\alpha, \beta}$.  Since strong convergence implies weak convergence, we have $\liminf_{k} \|f_{n_k}\|_{\gamma} \geq  \|R_{-\alpha}f\|_{\gamma}\,,$ and finally,  
\begin{align*}
    \|R_{-\alpha}f\|_{\alpha,\beta,\gamma} = |R_{-\alpha}f|_{\alpha,\beta} +  \|R_{-\alpha}f\|_{\gamma} &\leq \liminf_{k}| f_{n_k}|_{\alpha,\beta} + \liminf_{k} \|f_{n_k}\|_{\gamma} \\ &\leq \liminf_{k} \left(| f_{n_k}|_{\alpha,\beta}+ \|f_{n_k}\|_{\gamma}\right)= \liminf_{k} \|f_{n_k}\|_{\alpha,\beta,\gamma} \leq 1\,\end{align*}
as claimed.
\end{proof}

\begin{rem}\label{rem:ToApplyHN}
In particular, the above implies that $\|\cdot\|_{\alpha,\beta,\gamma}$-bounded sequences have $\|\cdot\|_{\bar\gamma}$-Cauchy subsequences.
\end{rem}

\subsection{Multiplication in \texorpdfstring{$\V$}{V}}\label{sec:Spaces}

\subsubsection{Multiplication by \texorpdfstring{$e^{is\chi}$}{exp isChi}}\label{sec:BanachSpaces} In this section, we prove some properties of multiplication by $e^{is\chi}$ in  $\V$ that are necessary for our proofs. 

Observe that the spaces $\V$, as opposed to spaces usually used in ergodic theory such as $L^\infty$, $\BV[0,1]$ or $C^1[0,1]$, are not Banach algebras. Hence, $s \mapsto \wh \psi_{is} \in \cL(\V)$ may not be continuous. The following lemma will allow us to establish its continuity as a function from $\reals$ to $\cL(\mathsf V_{\alpha_1,\beta_1, \gamma_1},\mathsf V_{\alpha_2,\beta_2, \gamma_2})$ for some good choices of indices. 

\begin{lem}\label{lem:NormEst}
Suppose $g \in \mathsf V_{\alpha_1,\beta_1, \gamma_1}$, $h \in \mathsf V_{\alpha_2,\beta_2, \gamma_2}$ and $\alpha_3=\alpha_1+\alpha_2$, $\beta_3 \leq \,\min\{\beta_1\,,\,\beta_2\}$ and $\gamma_3 \leq (\gamma_1^{-1}+\gamma_2^{-1})^{-1}$. Then,
$$  \|gh\|_{\alpha_{3},\beta_{3},\gamma_{3}} \lesssim  \|g\|_{\alpha_{1},\beta_{1},\gamma_{1}}\|h\|_{\alpha_{2},\beta_{2},\gamma_{2}}\,$$
with the proportionality constant independent of $g$ and $h$  but dependent on $\alpha_j, \beta_j,\gamma_j,\, j =1,2,3$. 
\end{lem}
\begin{proof} First, suppose $g$ and $h$ are real valued. Then
 \begin{equation}\label{eq: pos neg part}
   \osc(R_{\alpha}u, B_{\eps}(x))=\osc(R_{\alpha}u_-, B_{\eps}(x))+\osc(R_{\alpha}u_+, B_{\eps}(x))\,.  
 \end{equation}
 
By applying \cite[Proposition 3.2 (iii)]{BS} to the positive and negative parts of $g$,
\begin{align*}
    &\osc(R_{ \alpha_3} (gh),B_\eps(x))\\ &\quad\quad=\osc(R_{ \alpha_3} (g_+-g_-)h,B_\eps(x)) \\ &\quad\quad= \osc(R_{\alpha_{1}}(g_+-g_-) \cdot R_{\alpha_{2}}h, B_\eps(x)) \\
    &\quad\quad\leq \osc(R_{\alpha_{1}}g_+ \cdot R_{\alpha_{2}}h, B_\eps(x))+\osc(R_{\alpha_{1}}g_- \cdot R_{\alpha_{2}}h, B_\eps(x))\\
    &\quad\quad\leq \sum_{r = \pm } \Big(\osc(R_{\alpha_{1}}g_r, B_\eps(x))\cdot \esssup |R_{\alpha_{2}}h| + \osc(R_{\alpha_{2}}h, B_\eps(x))\cdot \esssup |R_{\alpha_{1}}g_r |\Big) \\ 
    &\quad\quad\leq  \osc(R_{\alpha_{1}}g, B_\eps(x))\esssup |R_{\alpha_{2}}h|+ 2 \cdot \osc(R_{\alpha_{2}}h, B_\eps(x))\esssup |R_{\alpha_{1}}g |\,. 
\end{align*}
If $g$ is complex valued, using the definition of $\osc$, we have 
\begin{align*}
    &\osc(R_{\alpha_3} (gh),B_\eps(x)) \\ &\quad\leq  \osc(R_{\alpha_{1}}g, B_\eps(x))\esssup |R_{\alpha_{2}}h|+ 2 \cdot \osc(R_{\alpha_{2}}h, B_\eps(x))(\esssup |R_{\alpha_{1}}\Re g |+\esssup |R_{\alpha_{1}}\Im g |)\,,
     \\ &\quad\leq  \osc(R_{\alpha_{1}}g, B_\eps(x))\esssup |R_{\alpha_{2}}h|+ 2 \sqrt{2}\cdot \osc(R_{\alpha_{2}}h, B_\eps(x)) \esssup |R_{\alpha_{1}} g |\,.
\end{align*}
If $h$ is not real valued, repeating the argument for the real and imaginary parts of $h$, we obtain 
\begin{align*}
    &\osc(R_{\alpha_3} (gh),B_\eps(x)) \\ 
    &\quad\leq  2\sqrt{2} \cdot \osc(R_{\alpha_{1}}g, B_\eps(x))\esssup |R_{\alpha_{2}}h|+ 2 \sqrt{2}\cdot \osc(R_{\alpha_{2}}h, B_\eps(x)) \esssup |R_{\alpha_{1}} g |\,.
\end{align*}

Now, we use the inclusion of $L^{\infty}$ in $\textsf{V}_{0,\beta_r,1}$ where $r=1,2$ to conclude that
\begin{align*}
    &\int \osc(R_{\alpha_3}(gh),B_\eps(x)) \, \mathrm{d}\lambda_I(x) \\ 
    &\qquad\lesssim \int \osc(R_{\alpha_{1}}g, B_\eps(x))\, \mathrm{d}\lambda_I(x)\cdot \|R_{\alpha_2}  h\|_{0,\beta_2,1}+  \int \osc(R_{\alpha_{2}}h, B_\eps(x)) \,\mathrm{d}\lambda_I(x) \cdot \|R_{\alpha_1} g\|_{0,\beta_1,1} \\
    &\qquad\lesssim \eps^{\beta_1}|g|_{\alpha_1,\beta_1} (|h|_{\alpha_2,\beta_2}+\|R_{\alpha_{2}}h\|_{1} )+  \eps^{\beta_2}|h|_{\alpha_2,\beta_2} (|g|_{\alpha_1,\beta_1}+\|R_{\alpha_{1}}g\|_{1} ).
\end{align*}
This gives us that for all $\eps \in (0,1]$,
\begin{align*}
&\eps^{-\beta_3}\int \osc(R_{\alpha_3}(gh),B_\eps(x)) \, \mathrm{d}\lambda_I(x)  \\ &\qquad\lesssim 
|g|_{\alpha_1,\beta_1}|h|_{\alpha_2,\beta_2}+|g|_{\alpha_1,\beta_1}\|h\|_{\gamma_2}+|h|_{\alpha_2,\beta_2} |g|_{\alpha_1,\beta_1}+|h|_{\alpha_2,\beta_2} \|g\|_{\gamma_1}\,.
\end{align*}
Taking the supremum over $\eps$ and combining with $\|gh\|_{ \gamma_3} \leq \|g\|_{\gamma_1} \|h\|_{\gamma_2}$ implies the result. 
\end{proof}

Due to the linearity of the operator $\wh\psi$, in order to show regularity of $s \mapsto \wh \psi_{is}=\psi(e^{is\chi} \times \cdot\, )$, it is enough to show the regularity of the one parameter group of multiplication operators $s \mapsto e^{is\chi}  \times  \cdot\,$. Our next lemma provides general conditions that guarantees this.

\begin{lem}\label{lem:MultOp}
Let $0 \leq \alpha_0,
\beta \leq 1 $ and $\gamma_0 \geq 1$. For each $s \in \reals$, consider the multiplication operator, $H_s(\cdot) = e^{is\chi}\times \cdot\,,$ on $\VVc{\alpha_0}{\beta}{\gamma_0}$. 
\begin{enumerate}[itemsep=3pt]
    \item[$(1)$] Suppose there is $\bar\beta \geq \beta$ such that, for all $s \in \reals$,  $|e^{is\chi}|_{0 ,\bar\beta}  < \infty$. Then, for all $s \in \reals$, $H_s \in \cL(\VVc{\alpha_0}{\beta}{\gamma_0})$\,.
    \item[$(2)$] Suppose$,$ in addition to the conditions in $(1),$ there exists $0<\alpha^*< \beta$ such that 
    \begin{equation}\label{eq:ContCond}
        \lim_{s \to 0}|1-e^{is\chi}|_{\alpha^*,\beta} = 0.
    \end{equation} Put $ \alpha_1 = \alpha_0+\alpha^*$ and $\gamma_1 \leq \gamma_0$. Then $s \mapsto H_s \in \cL(\VVc{\alpha_0}{ \beta}{\gamma_0}, \VVc{ \alpha_1}{ \beta}{\gamma_1})$ is continuous.
    \item[$(3)$]  Suppose, in addition to the conditions in $(1)$ and $(2)$, there exist $ 0<\alpha^{**}<\beta$ and $\gamma \geq 1$ such that  \begin{equation}\label{eq:DiffCond}
       \lim_{s \to 0}\left|\frac{e^{is\chi}-1-is\chi}{s}\right|_{\alpha^{**},\beta}=0\,\quad\text{and}\quad \|\chi\|_{\gamma}<\infty\,.  
        \end{equation} Put $\alpha_2=\alpha_0+\max\{\alpha^*, \alpha^{**}\}$ and $\gamma_2 \leq (\gamma_1^{-1}+\gamma^{-1})^{-1}$. Then, the function $s \mapsto H_s \in \cL(\VVc{\alpha_0}{\beta}{\gamma_0},\VVc{\alpha_2}{\beta}{\gamma_2})$ is differentiable with the derivative $$H^\prime_s (\cdot) = (i\chi)e^{is\chi}\,\times \,\cdot\,\,\,.$$
    \item[$(4)$] Suppose, the conditions in $(1), (2)$ and $(3)$ are true. Put $\alpha_3=\alpha_2+\alpha^*$ and $\gamma_3 \leq \gamma_2$. Then $s \mapsto H_s \in \cL(\VVc{\alpha_0}{ \beta}{\gamma_0}, \VVc{\alpha_3}{ \beta}{\gamma_3})$ is continuously differentiable. 
\end{enumerate}
\end{lem}

\begin{rem}
It would be possible to have some more flexibility on the parameter $\beta$ and change it for different spaces. However, we only use the version of the lemma as stated which also keeps a simpler notation. 
\end{rem}

\begin{proof}[Proof of Lemma \ref{lem:MultOp}]\,

\noindent \underline{Proof of $(1)$:}\\
We note that for all $g \in \VVc{\alpha_0}{\beta}{\gamma_0}$, $\|H_s(g)\|_{\gamma_0} = \|g\|_{\gamma_0}$ and due to \cite[Proposition 3.2 (iii)]{BS},
\begin{align*}
    \osc(R_{\alpha_0}(e^{is\chi}g), B_\eps(x)) 
    &\leq  \osc(R_{\alpha_0}(e^{is\chi}g_+), B_\eps(x))+\osc(R_{\alpha_0}(e^{is\chi}g_-), B_\eps(x))\\
    &\lesssim  \osc(R_{\alpha_{0}}g, B_\eps(x))+ \osc(e^{is\chi}, B_\eps(x))\cdot \esssup (|R_{\alpha_{0}}g|)\\ 
    &\lesssim  \osc(R_{\alpha_{0}}g, B_\eps(x))+ \osc(e^{is\chi}, B_\eps(x)) \|g\|_{\alpha_0,\beta,\gamma_0} \text{ and }
    \\
    \eps^{-\beta}\osc(R_{\alpha_0}(e^{is\chi}g), B_\eps(x)) &\lesssim  \eps^{-\beta}\osc(R_{\alpha_{0}}g, B_\eps(x))+ \eps^{- \beta} \osc(e^{is\chi}, B_\eps(x))\|g\|_{\alpha_0,\beta,\gamma_0}.
\end{align*}
The first $\lesssim$ is due to adding up the positive and negative part of $g$  the second is due to the inclusion $\VVc{0}{\beta}{\gamma_0} \hookrightarrow L^\infty$. 
Integrating and taking the supremum over $\eps$, we have
\begin{equation*}
    |H_s(g)|_{\alpha_0,\beta} \lesssim |g|_{\alpha_0,\beta}+|e^{is\chi}|_{0,\beta} \|g\|_{\alpha_0,\beta,\gamma_0} 
\end{equation*}
which gives
\begin{equation}\label{eq:ExpNormBound}
    \|H_s(g)\|_{\alpha_0,\beta,\gamma_0}  \leq (1+|e^{is\chi}|_{0, \beta})   \|g\|_{\alpha_0,\beta,\gamma_0}.
\end{equation}
Therefore, for all $s$, $H_s$ maps $\VVc{\alpha_0}{\beta}{\gamma_0}$ to itself, and is a bounded linear operator on $\VVc{\alpha_0}{\beta}{\gamma_0}$.\vs

\noindent \underline{Proof of $(2)$:}\\
We note that, 
$ H_t g - H_s g = (\Id-H_{s-t}) H_t g$ and if $g\in \VVc{\alpha_0}{\beta}{\gamma_0}$ then $H_tg \in \VVc{\alpha_0}{\beta}{\gamma_0}$. Hence, due to \Cref{lem:NormEst}, it is enough to prove that 
$$\lim_{s \to 0}\|\Id-H_{s}\|_{\VVc{\alpha_0}{\beta}{\gamma_0}\to \VVc{\alpha_1}{ \beta}{\gamma_1}} = 0.$$ 
To this end, let $g \in \VVc{\alpha_0}{\beta}{\gamma_0}$ be such that $\|g\|_{\alpha_0,\beta,\gamma_0}\leq 1$. Then,
$$\lim_{s\to 0}\|(\Id-H_{s})g\|^{\gamma_1}_{\gamma_1}=\lim_{s\to 0}\int |(1-e^{is\chi})g|^{\gamma_1}\, \, \mathrm{d}\lambda_I= 0 $$
by the dominated convergence theorem. 
Moreover, by \cite[Proposition 3.2 (iii)]{BS}
\begin{align*}
    &\osc( R_{\alpha_1}(\Id-H_{s})g, B_\eps(x) ) \\&\quad=  \osc(R_{\alpha^*}(1-e^{is\chi})R_{\alpha_0} g, B_\eps(x) ) \\
    &\quad \lesssim \osc(R_{\alpha_0}  g, B_\eps(x) ) \cdot \esssup |R_{\alpha^*}(1-e^{is\chi})| + \osc(R_{\alpha^*}(1-e^{is\chi}), B_\eps(x) ) \cdot \esssup |R_{\alpha_0} g|\,,
\end{align*}
 where $\lesssim$ is due to the fact that we have to consider the positive and negative part of $g$ separately.
Because $\VVc{0}{\beta}{1} \hookrightarrow L^\infty$, we have 
\begin{align*}
     \eps^{-\beta}\osc( R_{\alpha_1}(\Id-H_{s})g, B_\eps(x) )
     &\lesssim \eps^{-\beta}\osc(R_{\alpha_0}  g, B_\eps(x) ) (|1-e^{is\chi}|_{\alpha^*,\beta}+\| R_{\alpha^*}(1-e^{is\chi})\|_1)\\
     &\qquad+\eps^{-\beta} \osc(R_{\alpha^*}(1-e^{is\chi}), B_\eps(x) )\|g\|_{\alpha_0,\beta,\gamma_0}.
\end{align*}
Integrating, taking the sup over $\eps$, and finally, using $\|g\|_{\alpha_0,\beta,\gamma_0} \leq 1$, we get
\begin{align*}
    |(\Id-H_{s})g|_{\alpha_1,\beta}&\lesssim  |1-e^{is\chi}|_{\alpha^*,\beta}+\|R_{ \alpha^*}(1-e^{is\chi})\|_1.
\end{align*}
By the bounded convergence theorem $\lim_{s \to 0}\|R_{\alpha^*}(1-e^{is\chi})\|_1=0$. Therefore,  
$$ \lim_{s \to 0}|(\Id-H_{s})g|_{\alpha_1,\beta} = 0.$$
Hence, we have the continuity of $s \mapsto H_s$.\vs

\noindent \underline{Proof of $(3)$:}\\
First, we show that for all $g \in \mathsf V_{\alpha_1,\beta_1, \gamma_1}$ such that $  \|g\|_{\alpha_{1},\beta_{1},\gamma_{1}} \leq 1$,
$$\lim_{h \to 0}\Big\|\Big(\frac{H_{s+h}-H_s-H^\prime_s h}{h}\Big)g\Big\|_{\alpha_2, \beta, \gamma_2} = \lim_{h \to 0}\Big\|\Big(\frac{H_{h}-\Id-i\chi h}{h}\Big)H_sg\Big\|_{\alpha_2, \beta, \gamma_2} = 0.$$
Due to \Cref{lem:NormEst}, it is enough to show that 
$$\lim_{h \to 0}\Big\|\Big(\frac{H_{h}-\Id-i\chi h}{h}\Big)\One\Big\|_{\alpha^{**}, \beta, \gamma} = \lim_{h \to 0}\Big\|\frac{e^{ih\chi}-1-i\chi h}{h}\Big\|_{\alpha^{**}, \beta, \gamma}  =0.$$ From the dominated convergence theorem, we have
$$\lim_{h \to 0}\Big\|\frac{e^{ih\chi}-1-i\chi h}{h}\Big\|_{\gamma}  =0.$$
The assumption \eqref{eq:DiffCond} completes the proof of differentiability.

Finally, picking $h\neq 0$ sufficiently close to 0, 
applying the estimate in part (1), part (2) with $\gamma_1=\gamma_0$, and \Cref{lem:NormEst}, we note that for all $g \in \VVc{\alpha_1}{\beta}{\gamma_1}$ and for all $s\,,$
\begin{align*}
    \|H'_s(g)\|_{\alpha_2,\beta, \gamma_2} &=\left\|\frac{e^{ih\chi}-1-ih\chi}{h} e^{is\chi}g +\frac{1}{h}(1-e^{ih\chi})e^{is\chi}g\right\|_{\alpha_2,\beta, \gamma_2}\\
    &\leq  \left\|\frac{e^{ih\chi}-1-ih\chi}{h} e^{is\chi}g \right\|_{\alpha_2,\beta,\gamma_2}+\frac{1}{h}\|(1-e^{ih\chi})e^{is\chi}g\|_{\alpha_2,\beta,\gamma_2} \\
    &\leq \left\|\frac{e^{ih\chi}-1-ih\chi}{h}\right\|_{\alpha^{**},\beta,\gamma}\| H_s(g)\|_{\alpha_0,\beta,\gamma_0} + \frac{1}{h}\|(1-e^{ih\chi})\|_{\alpha^*,\beta,\gamma}\|H_s(g)\|_{\alpha_0,\beta,\gamma_0}\\
    & \lesssim \left(\left\|\frac{e^{ih\chi}-1-ih\chi}{h}\right\|_{\alpha^{**},\beta,\gamma}+ \|(1-e^{ih\chi})\|_{\alpha^*,\beta,\gamma}\right)(1+|e^{is\chi}|_{0,\beta})\|g\|_{\alpha_0,\beta,\gamma_0}
\end{align*} 
So, $H^\prime_s$ is, in fact, a bounded linear operator in $ \cL(\VVc{\alpha_0}{\beta}{\gamma_0},\VVc{\alpha_2}{\beta}{\gamma_2})$. \vs

\noindent\underline{Proof of $(4)$:}
Since $\VVc{\alpha_2}{\beta}{\gamma_2} \hookrightarrow \VVc{\alpha_3}{ \beta}{\gamma_3}$, we have that $s\to H_s \in \cL(\VVc{\alpha_0}{\beta}{\gamma_0},\VVc{\alpha_3}{\beta}{\gamma_3})$ is differentiable. So, we need to check whether $s \to H'_s$ is continuous. Note that for all $g \in \VVc{\alpha_0}{ \beta}{\gamma_0}$ and for all $s>0$, $H'_s(g) \in \VVc{\alpha_2}{ \beta}{\gamma_2}$ and for all $h>0$
\begin{align*}
    \|(H'_{s+h}-H'_s)g\|_{\alpha_3,\beta,\gamma_3} &=\|(e^{ih\chi}-1)H'_s(g)\|_{\alpha_3,\beta,\gamma_3}\\ &
     \lesssim \|(H_h-H_0)\One\|_{\alpha^*,\beta,\gamma_0}\|H'_s(g)\|_{\alpha_2,\beta,\gamma_2} \to 0,
\end{align*}
as $h \to 0$ due to part (2). Hence, we have the continuity of the derivative. 
\end{proof}

\subsubsection{Sufficient conditions for \texorpdfstring{\Cref{lem:MultOp}}{Lemma 4.2}}\label{sec:Suff4Lemma}

We limit our scope by providing sufficient conditions for the assumptions in \Cref{lem:MultOp}. 

\begin{lem}\label{rem:RegObs1}
Let $\bar\beta>0$. Suppose $\chi$ is continuous and the right and left derivatives of $\chi$ exist  on $\mathring{I}$. If there exists a constant $b\in [0, 1/\bar\beta)$ such that \begin{align}
  \max\{|\chi'(x+)|, |\chi'(x-)|\}\lesssim x^{-b}(1-x)^{-b}\label{eq: cond 3 cond}
 \end{align}
 then $$|e^{is\chi}|_{0,\bar\beta}<\infty$$ holds for all $s>0$.
\end{lem}

\begin{proof}
We have
\begin{align*}
 |e^{i s\chi}|_{0,\bar\beta}
 &\leq \sup_{\eps\in (0,\eps_0]} \int_0^{1/2}\frac{ \osc(e^{i s\chi}, B_{\eps}(x))}{\eps^{\bar\beta}}\,\mathrm{d}\lambda_I(x)+\sup_{\eps'\in (0,\eps_0]} \int_{1/2}^1\frac{\osc(e^{i s\chi}, B_{\eps'}(x))}{\eps'^{\bar\beta}}\,\mathrm{d}\lambda_I(x).
\end{align*}
We will only estimate the first summand as the estimation of the second follows analogously.
Using the definition $\osc(h,A)=\osc(\Re h,A)+\osc(\Im h,A)$
we note that for any measurable set $A$ we have
\begin{align}\label{eq:oscBnd}
    \osc\left(e^{is\chi}, A\right)\leq \min\{4, 4s/\pi \osc(\chi, A)\}.
\end{align}
By \eqref{eq: cond 3 cond} there exists $C>0$ such that for all $\eps>0$ and all $x\in [\eps,1/2]$ we have
\begin{align*}
 \osc(e^{i s\chi}, B_{\eps}(x))
 &\leq \frac{8|s|\eps}{\pi}\sup_{y\in B_{\eps}(x)}\max\{|\chi'(y+)|, |\chi'(y-)|\}
 \leq \frac{8C|s|\eps}{\pi} (x-\eps)^{-b}.
\end{align*}
We have that $8C|s|\eps/\pi (x-\eps)^{-b}\leq 4$ if and only if
\begin{align*}
 x\leq \left(\frac{2C|s|\eps}{\pi}\right)^{1/b}+\eps=: K_\eps.
\end{align*}
Hence, we split the integral on $[0,1/2]$ into two, one on $[0,K_\eps]$ and the other on $[K_\eps,1/2]$. For the first range, we use the first bound in \eqref{eq:oscBnd} and for the second range, we use the second bound. Then,
\begin{align}
 \sup_{\eps\in (0,\eps_0]} \int_0^{1/2}\frac{ \osc(e^{i s\chi}, B_{\eps}(x))}{\eps^{\bar\beta}}\,\mathrm{d}\lambda_I(x)
 &\leq \sup_{\eps\in (0,\eps_0]} \left( 4{K_\eps}\eps^{-\bar\beta}+\int_{K_\eps}^{1/2} \frac{8C|s|\eps^{1-\bar\beta}}{\pi} (x-\eps)^{-b}\mathrm{d}\lambda_{I}(x)\right)\notag\\
 &\leq \sup_{\eps\in (0,\eps_0]} 4{K_\eps}\eps^{-\bar\beta}+\sup_{\eps\in (0,\eps_0]}\int_{K_\eps}^{1/2} \frac{8C|s|\eps^{1-\bar\beta}}{\pi} (x-\eps)^{-b}\mathrm{d}\lambda_{I}(x).\label{eq: estim osc}
\end{align}
For the first summand, we have
\begin{align*}
 \sup_{\eps\in (0,\eps_0]} 4{K_\eps}\eps^{-\bar\beta}
 &\leq  8\sup_{\eps\in (0,\eps_0]} \max\left\{\left(\frac{2C|s|}{\pi}\right)^{1/b}{\eps}^{1/b-\bar\beta},{\eps}^{1-\bar\beta} \right\}<\infty,
\end{align*}
which follows from the fact that $b<1/\bar\beta$ and $\bar\beta\leq 1$.
For the second summand of \eqref{eq: estim osc}, we have
\begin{align*}
&\sup_{\eps\in (0,\eps_0]}\int_{K_\eps}^{1/2} \frac{8C|s|\eps^{1-\bar\beta}}{\pi} (x-\eps)^{-b}\mathrm{d}\lambda_{I}(x)\\
 &\qquad\leq \frac{8C|s|}{\pi}\sup_{\eps\in (0,\eps_0]} \eps^{1-\bar\beta}\int_{\left(\frac{2Cs\eps}{\pi}\right)^{1/b}}^{1/2}  x^{-b}\mathrm{d}\lambda_{I}(x)\\
  &\qquad\leq \begin{cases}
 \frac{8C|s|}{\pi}\sup_{\eps\in (0,\eps_0]} \frac{\eps^{1-\bar\beta}}{|1-b|} \max\left\{\frac{1}{2},\left(\frac{2C|s|\eps}{\pi}\right)^{1/b} \right\}^{1-b} & b \neq 1\, \\
 \vspace{-5pt} & \\
 \frac{8C|s|}{\pi}\sup_{\eps\in (0,\eps_0]} \eps^{1-\bar\beta} \log \left(\frac{\pi}{2C|s|\eps} \right) & b = 1\,
 \end{cases} \\
  &\qquad= \frac{8C|s|}{\pi}\max\left\{\frac{\eps_0^{1-\bar\beta}}{2^{1-b}|1-b|},\left(\frac{2C|s|}{\pi}\right)^{ 1/b-1}\frac{\eps_0^{1/b-\bar\beta}}{|1-b|}, \eps_0^{1-\bar\beta} \log\left(\frac{\pi}{2C|s|\eps_0} \right)\right\}<\infty\,,
\end{align*}
which again follows from the fact that $\bar\beta\leq 1$ and $b<1/\bar\beta$.
\end{proof}

\begin{rem}
The above lemma combined with \Cref{cor:RegOpr} gives a sufficient condition on $\chi$ for the operator $H_s$, and hence, $\wh\psi_{is}$ to be a bounded linear operator on $V_{\alpha,\beta,\gamma}$ for all $\alpha \geq 0$, $\beta \leq \bar\beta$ and $\gamma \geq 1$. 
\end{rem}

Next, we state a lemma that gives sufficient condition on $\chi$ for the operator valued function $s \mapsto H_s$, and hence, $s \mapsto \wh\psi_{is}$ to be continuous.

\begin{lem}\label{rem:RegObs3}
Suppose  $|\chi|_{\alpha,\beta}<\infty$ with $0 \leq \alpha \leq \beta  < 1/(1+\alpha)$ and there exists $b\in [0,1/\beta)$ such that \eqref{eq: cond 3 cond} holds. Then,  for all $\alpha^*\in(0,1)$ 
$$\lim_{s\to 0}|1-e^{is\chi}|_{\alpha^*,\beta}=0\,.$$
\end{lem}

\begin{proof}
We will do the calculation only for the real part $\Re(1-e^{is\chi})=1-\cos(s\chi)$ and the calculations for the imaginary part $\Im(1-e^{is\chi})= -\sin(s\chi)$ follow analogously and we mention these estimates briefly.  Furthermore, we use the splitting of the positive and negative part as in \eqref{eq: pos neg part}. Also,  since $\Re(1-e^{is\chi})_{-}=0,$ it does not contribute to the estimates. 

For $\delta\in (0,\eps_0)$ to be specified later depending on $\eps$ and $s$, we have
\begin{align}
  |\Re(1-e^{is\chi})_{+}|_{\alpha^*,\beta}
  &=\sup_{\eps\leq \eps_0}\frac{\int \osc(R_{\alpha^*}\Re(1-e^{is\chi})_{+}, B_{\eps}(x))\,\mathrm{d}\lambda_I(x)}{\eps^{\beta}}\notag\\
   &\leq \sup_{\eps\leq \eps_0}\frac{\int \osc(R_{\alpha^*}\Re(1-e^{is\chi})_{+}\One_{[0,\delta+\eps]}, B_{\eps}(x))\,\mathrm{d}\lambda_I(x)}{\eps^{\beta}}\label{eq: 1-e alphabeta 1}\\
   &\qquad + \sup_{\eps\leq \eps_0}\frac{\int \osc(R_{\alpha^*}\Re(1-e^{is\chi})_{+}\One_{(\delta+\eps,1-\delta-\eps)}, B_{\eps}(x))\,\mathrm{d}\lambda_I(x)}{\eps^{\beta}}\label{eq: 1-e alphabeta 2}\\
   &\qquad\qquad+\sup_{\eps\leq \eps_0}\frac{\int \osc(R_{\alpha^*}\Re(1-e^{is\chi})_{+}\One_{[1-\delta-\eps,1]}, B_{\eps}(x))\,\mathrm{d}\lambda_I(x)}{\eps^{\beta}},\label{eq: 1-e alphabeta 3}
 \end{align}
 where we assume that $s$ and $\eps_0$ are so small that $\delta+\eps<1-\delta-\eps$.

 We start by estimating the middle summand \eqref{eq: 1-e alphabeta 2}.  \cite[Proposition 3.2(ii)]{BS} yields
 \begin{align}
  \MoveEqLeft\osc(R_{\alpha^*}\Re(1-e^{is\chi})_{+}\One_{(\delta+\eps,1-\delta-\eps)},B_{\eps}(x))\notag\\
  &\leq \osc(R_{\alpha^*}\Re(1-e^{is\chi})_{+}, (\delta+\eps,1-\delta-\eps)\cap B_{\eps}(x))\One_{(\delta+\eps,1-\delta-\eps)}(x) \label{eq: osc with indicator function}
  \\ &\qquad+2 \left[\esssup_{(\delta+\eps,1-\delta-\eps)\cap B_{\eps}(x)}R_{\alpha^*}\Re(1-e^{is\chi})_{+}\right]\One_{B_{\eps}((\delta+\eps,1-\delta-\eps))\cap B_{\eps}((\delta+\eps,1-\delta-\eps)^c)}(x)\,.\notag
 \end{align}
 We first investigate the first summand of \eqref{eq: osc with indicator function}. For the following, we set 
 \begin{equation}
  D(\delta,\eps,x):=(\delta+\eps,1-\delta-\eps)\cap B_{\eps}(x).\label{eq: def D}
 \end{equation}
For $x\in (\delta+\eps, 1-\delta-\eps)$, 
 \begin{align}
  &\osc( R_{\alpha^*}(1-\cos(s\chi)),D(\delta,\eps,x)) \leq 2 \eps \sup_{D(\delta,\eps,x)} [R_{\alpha^*}(1-\cos(s\chi))]'\notag\\
  &\qquad\leq 
  2 \eps \left[\sup_{D(\delta,\eps,x)} |(R_{\alpha^*}\One)'|\,(1-\cos(s\chi))+\sup_{D(\delta,\eps,x)} (R_{\alpha^*}\One)\, |(1-\cos(s\chi))'|\right].
  \label{eq: middle term}
 \end{align}
Both of the above calculations follow analogously for the imaginary part with $|\sin (s\chi)|$ instead of $1-\cos(s\chi)$.

 We set $\delta=\delta(\eps,s)= \eps^{\kappa}\cdot|s|^\iota$ 
  with $\kappa \in (0,1)$ and $\iota >0$ to be specified later. 
  Since $|\chi|_{\alpha,\beta}<\infty$ implies that $R_{\alpha}\chi$ is essentially bounded, we can conclude that there exists $K(\chi)\in(0,\infty)$ such that
  $|\chi(x)|\leq K(\chi)\cdot x^{-\alpha}(1-x)^{-\alpha}$ almost everywhere. 
   Recall that there is $C>0$ such that $\max\{|1- \cos (x)\}|, |\sin (x)|\} \leq C|x|\,.$ Combining this with $(R_{\alpha^*}\One)' = \alpha^*( x^{\alpha^*-1}(1-x)^{\alpha^*}+ x^{\alpha^{*}}(1-x)^{\alpha^*-1})\,,$ we have
\begin{equation}\label{eq: 1-cos estim}
 \sup_{D(\delta,\eps,x)} |(R_{\alpha^*}\One)'|\max\{1-\cos(s\chi)\,,|\sin (s\chi)|\}\, \lesssim  \frac{|s|}{ (x-\eps)^{1+\alpha - \alpha^*}} \,,
\end{equation}
when $x \leq 1/2$. The estimates for $x\geq 1/2$ follows from replacing $(x-\eps)$ by $(1-x+\eps)$, and the final estimates remain unchanged. So, we restrict our attention to the former case. 

It follows that
  \begin{align}
  \MoveEqLeft\lim_{s\to0}\sup_{\eps\in (0,\eps_0]}\frac{1}{\eps^\beta} \int_{\delta+\eps}^{1/2} 2\eps\sup_{D(\delta,\eps,x)}  \Big(|(R_{\alpha^*}\One)'|(1-\cos(s\chi)) \Big)\One_{(\delta,1-\delta)}(x)\,\mathrm{d}\lambda_I(x)\notag\\
 &\lesssim \lim_{s\to0}|s|\sup_{\eps\in (0,\eps_0]}\eps^{1-\beta} 
 \int_{\delta+\eps}^{1/2} (x-\eps)^{\alpha^*-1-\alpha}\, \mathrm{d}\lambda_I(x) \notag\\
 &\lesssim \lim_{s\to0}|s|\sup_{\eps\in (0,\eps_0]}\eps^{1-\beta} 
 \int_{\delta}^{1/2-\eps} x^{\alpha^*-1-\alpha}\, \mathrm{d}\lambda_I(x) \notag\\
 &\lesssim \begin{cases}
 \lim_{s\to0} \eps_0^{1-\beta+\kappa(\alpha^*-\alpha)}\lim_{s\to0}|s|^{1+\iota\left(\alpha^*-\alpha\right)} 
 & \alpha^* < \alpha \\
 \eps_0^{1-\beta}(|\log(1/2-\eps_0)|+ \kappa |\log(\eps_0)|) \lim_{s\to0} |s| + \iota\eps_0^{1-\beta}\lim_{s\to0}  |s|\, |\log|s||  &\alpha^*=\alpha\\
 \eps_0^{1-\beta} \lim_{s\to 0}|s| & \alpha^*>\alpha
 \end{cases}\notag\\
 &=0\,
\label{eq: middle term 1}
 \end{align}
provided that under the condition $\alpha^*<\alpha$ we have
\begin{equation}\label{eq:iotakappa1}
\begin{aligned}
1-\beta+\kappa(\alpha^*-\alpha)>0 &\iff\kappa < (1-\beta)/(\alpha- \alpha^*)\,, \\
\iota\left(\alpha^*-\alpha\right)+1>0 &\iff \iota < 1/(\alpha - \alpha^*)\,.
\end{aligned} 
\end{equation}
 Analogously, under the same conditions,
\begin{align*}
 \lim_{s\to0}\sup_{\eps\in (0,\eps_0]}\frac{1
 }{\eps^\beta}\int 2\eps
 \sup_{D(\delta,\eps,x)} \Big(|(R_{\alpha^*}\One)'\,|(\sin (s\chi))_{\pm}\Big)\One_{(\delta+\eps,1-\delta-\eps)}(x)\,\mathrm{d}\lambda_I(x) = 0.
\end{align*}

 To estimate the second summand of \eqref{eq: middle term}, 
    we use $(1-\cos(s\chi))' = \sin(s \chi) \cdot s \chi'\,,$ $(\sin(s\chi))' = \cos(s \chi) \cdot s \chi'\,,$ $|\cos(s\chi)|\leq 1$, $|\sin (s\chi)| \leq 1 $, and our assumption about $\chi'$. Then, we have 
 \begin{align*}
\sup_{D(\delta,\eps,x)} \max\left\{(R_{\alpha^*}\One) \,|(1-\cos(s\chi))'|\,,(R_{\alpha^*}\One) |(\sin(s\chi)^\pm)'|\right\}
    &  \lesssim 
    \begin{cases}
    |s|(x-\eps)^{\alpha^*-b} & \alpha^*<b \\
    |s| \cdot 1 & \alpha^*\geq b
    \end{cases}
\end{align*}
for $x \leq 1/2$. Also, note that for $x \leq 1/2$ and the estimate for $x \geq 1/2$ is the same with $(x-\eps)$ replaced by $(1-x+\eps)$. 
 Thus, if $\alpha^*< b$
 \begin{align}\label{eq: middle term 1-1}
  \MoveEqLeft\lim_{s\to0}\sup_{\eps\in (0,\eps_0]}\frac{1}{\eps^\beta}\int_{\delta+\eps}^{1/2}   2\eps
   \sup_{D(\delta,\eps,x)} \Big((R_{\alpha^*}\One ) \, |(1-\cos(s\chi))'|\Big)\One_{(\delta+\eps,1-\delta-\eps)}(x)\,\mathrm{d}\lambda_I(x)
 \notag\\ &\lesssim \lim_{s\to 0}\sup_{\eps\in (0,\eps_0]} \eps^{1-\beta
 }|s|\int_{\delta}^{1/2-\eps} {x^{\alpha^*-b}} \mathrm d\lambda_I(x)
 \notag \\
 &\lesssim \begin{cases}
 \eps_0^{1-\beta+\kappa(1+\alpha^*-b)}\lim_{s\to0} |s|^{1+\iota(1+\alpha^*-b)} & b > 1+\alpha^* \\ 
 \eps_0^{1-\beta}(|\log(1/2-\eps_0)|+\kappa|\log(\eps_0)|)\,\lim_{s\to 0}|s|+ \iota \eps_0^{1-\beta}\lim_{s\to 0} |s|\, |\log |s||&  b=1+\alpha^*\\
 \eps_0^{1-\beta}\,\lim_{s\to0} |s| & b < 1+\alpha^*
 \end{cases}\notag\\
  &= 0\,,\end{align}
where, in the case of $b>1+\alpha^*$, we have assumed that \begin{equation}\label{eq:iotakappa2}  \begin{aligned} 1-\beta+\kappa(1+\alpha^*-b)>0 &\iff\kappa <(1-\beta)/(b-1-\alpha^*)\,, \\ 1+\iota(1+\alpha^*-b)>0 &\iff \iota < 1/(b-1-\alpha^*)\,. \end{aligned}  \end{equation}
The $\alpha^* \geq b$ case is similar to the $b < 1+ \alpha^*$ case above. 
 Analogously, under the same assumptions on $\kappa$ and $\iota$, 
 we obtain
  \begin{align*}
  \MoveEqLeft\lim_{s\to0}\sup_{\eps\in (0,\eps_0]}\frac{1
 }{\eps^\beta} \int  2\eps \sup_{D(\delta,\eps,x)} \Big(R_{\alpha^*}\One\, \left|\left(\sin(s\chi)_{\pm}\right)'\right|\Big) 
 \One_{(\delta+\eps,1-\delta-\eps)}(x)\,\mathrm{d}\lambda_I(x)
  = 0\,.
 \end{align*}
Hence, combining \eqref{eq: middle term 1} and \eqref{eq: middle term 1-1}, we can conclude
 \begin{align}
  \lim_{s\to0}\sup_{\eps\in (0,\eps_0]}\frac{1
 }{\eps^\beta}\int \osc\left(R_{\alpha^*}\Re (1-e^{is\chi})_{ +}, D(\delta,\eps,x)\right)\One_{(\delta,1-\delta)}(x)\,\mathrm{d}\lambda_I(x)=0.\label{eq: estim middle term}
\end{align}
Also, the analogous result for the imaginary part, $\Im (1-e^{is\chi})_{\pm},$ follows.

Next, we will estimate the second summand in \eqref{eq: osc with indicator function}.
We note that
\begin{align*}
  B_{\eps}((\delta+\eps,1-\delta-\eps))\cap B_{\eps}((\delta+\eps,1-\delta-\eps)^c)=B_{\eps}(\delta+\eps)\cup B_{\eps}(1-\delta-\eps)\, 
\end{align*}
and hence,
\begin{align}\label{eq:InterEpsBall}
    \One_{ B_{\eps}((\delta+\eps,1-\delta-\eps))\cap B_{\eps}((\delta+\eps,1-\delta-\eps)^c)} = \One_{ B_{\eps}(\delta+\eps)\cup B_{\eps}(1-\delta-\eps)}.
\end{align}
It follows that  
\begin{align}
\sup_{D(\delta,\eps,x)}R_{\alpha^*}(1-\cos(s\chi))
\lesssim \begin{cases}
    |s|(x-\eps)^{\alpha^*-\alpha}&\alpha^*<\alpha\\
    |s|(x+\eps)^{\alpha^*-\alpha}&\alpha^*\geq \alpha.
\end{cases} \label{eq: 1-cos estim a}   
\end{align}
 Due to the symmetry around $x=1/2$, we obtain
  \begin{align}
  \MoveEqLeft\lim_{s\to0}\sup_{\eps\in (0,\eps_0]}\frac{1
 }{\eps^\beta}  \int  \sup_{D(\delta,\eps,x)}R_{\alpha^*}\Re(1-e^{is\chi})_{ +}\cdot\One_{ B_{\eps}((\delta+\eps,1-\delta-\eps))\cap B_{\eps}((\delta+\eps,1-\delta-\eps)^c)}{(x)}\,\mathrm{d}\lambda_I(x) \notag\\
 & = \lim_{s\to0}\sup_{\eps\in (0,\eps_0]}\frac{
 1}{\eps^\beta}  \bigg( \int_{\delta}^{\delta+2\eps} +\int_{1-\delta-2\eps}^{1-\delta}\bigg)\sup_{D(\delta,\eps,x)}R_{\alpha^*}(1-\cos(s\chi))\,\mathrm{d}\lambda_I(x)\notag
 \notag\\
 &\lesssim \lim_{s\to0}\sup_{\eps\in (0,\eps_0]} |s|\eps^{-\beta} \int_{\delta}^{\delta+2\eps}  \max\{(\delta+\eps)^{\alpha^*-\alpha}, (\delta+2\eps)^{\alpha^*-\alpha}\} \,\mathrm{d}\lambda_I(x)\notag\\
 &\lesssim\begin{cases}
 \lim_{s\to0} \eps_0^{1-\beta-\kappa(\alpha-\alpha^*)}|s|^{1-\iota(\alpha-\alpha^*)} & \alpha^* < \alpha\\ 
 \lim_{s\to0} \eps_0^{1-\beta}|s| & \alpha^* \geq \alpha
 \end{cases} \nonumber\\
 &=0 \label{eq: middle term 2}
\end{align}
 where, in the case of $\alpha^* < \alpha$, we assume that
\begin{equation}\label{eq:iotakappa3}
   \begin{aligned}
        1-\beta-\kappa(\alpha-\alpha^*)>0 &\iff \kappa < (1-\beta)/(\alpha-\alpha^*) \,,\\
        1-\iota(\alpha-\alpha^*)>0 &\iff \iota < 1/(\alpha-\alpha^*)\,.
    \end{aligned}
\end{equation}
Combining this with \eqref{eq: osc with indicator function} and \eqref{eq: estim middle term} yields that
the summand \eqref{eq: 1-e alphabeta 2} tends to zero for $s\to0$ and the same is true for the imaginary part, $\Im (1-e^{is\chi})_\pm,$ because the same assumptions on $\kappa $ and $\iota$ along with $|\sin(x)|\lesssim |x|$ and \eqref{eq: 1-cos estim a} yield  
$$
\lim_{s\to0}\sup_{\eps\in (0,\eps_0]}\frac{1
 }{\eps^\beta}\int \osc\left( (\sin(s\chi))_{\pm}, D(\delta,\eps,x)\right) \One_{(\delta+\eps,1-\delta-\eps)}(x)\,\mathrm{d}\lambda_I(x)=0\,,$$
$$\lim_{s\to0}\sup_{\eps\in (0,\eps_0]}\frac{1
 }{\eps^\beta} \int  \sup_{D(\delta,\eps,x)}R_{\alpha^*}(\sin(s\chi))_{\pm}\cdot\One_{B_{\eps}((\delta+\eps,1-\delta-\eps))\cap B_{\eps}((\delta+\eps,1-\delta-\eps)^c)}\,\mathrm{d}\lambda_I(x) = 0\,.$$

Finally, we investigate into the first summand \eqref{eq: 1-e alphabeta 1}. As the calculation for the summand \eqref{eq: 1-e alphabeta 3} is very similar, we will only give the details for \eqref{eq: 1-e alphabeta 1}.  We split the integral into
\begin{align}
  \MoveEqLeft\lim_{s\to 0}\sup_{\eps\leq \eps_0}\frac{1}{\eps^{\beta}}\int \osc(R_{\alpha^*}\Re(1-e^{is\chi})\One_{[0,{\delta+\eps}]},B_{\eps}(x))\,\mathrm{d}\lambda_I(x)\notag\\
  &=\lim_{s\to 0}\sup_{\eps\leq \eps_0}\frac{1}{\eps^{\beta}}\left(\int_{{[0,\delta)}}+ \int_{[{\delta,\delta+2\eps}]}\right)\osc(R_{\alpha^*}(1-\cos(s\chi))\One_{[0,{\delta+\eps}]},B_{\eps}(x))\,\mathrm{d}\lambda_I(x).\label{eq: estim near zero0}
\end{align} 
 For the first summand of \eqref{eq: estim near zero0}, 
 we write 
 \begin{equation}\label{eq: bar D}
  \bar D(\delta,\eps,x):=[0,\delta+\eps]\cap B_{\eps}(x)   
 \end{equation} and
 we note that $\Re(1-e^{is\chi})\in (0,2)$ and
 \begin{align}
  \osc(R_{\alpha^*}{ (1-\cos(s\chi))}\One_{[0,\delta+\eps]}, B_{\eps}(x))
 &\leq 2 \cdot {\sup}_{\bar D(\delta,\eps,x)}R_{\alpha^*}\One
  \leq 2R_{\alpha^*}\One(x+\eps)
 \leq 2 (x+\eps)^{\alpha^*}.\label{eq: estim near zero}
 \end{align}
Now, we have
 \begin{align}\label{eq: estim near zero1} 
   \MoveEqLeft\lim_{s\to 0}\sup_{\eps\leq \eps_0}\frac{1}{\eps^{\beta}}\int_{0}^{\delta} \osc(R_{\alpha^*}{(1-\cos(s\chi))}\One_{[0,\delta+\eps]}, B_\eps(x))\,\mathrm{d}\lambda_I(x)=0
\end{align}
under the condition 
\begin{equation}\label{eq: iotakappa3}
\iota>0\quad\text{ and }\quad \kappa>\frac{\beta}{1+\alpha^*}\end{equation} 
due to the following Sub-lemma. 

\begin{sublem}\label{lem:KeyEst1}
 Define 
\begin{align*}
    \Theta_1 &:= \lim_{s\to 0}\sup_{\eps\leq \eps_0}\frac{1}{\eps^{\beta}}\int_0^\delta \osc(R_{\alpha^*}\Re(1-e^{is\chi})_{+}, \bar D(\delta,\eps,x))
    \,\mathrm{d}\lambda_I(x)
\end{align*}
and 
\begin{align*}
   \Theta_2 &:= \lim_{s\to 0}\sup_{\eps\leq \eps_0}\frac{1}{s\eps^{\beta}}\int_0^\delta\osc(R_{\alpha^*}\Re(1-e^{is\chi})_{+}, \bar D(\delta,\eps,x))
   \,\mathrm{d}\lambda_I(x)
\end{align*}
with $\alpha^*\geq 0$, $\delta =\delta(\eps,s) = \eps^\kappa s^\iota$ where $\iota,\kappa>0$ and with $\bar D$ as in \eqref{eq: bar D}. Suppose $|\chi|_{\alpha, \beta}<\infty$ with $0\leq \alpha<\beta\leq 1$.
   \begin{enumerate}[leftmargin=*]
       \item[$(1)$] If
\begin{equation}\label{eq:kappaiotaC0}
    \iota >0\quad\text{and}\quad  \kappa \geq \frac{\beta}{\alpha^*+1},
\end{equation}
then $\Theta_1=0\,.$
\item[$(2)$] If 
\begin{equation}\label{eq:kappaiotaC2}
    \iota >1\quad\text{and}\quad \kappa \geq \frac{\beta}{\alpha^*+1},
\end{equation}
then $\Theta_2=0.$
   \end{enumerate}
\end{sublem}

\begin{proof}[Proof of Sub-Lemma \ref{lem:KeyEst1}]
Without loss of generality, we assume that $s>0$. Note that due to \eqref{eq: estim near zero}, we have
\begin{align*}
  \Theta_1 &\leq \lim_{s\to 0}\sup_{\eps\leq \eps_0}\frac{\int_0^{\delta} 2 (x+\eps)^{\alpha^*}\,\mathrm{d}\lambda_I(x)}{\eps^{\beta}}\notag\\
  &= \lim_{s\to 0}\sup_{\eps\leq \eps_0}\frac{2 \left(\left({ \delta}+\eps\right)^{\alpha^*+1}- \eps^{\alpha^*+1}\right)}{(\alpha^*+1)\eps^{\beta}}= \lim_{s\to 0} L(s,\eps_0)
\end{align*}
where 
\begin{align*}
    L(s,\eps_0)&:=\frac{2}{(\alpha^*+1)}\sup_{\eps\leq \eps_0} J(s,\eps)\,,
\end{align*}
and 
\begin{align*}
    J(s,\eps) :=\frac{ \left({ \delta}+\eps\right)^{\alpha^*+1}- \eps^{\alpha^*+1}}{\eps^{\beta}} = 
 (\eps^\kappa s^\iota+\eps)^{\alpha^*+1}\eps^{-\beta}-\eps^{\alpha^*+1-\beta}\,.
\end{align*}

First, we note that
\begin{align}
    (\eps^\kappa s^\iota+\eps)^{\alpha^*+1}\eps^{-\beta}= (\eps^{\kappa-\beta/(\alpha^*+1)}s^{\iota}+\eps^{1-\beta/(\alpha^*+1)})^{\alpha^*+1}\,,\label{eq: J(s eps)}
\end{align}
and hence, for $J(s,\eps)$ to not blow up near $\eps=0$, we should have \eqref{eq:kappaiotaC0}. Due to the first inequality in \eqref{eq:kappaiotaC0} and \eqref{eq: J(s eps)}, we have for given $s>0$ that $\sup_{\eps\leq \eps_0} J(s,\eps)= J(s, \eps_0)$ and thus $J(0,\eps):=\lim_{s\to 0}J(s,\eps)=0$ for all $\eps$. Therefore, under the assumption \eqref{eq:kappaiotaC0}, $\Theta_1 = 0$ as claimed because
$$ \Theta_1 \leq \lim_{s\to 0} L(s,\eps_0) = \frac{2}{\alpha^*+1} \lim_{s\to 0} J(s,\eps_0) = 0\,.$$

Now, using \eqref{eq: estim near zero} and l'Hôpital's rule, we obtain
\begin{align*}
   \Theta_2 &= \lim_{s\to 0}\sup_{\eps\leq \eps_0}\frac{\int \osc(R_{\alpha}\Re(1-e^{is\chi})_{+}, \bar D(\delta,\eps,x))\One_{[0,\delta]}(x)\,\mathrm{d}\lambda_I(x)}{s\eps^{\beta}}\notag\\
  &\leq \lim_{s\to 0}\frac{1}{s}\sup_{\eps\leq \eps_0}\frac{2 \left(\left({ \delta}+\eps\right)^{\alpha^*+1}- \eps^{\alpha^*+1}\right)}{(\alpha^*+1)\eps^{\beta}}= \frac{\mathrm{d}}{\mathrm{d}s} L(s,\eps_0) \Big|_{s=0}. 
\end{align*}

 We note that the last equality follows by the above calculation, namely that $\sup_{\eps\leq \eps_0} J(s,\eps)= J(s, \eps_0)$ holds because of $1+\alpha^*>\beta$ and the additional conditions $\iota>0$ and $\kappa\geq \beta/(\alpha^*+1)$.

Next, taking the derivative of $L$ wrt $s$, we obtain
\begin{align*}
   \frac{\mathrm d}{\mathrm ds} L(s,\eps_0) = \frac{2}{\alpha^*+1}\frac{ \mathrm d}{\mathrm ds}(\eps_0^\kappa s^\iota+\eps_0)^{\alpha^*+1}\eps_0^{-\beta} = 2\iota(\eps^\kappa_0s^\iota+\eps_0)^{\alpha^*}\eps^{-\beta}_0\eps^{\kappa}_0 s^{\iota-1}.
\end{align*}
Note that for $\Theta_2 = 0$ we should have $\frac{\mathrm d}{\mathrm ds} L(s,\eps_0)\big|_{s=0} = 0$ and this is true, if
$\iota > 1$.
Therefore, under \eqref{eq:kappaiotaC2}, we have that $\Theta_2 = 0$ as claimed.
\end{proof}

Next, in order to estimate the second summand of \eqref{eq: estim near zero0}, 
we first note that for $x \in [\delta,\delta+2\eps],$
\begin{equation}\label{eq:Near0SupEst}
    \sup_{\bar D(\delta,\eps,x)} R_{\alpha^*}(1-\cos (s\chi))\One_{[\delta,\delta+\eps]}\lesssim |s|\sup_{y \in \bar D(\delta,\eps,x)\cap [\delta, \delta+2\eps]} y^{\alpha^*-\alpha} \leq \begin{cases}
 |s| & \alpha^* \geq  \alpha \\
 |s|\cdot \delta^{\alpha^*-\alpha} & \alpha^* < \alpha\,. \end{cases}
\end{equation}
Hence,
\begin{align}
  \MoveEqLeft
  \lim_{s\to 0}\sup_{\eps\leq \eps_0}\frac{1}{\eps^{\beta}} \int_{[{\delta,\delta+2\eps}]}\osc\left(R_{\alpha^*}(1-\cos(s\chi))\One_{[\delta,\delta+\eps]},B_{\eps}(x)\right)\,\mathrm{d}\lambda_I(x)\notag\\
  &\lesssim\lim_{s\to 0}\sup_{\eps\leq \eps_0}\frac{1}{\eps^{\beta}} \int_{[\delta,\delta+2\eps]}  \sup_{\bar D(\delta,\eps,x)} \left(R_{\alpha^*}(1-\cos (s\chi))\One_{[\delta,\delta+\eps]}\right)\,\mathrm{d}\lambda_I(x) \notag\\
  &\lesssim  \begin{cases}\eps_0^{1-\beta}\lim_{s\to 0}  |s|
  &\alpha^* \geq  \alpha\\
  \eps_0^{1-\beta+\kappa(\alpha^*-\alpha)}\lim_{s\to 0}|s|^{1+\iota(\alpha-\alpha^*)}
  &\alpha^*< \alpha\end{cases}\notag\\
  &=0
\label{eq: cond 1 last estim}
 \end{align}
 provided that, in the case of $\alpha^*< \alpha ,$ 
\begin{equation}\label{eq: iotakappa4}
  \begin{aligned}
1-\beta+\kappa\left(\alpha^*- \alpha\right)>0 &\iff\kappa < (1-\beta)/(\alpha-\alpha^*)\,, \\
1+\iota(\alpha^*- \alpha)>0&\iff \iota< 1/(\alpha-\alpha^*).
\end{aligned}  
\end{equation}
Next, by \cite[Prop.~3.2(ii)]{BS} we have for $x\in (\delta, \delta+ 2\eps]$
\begin{align*}
    &\osc\left(R_{\alpha^*}(1-\cos(s\chi))\One_{[0,\delta]},B_{\eps}(x)\right)\\ 
    &\qquad\leq \osc\left(R_{\alpha^*}(1-\cos(s\chi)),B_{\eps}(x)\cap [0,\delta]\right)\One_{[0,\delta]}(x)
    +2\esssup_{B_{\eps}(x)\cap [0,\delta]} R_{\alpha^*}(1-\cos(s\chi))\One_{ [\delta-\eps \vee 0, \delta+\eps]}(x)\\
    &\qquad \leq 0+2\esssup_{[\delta-\eps \vee 0,\delta]} R_{\alpha^*}\One
     \leq 2 \delta^{\alpha^*}.
\end{align*}
Hence, 
\begin{align*}
    \MoveEqLeft
  \lim_{s\to 0}\sup_{\eps\leq \eps_0}\frac{1}{\eps^{\beta}} \int_{[{\delta,\delta+2\eps}]}\osc\left(R_{\alpha^*}(1-\cos(s\chi))\One_{[0,\delta]},B_{\eps}(x)\right)\,\mathrm{d}\lambda_I(x)\notag\\  
  &\lesssim \lim_{s\to 0}\sup_{\eps\leq \eps_0}\frac{1}{\eps^{\beta}} \int_{[{\delta,\delta+2\eps}]}\delta^{\alpha^*}\,\mathrm{d}\lambda_I(x)
  \lesssim \lim_{s\to 0}\sup_{\eps\leq \eps_0} \eps^{1-\beta+\kappa\alpha^*}s^{\iota\alpha^*}=0
\end{align*}
under \eqref{eq: iotakappa3}. This together with \eqref{eq: cond 1 last estim} imply
\begin{align*}
  \lim_{s\to 0}\sup_{\eps\leq \eps_0}\frac{1}{\eps^{\beta}} \int_{[{\delta,\delta+2\eps}]}\osc\left(R_{\alpha^*}(1-\cos(s\chi))\One_{[0,\delta+2\eps]},B_{\eps}(x)\right)\,\mathrm{d}\lambda_I(x)=0.
\end{align*}
Combining this with \eqref{eq: estim near zero0} and \eqref{eq: estim near zero1} implies that \eqref{eq: 1-e alphabeta 1} tends to zero for $s$ tending to zero. The same is true for the imaginary part, $\Im (1-e^{is\chi})_{\pm},$ as $\Im (1-e^{is\chi})_{\pm} \leq 1.$ 

 Finally, we discuss here possible values of $\alpha^*$ and the implicit requirements  on $b$ that ensure the existence of $\iota>0$ and $\kappa>0$ used in the proof. There are four cases.

Note that, in the case of $\alpha^*< \alpha$ and $b>1+\alpha^*\,,$ under \eqref{eq:iotakappa1}, \eqref{eq:iotakappa2}, \eqref{eq:iotakappa3}, \eqref{eq: iotakappa3} and \eqref{eq: iotakappa4}, we have
\begin{equation*}
    \begin{aligned}
        \frac{\beta}{\alpha^*+1} <\, &\kappa < \min\left\{{\frac{1-\beta}{\alpha-\alpha^*}, \frac{1-\beta}{b-1-\alpha^*}, 1}\right\}, \\
         0 <\, &\iota < \min\left\{{\frac{1}{\alpha-\alpha^*}, \frac{1}{b-1-\alpha^*}}\right\}.
    \end{aligned}
\end{equation*}
First, we see that the conditions on $\iota$ are always fulfilled, because
$$0<\frac{1}{\alpha-\alpha^*} \quad\text{and}\quad 0<\frac{1}{b-1-\alpha^*}\,.$$
Similarly, considering the inequalities that guarantee the existence of $\kappa$, we have  
$$\alpha> \alpha^* > \max\left\{\alpha\beta+\beta-1, \beta b -1\right\}$$
is necessary and sufficient.
Note that due to $\beta <\min\{ 1/b,1/(\alpha +1)\}$ we have 
 $\alpha\beta+\beta-1 < 0\,,$ and also, 
${\beta b -1}< 0.$ So, $0 < \alpha^* <\min\{\alpha, b-1\}$  which  is equivalent to $\alpha^*<\alpha$ and $b>1+\alpha^*$.

In the case of $ \alpha^* < \alpha$ and $b \leq 1+\alpha^*\,,$ \eqref{eq:iotakappa2} poses no restrictions. So, under \eqref{eq:iotakappa1}, \eqref{eq:iotakappa3}, \eqref{eq: iotakappa3} and \eqref{eq: iotakappa4} we have $b-1 < \alpha^* < \alpha$ and $b \leq 1+\alpha$  which is equivalent to our assumptions $\alpha^* < \alpha$ and $b \leq 1+\alpha^*$. 

In the case of $ \alpha^* \geq  \alpha\,,$ \eqref{eq:iotakappa1}, \eqref{eq:iotakappa3}, and \eqref{eq: iotakappa4} pose no restrictions. So, when $b<1+\alpha^*\,,$ we have $ \alpha^* > \max\{\alpha\,, b-1 \}$ and when 
$b>1+\alpha^*$, we have $\alpha < \alpha^* < b-1$ and $b>1+\alpha$ and we do not  obtain any additional restrictions either.
\end{proof}

The next lemma of this section gives a sufficient condition on $\chi$ for the operator valued function $s \mapsto H_s$, and hence, $s \mapsto \wh\psi_{is}$ to be differentiable.

\begin{lem}\label{rem:RegObs4}
Suppose  $|\chi|_{\alpha,\beta}<\infty$ with $0\leq \alpha \leq \beta < 1/(1+\alpha)$ and there exists $b\in [0,1/\beta)$ such that \eqref{eq: cond 3 cond} holds. 
Then, for all $\alpha^* > \min\{2\alpha, \max\{\alpha, \alpha+b-2\}\}$   
we have
$$\lim_{s\to 0}\left|\frac{e^{is\chi}-1-is\chi}{s}\right|_{\alpha^*, \beta}=0\,.$$
\end{lem}

\begin{proof}
The proof follows very similar to the proof of the previous lemma and we will stick to the same notation. 
Again, we will do the calculations only for the  non-negative real part, only noting some differences for the imaginary part.
We have
\begin{align*}
 \osc\left(R_{\alpha^*} \Re \left(\frac{e^{is\chi}-1-is\chi}{s}\right),B_{\eps}(x)\right)
 &=\frac{1}{s}\osc\left(R_{\alpha^*}(1-\cos(s\chi)),B_{\eps}(x)\right)
\end{align*}
and
\begin{align*}
 \osc\left(R_{\alpha^*}\Im \left(\frac{e^{is\chi}-1-is\chi}{s}\right),B_{\eps}(x)\right)
 &=\frac{1}{s}\osc\left(R_{\alpha^*}(\sin(s\chi)-s\chi),B_{\eps}(x)\right).
\end{align*}
As in \eqref{eq: 1-e alphabeta 1} to \eqref{eq: 1-e alphabeta 3}, we have for $\delta\in (0,\eps_0)$ (to be specified later and depending on $s$ and $\eps$) that
 \begin{align}
  \left|\Re(e^{is\chi}-1-is\chi)\right|_{\alpha^*,\beta}
  &=\sup_{\eps\leq \eps_0}\frac{\int \osc(R_{\alpha^*}(1-\cos(s\chi)),B_{\eps}(x))\,\mathrm{d}\lambda_I(x)}{s\eps^{\beta}}\notag\\
   &\leq \sup_{\eps\leq \eps_0}\frac{\int \osc(R_{\alpha^*}(1-\cos(s\chi))\One_{[0,\delta+\eps]}, B_{\eps}(x))\,\mathrm{d}\lambda_I(x)}{s\eps^{\beta}}\label{eq: 1-e alphabeta 1a}\\
   &\qquad   + \sup_{\eps\leq \eps_0}\frac{\int \osc(R_{\alpha^*}(1-\cos(s\chi))\One_{(\delta+\eps,1-\delta-\eps)},B_{\eps}(x))\,\mathrm{d}\lambda_I(x)}{s\eps^{\beta}}\label{eq: 1-e alphabeta 2a}\\
   &\qquad\qquad+\sup_{\eps\leq \eps_0}\frac{\int \osc(R_{\alpha^*}(1-\cos(s\chi))\One_{[1-\delta-\eps,1]},B_{\eps}(x))\,\mathrm{d}\lambda_I(x)}{s\eps^{\beta}}\label{eq: 1-e alphabeta 3a}\,,
 \end{align}
and similarly, for the imaginary part. 

Now, we start by estimating the middle term \eqref{eq: 1-e alphabeta 2a}, and as in \eqref{eq: osc with indicator function}, we use \cite[Proposition 3.2(ii)]{BS} to obtain
 \begin{align}
  \MoveEqLeft\osc(R_{\alpha^*}(1-\cos(s\chi))\One_{(\delta+\eps,1-\delta-\eps)},B_{\eps}(x))\notag\\
  &\leq \osc(R_{\alpha^*}(1-\cos(s\chi)), D(\delta,\eps,x))\One_{(\delta+\eps,1-\delta-\eps)}(x)\label{eq: osc with indicator function deriv}
  \\&\qquad+2 \left[\sup_{D(\delta,\eps,x)}R_{\alpha^*}(1-\cos(s\chi))\right]\One_{B_{\eps}((\delta+\eps,1-\delta-\eps))\cap B_{\eps}((\delta+\eps,1-\delta-\eps)^c)}(x)\,.\notag
 \end{align}
For $x\in (\delta+\eps, 1-\delta-\eps)$, 
 \begin{align}
  &\osc( R_{\alpha^*}(1-\cos(s\chi)),D(\delta,\eps,x)) \leq 2 \eps \sup_{D(\delta,\eps,x)} |[R_{\alpha^*}(1-\cos(s\chi))]'|\notag\\
  &\qquad\leq 2 \eps \left[\sup_{D(\delta,\eps,x)} |(R_{\alpha^*}\One)'|\,(1-\cos(s\chi))+\sup_{D(\delta,\eps,x)} (R_{\alpha^*}\One) |(1-\cos(s\chi))'|\right].
  \label{eq: osc Ralpha cos} 
 \end{align}
Both of the above calculations follow analogously for the imaginary part.

For the following, as in the previous proof, we set  $\delta=\delta(\eps,s)=\eps^{\kappa}\cdot |s|^\iota$ with $\kappa \in (0,1)$, $\iota>0$ and recall that there is $C>0$ such that $\max\{|1- \cos (x)|, |\sin (x) - x|\} \leq C|x|^2$. The latter fact and  $(R_{\alpha^*}\One)' = \alpha^*( x^{\alpha^*-1}(1-x)^{\alpha^*}+ x^{\alpha^*}(1-x)^{\alpha^*-1})\,,$ imply that
\begin{equation}\label{eq:Strong 1-cos estim}
\sup_{D(\delta,\eps,x)} |(R_{\alpha^*}\One)'|\cdot\max\left\{1-\cos(s\chi)\,,\left|\sin (s\chi)-(s\chi)\right|\right\}\, \lesssim  \frac{|s|^2}{ (x-\eps)^{1+2\alpha - \alpha^*} } \,,
\end{equation}
when $x \leq 1/2$. The estimates for $x\geq 1/2$ follows from replacing $(x-\eps)$ by $(1-x+\eps)$, and the final estimates remain unchanged. So, we restrict our attention to the former case. 

This implies that the contribution of the first term in \eqref{eq: osc Ralpha cos} is 
\begin{align*}
  \MoveEqLeft\lim_{s\to0}\sup_{\eps\in (0,\eps_0]} \frac{1}{|s|\eps^\beta} {\int_{\delta+\eps}^{1/2} 2\eps{ \sup_{D(\delta,\eps,x)}  |(R_{\alpha^*}\One)'|\,(1-\cos (s\chi))}\One_{(\delta+\eps,1-\delta-\eps)}(x)\,\mathrm{d}\lambda_I(x)}\notag\\
  &\lesssim \lim_{s\to0} |s| \sup_{\eps\in (0,\eps_0]}\eps^{1-\beta}\int_{\delta}^{1/2-\eps} {x^{\alpha^*-1-2\alpha}}\,\mathrm{d}\lambda_I(x)\notag\\
 &\lesssim \begin{cases}
 \eps_0^{1-\beta+\kappa(\alpha^*-2\alpha)}\lim_{s\to0}|s|^{1+\iota\left(\alpha^*-2\alpha\right)}
 =0\,,& \alpha^* < 2\alpha \\
  \eps_0^{1-\beta}(|\log(1/2-\eps_0)|+ \kappa |\log(\eps_0)|) \lim_{s\to0} |s| + \iota\eps_0^{1-\beta}\lim_{s\to0}  |s|\, |\log|s||
 &  \alpha^* = 2\alpha \\
 \eps_0^{1-\beta} \lim_{s\to 0}|s| & \alpha^* >  2\alpha
 \end{cases}\\
 &=0\,
\end{align*}
provided that, in the $\alpha^*<2\alpha$ case, 
\begin{equation}\label{eq:iotakappa4}
\begin{aligned}
     1-\beta+\kappa(\alpha^*-2\alpha)>0 &\iff \kappa< (1-\beta)/(2\alpha-\alpha^*)\,\\
     1+\iota(\alpha^*-2\alpha)>0 &\iff \iota< 1/(2\alpha -\alpha^*)\,,
\end{aligned}\,
\end{equation}
and similarly,
\begin{align*}
  \MoveEqLeft\lim_{s\to0}\sup_{\eps\in (0,\eps_0]}\frac{1}{s\eps^\beta}\int  2\eps \sup_{D(\delta,\eps,x)} {|(R_{\alpha^*}\One)'|\,}(\sin (s\chi)-s\chi)_{\pm}\One_{(\delta+\eps,1-\delta-\eps)}(x)\,\mathrm{d}\lambda_I(x)
  = 0\,.
\end{align*}

Next, we estimate the second summand of \eqref{eq: osc Ralpha cos}. 
Using $(1-\cos(s\chi))' = \sin(s \chi) \cdot s \chi'\,,$ $|\sin (s\chi)| \leq |s\chi|$, and our assumption about $\chi$ and $\chi'$ we have 

\begin{align*}
 \sup_{D(\delta,\eps,x)} (R_{\alpha^*}\One)\,| (1-\cos(s\chi))'|
    & \lesssim 
    \begin{cases}
    |s|^2 (x-\eps)^{\alpha^*-(\alpha+b)} & \alpha^*<\alpha+b \\
    |s|^2 \cdot 1 & \alpha^*\geq \alpha+b\,.
    \end{cases}
\end{align*}
Also, note that the estimate for $x \leq 1/2$ and for $x \geq 1/2$ are the same with $(x-\eps)$ replaced by $(1-x+\eps)$. 
 Thus, when $\alpha^*<\alpha+b$
 \begin{align*}
  \MoveEqLeft\lim_{s\to0}\sup_{\eps\in (0,\eps_0]}\frac{\eps}{s\eps^\beta}\int_{\delta+\eps}^{1/2}  
  \sup_{D(\delta,\eps,x)} (R_{\alpha^*}\One)\,| (1-\cos(s\chi))'|\One_{(\delta+\eps,1-\delta-\eps)}(x)\,\mathrm{d}\lambda_I(x)
 \notag\\ &\lesssim \lim_{s\to0}\sup_{\eps\in (0,\eps_0]} \eps^{1-\beta}|s|\int_{\delta}^{1/2-\eps} {x^{\alpha^*-(\alpha+b)}} \mathrm d\lambda_I(x) \notag \\
 &\lesssim\begin{cases}
 \eps_0^{1-\beta+\kappa(1+\alpha^*-\alpha-b)}\lim_{s\to0} |s|^{1+\iota(1+\alpha^*-\alpha-b)} & \alpha+b > 1+\alpha^* \\
    \eps_0^{1-\beta}(|\log(1/2-\eps_0)|+ \kappa |\log(\eps_0)|) \lim_{s\to0} |s| + \iota\eps_0^{1-\beta}\lim_{s\to0}  |s|\, |\log|s|| &  \alpha+b = 1+\alpha^* \\
 \eps_0^{1-\beta}\lim_{s\to0} |s| & \alpha+b < 1+\alpha^*
 \end{cases}\\
  &= 0\,.\end{align*}
where, in the case of $\alpha + b > 1 + \alpha^*$, we have assumed that \begin{equation}\label{eq:iotakappa5}  \begin{aligned} 1-\beta+\kappa(1+\alpha^*-\alpha-b)>0 &\iff\kappa <(1-\beta)/(\alpha+b-1-\alpha^*)\,, \\ 1+\iota(1+\alpha^*-\alpha-b)>0 &\iff \iota < 1/(\alpha+b-1-\alpha^*)\,. \end{aligned}  \end{equation}
 Analogously, under the same assumptions on $\kappa$ and $\iota$, 
 we obtain
  \begin{align*}
  \MoveEqLeft\lim_{s\to0}\sup_{\eps\in (0,\eps_0]}\frac{\eps
 }{s\eps^\beta} \int {\sup_{D(\delta,\eps,x)}R_{\alpha^*}\One(x) \left|\left((\sin(s\chi)-s\chi)_{\pm}\right)'\right|} 
 \One_{(\delta+\eps,1-\delta-\eps)}(x)\,\mathrm{d}\lambda_I(x)
  = 0\,
 \end{align*}
because $|(\sin(s\chi)-s\chi)'|= |\cos(s \chi) - 1| \cdot |s \chi'|\,$ and $|\cos(s \chi) - 1| \leq |s\chi|$. 

Next, we look at the second summand of \eqref{eq: osc with indicator function deriv}. Using \eqref{eq:InterEpsBall}, 
our assumption about $\chi$ and the symmetry around $x=1/2$, the corresponding integral over the second summand is dominated by
 \begin{align*}
  \MoveEqLeft\lim_{s\to0}\sup_{\eps\in (0,\eps_0]}\frac{
  2}{s\eps^\beta} { \left( \int_{\delta}^{\delta+2\eps}+\int_{1-\delta-2\eps}^{1-\delta}\right)} \sup_{D(\delta,\eps,x)}R_{\alpha^*}(1-\cos(s\chi))
 \,\mathrm{d}\lambda_I(x)\notag\\
 &\lesssim \lim_{s\to0}\sup_{\eps\in (0,\eps_0]} |s|\eps^{-\beta} \int_{\delta}^{\delta+2\eps} \max\{ (\delta+\eps)^{\alpha^* - 2\alpha}, (\delta+3\eps)^{\alpha^* - 2\alpha}\} \,\mathrm{d}\lambda_I(x)\notag\\
 &\lesssim|s|^{1-\iota(2\alpha-\alpha^*)} \lim_{s\to0} \eps_0^{1-\beta}|s| \\
 &=0\,.
 \end{align*}
Here, in the case of $\alpha^*<2\alpha$, we have to assume additionally that
 \begin{equation}\label{eq:iotakappa6}
  \begin{aligned}
1-\beta-\kappa(2\alpha-\alpha^*)>0 &\iff\kappa <(1-\beta)/(2\alpha-\alpha^*)\,, \\
1-\iota(2\alpha-\alpha^*)>0 &\iff \iota < 1/(2\alpha-\alpha^*)\,.
\end{aligned}  
\end{equation}
 Analogously, under the same assumptions on $\kappa$ and $\iota$, using our assumption about $\chi$\,, 
 we have
 \begin{align*}
 \MoveEqLeft\lim_{s\to0}\sup_{\eps\in (0,\eps_0]}\frac{ 2
 }{s\eps^\beta} {\left(\int_{\delta}^{\delta+2\eps}+\int_{1-\delta-2\eps}^{1-\delta}\right)} \sup_{D(\delta,\eps,x)}R_{\alpha^*}(\sin(s\chi)-s\chi)_{ \pm}\,\mathrm{d}\lambda_I(x) =0.
 \end{align*}

Finally, we investigate \eqref{eq: 1-e alphabeta 1a}. The estimations for \eqref{eq: 1-e alphabeta 3a} then follow analogously.  We split the integral as in \eqref{eq: estim near zero0}. 

For the first integral, due to \Cref{lem:KeyEst1}, we have
 \begin{align}
 \MoveEqLeft\lim_{s\to 0}\sup_{\eps\leq \eps_0}\frac{1}{s\eps^{\beta}}\int_{{[0,\delta)}}\osc(R_{\alpha^*}(1-\cos(s\chi))\One_{[0,{\delta+\eps}]},B_{\eps}(x))\,\mathrm{d}\lambda_I(x)=0\label{eq: estim near zero1 deriv}
 \end{align}
provided that 
\begin{equation}\label{eq:iotakappa7}
  \begin{aligned}
  \kappa(1+\alpha^*)-\beta>0&\iff \kappa>\beta/(1+\alpha^*)\,,\\
\iota-1>0&\iff \iota>1\,.
\end{aligned}  
\end{equation}
For the imaginary part, since we assumed $\alpha^* \geq \alpha\,,$ we can use the following estimate. 
\begin{align*}
    {\sup}_{\bar D(\delta,\eps,x)}|R_{\alpha^*}(\sin(s\chi)-s\chi) \One_{[0,{\delta+\eps}]}| &\lesssim |s| {\sup}_{\bar D(\delta,\eps,x)}|R_{\alpha^*}\chi \One_{[0,{\delta+\eps}]}| \\ &\lesssim |s| \, {\sup}_{\bar D(\delta,\eps,x)} R_{\alpha^*-\alpha}\One_{[0,{\delta+\eps}]} \lesssim |s|\, (x+\eps)^{\alpha^*-\alpha}.
\end{align*}
Then repeating the argument leading to \Cref{eq:kappaiotaC2} with $\alpha^*-\alpha$ replacing $\alpha^*$, we have that 
\begin{align}
 \MoveEqLeft\lim_{s\to 0}\sup_{\eps\leq \eps_0}\frac{1}{s\eps^{\beta}}\int_{{[0,\delta)}}\osc(R_{\alpha^*}(\sin(s\chi)-s\chi)\One_{[0,{\delta+\eps}]},B_{\eps}(x))\,\mathrm{d}\lambda_I(x) = 0  \notag \label{eq: estim near zero1 deriv img}
 \end{align}
 provided that 
\begin{equation}\label{eq:iotakappa7a}
  \begin{aligned}
  \kappa(1+\alpha^*-\alpha)-\beta>0&\iff \kappa>\beta/(1+\alpha^*-\alpha)\,,\\
\iota-1>0&\iff \iota>1\,.
\end{aligned}  
\end{equation}

For the second integral, as in \eqref{eq:Near0SupEst} but using \eqref{eq:Strong 1-cos estim} instead, we obtain for all $x \in (\delta, \delta+2\eps]$
\begin{equation*}
    \sup_{\bar D(\delta,\eps,x)} R_{\alpha^*}(1-\cos (s\chi))\One_{[\delta,\delta+\eps]}\lesssim s^2\sup_{y \in \bar D(\delta,\eps,x) \cap (\delta, \delta+2\eps]} y^{\alpha^*-2\alpha} \leq \begin{cases}
 s^2 & \alpha^* \geq  2\alpha \\
s^2\cdot \delta^{\alpha^*-2\alpha} & \alpha^* < 2\alpha\,. \end{cases}
\end{equation*}
Therefore,
 \begin{align*}
 \MoveEqLeft\lim_{s\to 0}\sup_{\eps\leq \eps_0}\frac{1}{s\eps^{\beta}}\int_{{(\delta,\delta+2\eps]}}\osc(R_{\alpha^*}(1-\cos(s\chi))\One_{[\delta,\delta+\eps]},B_{\eps}(x))\,\mathrm{d}\lambda_I(x)\notag\\
  &\lesssim  \begin{cases}
  \eps_0^{1-\beta}\lim_{s\to 0}  |s|
  &\alpha^* \geq  2\alpha\\
  \eps_0^{1-\beta+\kappa(\alpha^*-2\alpha)}\lim_{s\to 0}|s|^{1+\iota(\alpha-2\alpha^*)}
  &\alpha^*< 2\alpha\end{cases}\notag\\
  &=0 
 \end{align*}
 provided that, in the case of $\alpha^*< 2\alpha ,$ 
\begin{equation}\label{eq:iotakappa8}
  \begin{aligned}
1-\beta+\kappa\left(\alpha^*-2\alpha\right)>0 &\iff\kappa < (1-\beta)/(2\alpha-\alpha^*)\,, \\
1+\iota(\alpha^*-2\alpha)>0&\iff \iota< 1/(2\alpha-\alpha^*)\,.
\end{aligned}  
\end{equation}
Due to \cite[Prop.~3.2(ii)]{BS} and our assumption that $\alpha^* > \alpha$ we have for $x\in (\delta, \delta+ 2\eps]$,
\begin{align*}
    &\osc\left(R_{\alpha^*}(1-\cos(s\chi))\One_{[0,\delta]},B_{\eps}(x)\right)\\ 
    &\qquad\leq \osc\left(R_{\alpha^*}(1-\cos(s\chi)),B_{\eps}(x)\cap [0,\delta]\right)\One_{[0,\delta]}(x)
    +2\esssup_{B_{\eps}(x)\cap [0,\delta]} R_{\alpha^*}(1-\cos(s\chi))\One_{ [\delta-\eps \vee 0, \delta+\eps]}(x)\\
    &\qquad \leq 0+2|s|\esssup_{[\delta-\eps \vee 0,\delta]} R_{\alpha^*}\chi
     \leq 2|s| \delta^{\alpha^*-\alpha}.
\end{align*}
So, 
\begin{align*}
    \MoveEqLeft
  \lim_{s\to 0}\sup_{\eps\leq \eps_0}\frac{1}{s\eps^{\beta}} \int_{({\delta,\delta+2\eps}]}\osc\left(R_{\alpha^*}(1-\cos(s\chi))\One_{[0,\delta]},B_{\eps}(x)\right)\,\mathrm{d}\lambda_I(x)\notag\\  
  &\lesssim \lim_{s\to 0}\sup_{\eps\leq \eps_0}\frac{1}{\eps^{\beta}} \int_{({\delta,\delta+2\eps}]}\delta^{\alpha^*-\alpha}\,\mathrm{d}\lambda_I(x)
  \lesssim \lim_{s\to 0} \eps_0^{1-\beta+\kappa\alpha^*}s^{\iota(\alpha^*-\alpha)}=0
\end{align*}
under \eqref{eq:iotakappa7}. 
So we have
\begin{align*}
  \lim_{s\to 0}\sup_{s\eps\leq \eps_0}\frac{1}{\eps^{\beta}} \int_{[{\delta,\delta+2\eps}]}\osc\left(R_{\alpha^*}(1-\cos(s\chi))\One_{[0,\delta+2\eps]},B_{\eps}(x)\right)\,\mathrm{d}\lambda_I(x)=0.
\end{align*}

Finally, we discuss here values of $\alpha^*$ and implicit restrictions on $b$ that ensure the existence of $\iota>0$ and $\kappa \in (0,1)$ used in the proof. There are four 
key cases to consider.
\begin{enumerate}[leftmargin=*, itemsep=5pt]
    \item\label{en: cond1} $\alpha < \alpha^* < 2\alpha$ and $\alpha+b>1+\alpha^*$: 
   \noindent Under \eqref{eq:iotakappa4}, 
\eqref{eq:iotakappa5},  \eqref{eq:iotakappa6}, \eqref{eq:iotakappa7},  \eqref{eq:iotakappa7a} and \eqref{eq:iotakappa8}, 
we have
\begin{equation*}
    \begin{aligned}
        \frac{\beta}{\alpha^*-\alpha+1} <\, &\kappa < \min\left\{\frac{1-\beta}{2\alpha-\alpha^*}, \frac{1-\beta}{\alpha+b-1-\alpha^*}\right\} \\
         1 <\, &\iota < \min\left\{\frac{1}{2\alpha-\alpha^*}, \frac{1}{\alpha+b-1-\alpha^*}\right\}.        
    \end{aligned}
\end{equation*}
Considering the conditions for $\iota$, we have 
$\alpha^*>2\alpha-1$ and $ \alpha^* > \alpha+b-2$. Since $\alpha>2\alpha-1$, the former is automatic. Next, considering each of the two inequalities that guarantee the existence of $\kappa$, we obtain that
$$\alpha^* > \max\left\{\beta-1 + \beta\alpha +\alpha,\,\beta b+\alpha-1 \right\} = \alpha + \max\left\{\beta-1 + \beta\alpha ,\,\beta b -1 \right\}$$ is necessary.
Note that $\beta-1 + \beta\alpha < 0 $ 
and $\beta b-1 < 0$ because $\beta<\min\{1/b,1/(1+\alpha)\}$.  So $\alpha^* > \alpha$ is a sufficient choice. 
Combining everything, we have that
$$\max\{\alpha+b - 2, \alpha\} < \alpha^* < \min\{\alpha +b -1, 2\alpha\}$$
is sufficient. 

\item\label{en: cond2} $\alpha <  \alpha^* < 2\alpha$ and $\alpha+b <  1+\alpha^*$: 
\noindent\eqref{eq:iotakappa5} poses no extra restriction.  So, under \eqref{eq:iotakappa4}, \eqref{eq:iotakappa6}, \eqref{eq:iotakappa7}, \eqref{eq:iotakappa7a} and \eqref{eq:iotakappa8}, we have $\alpha^*  > \alpha$ as before. Hence, $$\max\{\alpha+b-1,\alpha\} < \alpha^*< 2\alpha$$ 
is sufficient. 

\item\label{en: cond3} $\alpha^*  > 2\alpha$ and $\alpha+b > 1+\alpha^*$:
\eqref{eq:iotakappa4}, \eqref{eq:iotakappa6}  and \eqref{eq:iotakappa8} pose no extra restrictions.
 Under 
\eqref{eq:iotakappa5},  \eqref{eq:iotakappa7} and \eqref{eq:iotakappa7a} we have 
 $\alpha^* < \alpha+b-1$ and $\alpha^* >  \beta b+\alpha-1$. Since $\beta b+\alpha-1 < \alpha < 2\alpha $, the latter is true. So,
$$ 2\alpha< \alpha^*  <  \alpha+b-1$$
is sufficient. 

\item\label{en: cond4} $\alpha^*  > 2\alpha$ and $\alpha+b <   1+\alpha^*$:
\eqref{eq:iotakappa5} is not relevant, and both  \eqref{eq:iotakappa7} and \eqref{eq:iotakappa7a} pose no extra restrictions. Hence, 
$$\max\{2\alpha,\alpha+b-1\} < \alpha^*$$
is sufficient. 
\end{enumerate}
We obtain from \eqref{en: cond1} and \eqref{en: cond2} that
$\max\{\alpha+b-2,\alpha\}<\alpha^*$ is sufficient if $\alpha^* <2\alpha$. From \eqref{en: cond3} and \eqref{en: cond4} we obtain that $\alpha^* > 2\alpha$ is sufficient if $\alpha^* > 2\alpha$. So, $$\alpha^* > \min\{2\alpha, \max\{\alpha, \alpha+b-2\}\}$$ is sufficient. 
\end{proof}

\begin{lem}\label{lem:abtoalphabeta}
Assume $\chi$ is continuous, the right and left derivatives of $\chi$ exist  on $\mathring{I}$, and there exist $a\geq 0\,,b>0$ such that
 \begin{align}\label{eq:bddness condition}
 |\chi(x)|\lesssim x^{-a}(1-x)^{-a}
 \quad\text{ and }\quad
  \max\{|\chi'(x+)|, |\chi'(x-)|\}\lesssim x^{-b}(1-x)^{-b}\,,
 \end{align}
 then
 $\|\chi\|_{\alpha,\beta,\gamma}<\infty$
 if
 \begin{align}\label{eq: cond alpha beta a b}
 \alpha &> a\,,\nonumber\\ 
 \beta < (1+ \alpha -a)/(b-a)\quad &\text{or}\quad b<a+1\, \quad\text{and}\\ 
 1 \leq \gamma &< 1/a \,.\nonumber
 \end{align}
\end{lem}

\begin{proof} 
 The first inequality of \eqref{eq:bddness condition} 
 implies $\chi \in L^\gamma$ with $1 \leq \gamma<1/a$.  
 
 For simplicity we assume $\chi$ is differentiable. Otherwise, at a point where $\chi$ is not differentiable, both one-sided derivatives will exist and the following estimates do hold for them.

Now, we proceed as in the proof of \Cref{rem:RegObs3}, however, with $\delta = \eps^\kappa$ to find the minimal $\alpha$ and maximal $\beta$ such that $|R_\alpha \chi|_{0,\beta}<\infty$. Set $g\coloneqq R_\alpha\chi$, then
 \begin{align*}
  g'(x)= \alpha(1-2x)R_{\alpha-1} \chi(x)+  R_\alpha\chi'(x).
 \end{align*}
Choose $\eps$ sufficiently small and split the domain into three parts, $[0,\eps^\kappa+\eps), (\eps^\kappa+\eps,1-\eps-\eps^\kappa)$ and $(1-\eps-\eps^\kappa,1]$. Due to the symmetry of the bounds, we only focus on $[0,1/2]$. 

On $(\eps^\kappa+\eps,1-\eps-\eps^\kappa)\,,$ we use \cite[Prop.~3.2(ii)]{BS} implying
\begin{align}
  &\osc\left(g\,\One_{(\eps^\kappa+\eps,1-\eps^\kappa-\eps)},B_{\eps}(x)\right)\nonumber\\ &\qquad\leq \osc(g, D(\eps^\kappa,\eps,x))\One_{(\eps^\kappa+\eps,1-\eps^\kappa-\eps)}(x)\,+\,2 \Big(\sup_{D(\eps^\kappa,\eps,x)} g \Big)\left(\One_{B_{\eps}(\eps^\kappa+\eps) \cup B_{\eps}(1- \eps^\kappa-\eps) }(x)\right)
  \label{eq: osc last estim}
 \end{align}
 with $D$ as in \eqref{eq: def D}.

For the following we set $\tilde\alpha = \min\{\alpha-b+1, \alpha-a\}$.
Then the contribution from the first term to $|R_\alpha\chi|_{0,\beta}$ is (up to a constant) bounded by
\begin{align*}
    \sup_{0 < \eps \leq \eps_0}\eps^{1-\beta}\int_{\eps^\kappa+\eps}^{1/2}\sup_{D(\eps^\kappa,\eps,x)} g'\,\mathrm{d}\lambda_I(x) &\lesssim \sup_{0 < \eps \leq \eps_0}\eps^{1-\beta}\int_{\eps^\kappa+\eps}^{1/2}(x-\eps)^{-a+\alpha-1} +(x-\eps)^{-b+\alpha}\,\mathrm{d}\lambda_I(x)\\
    &\lesssim \sup_{0 < \eps \leq \eps_0}\eps^{1-\beta}\int_{\eps^\kappa}^{1/2}x^{\bar\alpha-1}\,\mathrm{d}\lambda_I(x)\\\
    &\lesssim \begin{cases}
    \sup_{0 < \eps \leq \eps_0}\eps^{1-\beta+\kappa\bar \alpha} &  \tilde\alpha <  1\\
    \sup_{0 < \eps \leq \eps_0}\eps^{1-\beta}(\log(1/2)-\kappa\log(\eps)) & \bar \alpha=1\\
    \eps_0^{1-\beta} & \tilde\alpha >  1\,.
    \end{cases}
\end{align*}
In the $\tilde\alpha \leq 1$ case, we require that 
 \begin{equation}\label{eq: cond for alpha beta a b} 
1-\beta+\kappa\tilde \alpha >0 \iff \big(\kappa < (1-\beta)/(b-\alpha-1)\,\text{ or }b<\alpha+1\big)\,,
\end{equation}
where we have made use of the fact $\alpha > a$.
On the other hand, \eqref{eq: cond for alpha beta a b} is automatically fulfilled if $\tilde\alpha > 1$, so we do not  have to distinguish the cases anymore.

Since $\alpha > a\,$ the contribution from the second term in \eqref{eq: osc last estim} is bounded by 
\begin{align*}
 &\sup_{\eps\in (0,\eps_0]}\frac{1}{\eps^\beta} \int_{\eps^\kappa}^{\eps^{\kappa}+2\eps} \sup_{D(\eps^\kappa,\eps,x)}g
 \,\mathrm{d}\lambda_I(x)\lesssim\sup_{\eps\in (0,\eps_0]} \eps^{-\beta} \int_{\eps^{\kappa}}^{\eps^{\kappa}+2\eps} 
 1\,\, \,\mathrm{d}\lambda_I(x)\lesssim \eps_0^{1-\beta}\,. 
\end{align*}

Now, for $x \in [0,\eps^\kappa)$ we use the following estimate
\begin{align*}
    {\sup}_{\bar D(\eps^\kappa,\eps,x)}|g| &\lesssim  {\sup}_{\bar D(\eps^\kappa,\eps,x)}|R_{\alpha-a}\One_{[0,{\eps^\kappa+\eps}]}| \lesssim  (x+\eps)^{\alpha- a}
\end{align*}
 with $\bar D$ as in \eqref{eq: bar D}.
Following the argument in \Cref{lem:KeyEst1} with $\alpha - a$ replacing $\alpha^*$ and \textit{without} the $s \to 0$ limit but fixing $s=1$, we have, since $\alpha-a+1>\beta$ automatically holds, that 
\begin{align*}
 \sup_{\eps\leq \eps_0}\frac{1}{\eps^{\beta}}\int_{{[0,\delta)}}\osc(g,B_{\eps}(x))\,\mathrm{d}\lambda_I(x)   &\lesssim \sup_{\eps\leq \eps_0}\frac{2 \left(\left(\eps^\kappa+\eps\right)^{{\alpha - a}+1}- \eps^{{\alpha - a}+1}\right)}{({\alpha - a}+1)\eps^{\beta}} \\ &= \frac{2 \left(\left(\eps_0^\kappa+\eps_0\right)^{{\alpha - a}+1}- \eps_0^{\alpha - a+1}\right)}{(\alpha - a+1)\eps_0^{\beta}}
 \end{align*}
 provided that 
\begin{equation*}
\kappa(1+{\alpha - a})-\beta>0 \iff \kappa>\frac{\beta}{1+\alpha - a}\,.
\end{equation*}
 
So, together with \eqref{eq: cond for alpha beta a b} we require that there exists $\kappa$ such that
$$\frac{\beta}{1+{\alpha - a}} < \kappa < \frac{1-\beta}{b-\alpha-1}\quad\text{or}\quad b<\alpha+1\,.$$
This is true if and only if 
$$(b-a)\beta < 1+ \alpha -a \quad\text{or}\quad b<\alpha+1\,. $$
\end{proof}

\section{H{\"o}lder Continuity of \texorpdfstring{$\bar{R}_{j+1}$}{R j+1}}\label{sec:alphaHolder}

\begin{lem}\label{lem: Rj+1 Hoelder}
 For all $j=0,
 \dots, k-1$, let $\bar R_{j+1}:[c_j,c_{j+1}] \to \mathbb{R}$ be given by
 $$ \bar{R}_{j+1}=\frac{(R_\alpha\One) \circ  \psi_{j+1}}{R_\alpha\One} \,.$$ 
 Then $\bar{ R}_{j+1}$ is bounded and $\alpha$-H\"older continuous for all $j$.
\end{lem}
\begin{proof}

 Our strategy is to prove the following two steps:
 \begin{enumerate}[leftmargin=*]
  \item\label{en: proof strategy 1}  There exists $\delta_0>0$ such that $\bar R_1'$ is bounded on the interval $[0,c_1-\delta_0)$, $\bar R_{k+1}'$ is bounded on the interval $(c_{k}+\delta_0,1]$ and $\bar R_{j+1}'$, $j=1,\ldots, k-1$ is bounded on the interval $(c_{j}+\delta_0,c_{j+1}-\delta_0)$.
  \item 
  Since $ \bar R_{j+1}(c_j)=\bar R_{j+1}(c_{j+1})=0$ for $j=1,\ldots,k-1$, it is enough to show that there exists $C>0$ such that  $\bar R_{j+1}(c_{j}+\eps)\leq C\eps^{\alpha}$ and $\bar R_{j+1}(c_{j+1}-\eps)\leq C\eps^{\alpha}$, for all $\eps>0$. 
 \end{enumerate}
We have 
 \begin{align}\label{eq: bar R deriv}
   \bar R_{j+1}'(x)
  =\alpha \cdot\frac{\psi_{j+1}'(x)(1-2\psi_{j+1}(x))x(1-x)-\psi_{j+1}(x)(1-\psi_{j+1}(x))(1-2x)}{(\psi_{j+1}(x)(1-\psi_{j+1}(x)))^{1-\alpha}\left(x\left(1-x\right)\right)^{1+\alpha}}.
 \end{align}
  The numerator is bounded, and for $j=1,\ldots, k-2$, the denominator has zeros only at $c_{j}$ and $c_{j+1}$. So, we immediately get that $ \bar R_{j+1}'$ is bounded on $(c_{j}+\delta_0,c_{j+1}-\delta_0)$ . 
 
 We only have to further consider the cases $j=0$ and $j=k-1$. We have to show that $\bar R_1'(x)$ is bounded in a neighbourhood of $0$.  
 Since $\psi_1$ has a bounded second derivative, we can write $\psi_1(x)=\psi_1(0)+\psi'_1(0)x+\cO(x^2)=\psi'_1(0)x+\cO(x^2)$.
 This yields $(\psi_1(x)(1-\psi_1(x)))^{1-\alpha}\left(x\left(1-x\right)\right)^{1+\alpha}=\Omega(x^2)\,.$\footnote{$f(x)=\Omega(g(x))$ as $x\to 0$ if $\liminf_{x\to0}|f(x)|/g(x)>0\,.$} 
 On the other hand, by simply multiplying out we obtain 
 \begin{align*}
  \psi_1'(x)(1-2\psi_1(x))x(1-x)-\psi_1(x)(1-\psi_1(x))(1-2x)
  &=\cO(x^2)
 \end{align*}
 implying that $\lim_{x\to 0}\bar R_1(x)<\infty$.
 The calculation for $\lim_{x\to 1} \bar R_k(x)$ follows analogously.

 In order to analyse the behaviour for $x\to c_{j}$ and $x\to c_{j+1}$  with $x$ starting from $[c_j,c_{j+1}]$ we note that $\bar R_{j+1}'$ can be written as 
 \begin{align}
  \bar R_{j+1}'(x)
  &=\alpha \cdot\frac{\psi_{j+1}'(x)(1-2\psi_{j+1}(x))}{(\psi_{j+1}(x)(1-\psi_{j+1}(x)))^{1-\alpha}\left(x\left(1-x\right)\right)^{\alpha}}
  -\alpha \cdot\frac{\left(\psi_{j+1}(x)(1-\psi_{j+1}(x))\right)^{\alpha}(1-2x)}{\left(x\left(1-x\right)\right)^{1+\alpha}}.\label{eq: bar Rj+1'}
 \end{align} 
 The minuend tends to $\infty$ for $x\to c_{j}$ and {to $-\infty$} for $x\to c_{j+1}$ since $\psi_{j+1}(x)$ and $1-\psi_{j+1}(x)$ tend to zero, respectively, and the numerator remains bounded and is positive near $c_j$ and negative near $c_{j+1}$. The subtrahend is bounded on an interval $[\delta_0, 1-\delta_0]$. Thus, $\bar R_{j+1}'(x)$ tends to $\infty$ for $x\to c_{j}$ and {to $-\infty$ for} $x\to c_{j+1}$ except if $c_{j}=0$ or $c_{j+1}=1$. 
 
 Hence, we can conclude that $|\bar R_{j+1}(x)-\bar R_{j+1}(y)|\leq \bar R_{j+1}(c_{j}+|x-y|)-\bar R_{j+1}(c_{j})=\bar R_{j+1}(c_{j}+|x-y|)$ for $x,y\in [c_{j}, c_{j}+\delta_0]$ and $\delta_0>0 $ sufficiently small. 
 Similarly, we have $|\bar R_{j+1}(x)-\bar R_{j+1}(y)|\leq  \bar R_{j+1}(c_{j+1}-|x-y|)-\bar R_{j+1}(c_{j+1})=\bar R_{j+1}(c_{j+1}-|x-y|)$ for $x,y\in [c_{j+1}-\delta_0,c_{j+1}]$ and $\delta_0>0$ sufficiently small. 
 On the other hand, we have   
 \begin{align*}
  \bar R_{j+1}(c_j-\eps)
  &=\left(\frac{\psi_{j+1}(c_j-\eps)(1-\psi_{j+1}(c_j-\eps))}{(c_j-\eps)(1-c_j+\eps)}\right)^{\alpha}.
 \end{align*}
There exists $C_{j,\delta_0}>0$ such that
 \begin{align*}
  \left(\frac{\psi_{j+1}(c_j-\eps)}{(c_j-\eps)(1-c_j+\eps)}\right)^{\alpha}<C_{j,\delta_0}
 \end{align*}
uniformly for all  $\eps\in (0,\delta_0)$ and thus
\begin{align*}
   \bar R_{j+1}(c_j-\eps)&\leq C_{j,\delta_0}\left(\eta_+\eps\right)^{\alpha}.
\end{align*}
Similarly, we have 
 \begin{align*}
  \bar R_{j+1}(c_{j-1}+\eps)
  &=\left(\frac{\psi_{j+1}(c_{j-1}+\eps)(1-\psi_{j+1}(c_{j-1}+\eps))}{(c_{j-1}+\eps)(1-c_{j-1}-\eps)}\right)^{\alpha}
 \end{align*}
and there exists $\bar C_{j,\delta_0}>0$ such that
 \begin{align*}
  \left(\frac{1-\psi_{j+1}(c_{j-1}+\eps)}{(c_{j-1}+\eps)(1-c_{j-1}+\eps)}\right)^{\alpha}<\bar C_{j,\delta_0}
 \end{align*}
uniformly for all  $\eps\in (0,\delta_0)$ and thus
\begin{align*}
   \bar R_{j+1}(c_j-\eps)&\leq \bar C_{j,\delta_0}\left(\eta_+\eps\right)^{\alpha}.
\end{align*}
Setting $C=\max_j\max\{C_{\delta_0,j}, \bar{C}_{\delta_0,j}\}$ concludes the proof of the lemma.
\end{proof}

\begin{ack}T.S.~was supported by the Austrian Science Fund FWF:~P33943-N and by the MSCA Project	ErgodicHyperbolic - p.n.\ 101151185. Furthermore, she acknowledges the support of the Università di Pisa through the ``visiting fellows'' programme, as this work was partially done during her visit to the Dipartimento di Matematica, Università di Pisa. T.S.~would like to thank the Centro De Giorgi, Scuola Normale Superiore di Pisa for their hospitality during various visits. K.F.~was supported by the UniCredit Bank R\&D group through the `Dynamics and Information Theory Institute’ at the Scuola Normale Superiore di Pisa. K.F.~would like to thank the Centro De Giorgi for the excellent working conditions and research travel support, the Fields Institute and the Toronto Public Library System for letting him use their common spaces while visiting Toronto and the Universität Wien for the hospitality during research visits to meet T.S.. Both authors would like to thank Carlangelo Liverani for useful discussions.
\end{ack}

\typeout{}

\end{document}